%% file: higharX2.tex
\documentclass[12pt]{amsart}
\usepackage{amssymb,amscd}
\usepackage{verbatim}

\usepackage[dvipdfmx]{graphicx,color}
\usepackage{graphicx}

\newtheorem*{Whitney towers}{Theorem~\ref{Whitney towers}}
\newtheorem*{h-towers}{Theorems ~\ref{half} \& \ref{$(n)$-solvable}}

\newtheorem*{surgery curves}{Theorem~\ref{surgery curves}}
\newtheorem*{cg=0}{Theorem~\ref{vanish}}

\newtheorem{thm}{Theorem}[section]

\theoremstyle{definition}
\newtheorem{defn}[thm]{Definition}

\newtheorem{que}[thm]{Question}

\numberwithin{equation}{section}
\numberwithin{figure}{section}

\newcommand{\x}{\times}
\newcommand{\np}{\newpage}
\newcommand{\Z}{\mathbb{Z}}
\newcommand{\N}{\mathbb{N}}

\newcommand{\Q}{\mathbb{Q}}

\newcommand{\R}{\mathbb{R}}

\usepackage{here}   %zu wooku  \begin{figure}[H]

\def\yen{{\setbox0=\hbox{Y}Y\kern-.97\wd0\vbox{hrule height.lex width.98%
\wd0\kern.33ex\hrule height.lex width.98\wd0\kern.45ex}}}%%

\def\np{\newpage}

\begin{document}
\pagestyle{plain}
%{\color{red}{  

\title{Introduction to high dimensional knots}
\author{Eiji Ogasa\\}

%\thanks{pqr100pqr100@yahoo.co.jp\quad ogasa@mail1.meijigakuin.ac.jp}
%\thanks{Keywords:    
% {\bf PACS nos.} 11-25w, 11-25Uv.
%\newline MSC2000 57N10, 57N13, 57N15. }
\date{}

\begin{abstract} 
This is an introductory article on high dimensional knots for the beginners.  
High dimensional knot theory is an exciting field. 
It is a field of knot theory, which is one of topology and is connected with many ones.

In this article we use  few literal expressions, equations, functions, etc.  
We barely suppose that the readers have studied manifolds, homology theory, or topics beyond them.  

Is there a truly knotted high dimensional knot? 
We first answer this question.  
After that, we explain local moves on high dimensional knots and 
the projections of high dimensional knots.  
\end{abstract}   
\maketitle 

\setcounter{tocdepth}{3}
\tableofcontents

\section{Introduction}\label{introduction} 
This is an introductory article on high dimensional knots for the beginners.  
High dimensional knot theory is an exciting field. 
It is a field of knot theory, which is one of topology and is connected with many ones.

In this article we use  few literal expressions, equations, functions, etc.  
We barely suppose that the readers have studied manifolds, homology theory, or topics beyond them.  
We only suppose that the readers know 
the $n$-dimensional spheres,  
the $n$-dimensional balls (or discs),  
and 
the $n$-dimensional Euclidean space $\R^n$. 
Here, $n$ is any natural number.

\bigbreak
A 1-dimensional knot is a circle in $\R^3$ which does not touch itself. 
Do you feel that  
there is an unknotted knot and a truly knotted knot? 
Yes, it is true. See \S\ref{ichiichi}. 

\begin{figure}[H]
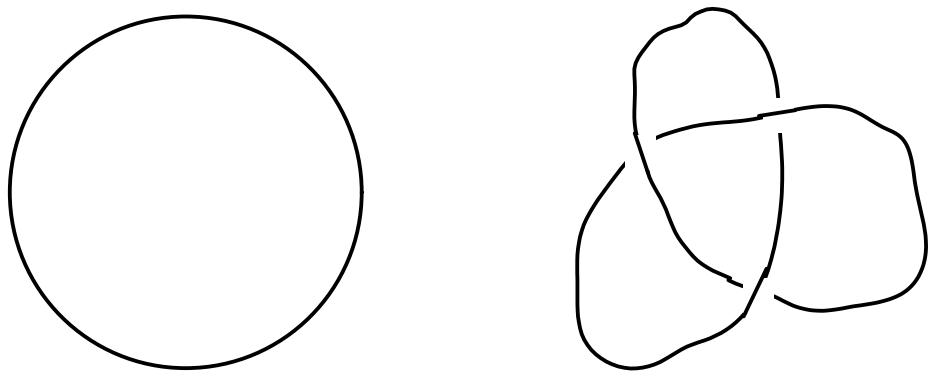

\bigbreak
\input 0Y207.tex  
\bigbreak 
\bigbreak\caption{{\bf 
There is an unknotted knot and a truly knotted knot.}\label{knotted}}   
\end{figure}

\bigbreak
A 2-dimensional knot is a sphere in $\R^4$ which does not touch itself.
Is there a truly knotted 2-dimensional knot? 
We answer this question in \S\ref{ichini}. 

\bigbreak
An $n$-dimensional knot is an $n$-sphere in $\R^{n+2}$ which does not touch itself ($\geqq3$). 
Is there a truly knotted $n$-dimensional knot? 
We answer this question in \S\ref{ichisan}.

\begin{figure}[H]
\bigbreak
\hskip15mm     \includegraphics[width=10cm]{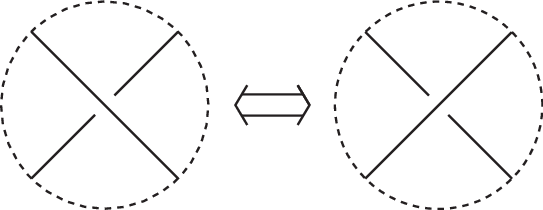}   
\bigbreak\caption{{\bf The crossing change of a 1-knot}\label{Alabama}}   
\bigbreak 
\end{figure}

\np
Any nontrivial 1-dimensional knot is changed into the unknot by a sequence of 
crossing-changes drawn in Figure \ref{Alabama}. 

%H 記述した部分  usepacckage使ってるとき。こっちが強力
%h 記述した部分 
%t ページの上部 
%b ページの下部 
%p 独立したページ 

\begin{figure}[H]
\bigbreak
\hskip15mm  \includegraphics[width=10cm]{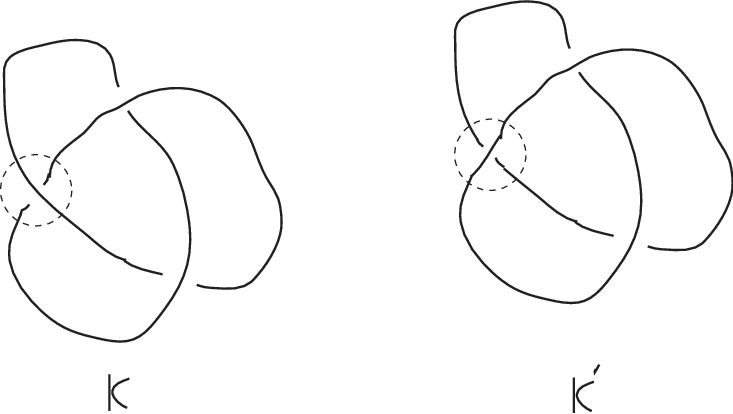}   
\bigbreak
\caption{{\bf $K$ is obtained from $K'$ by one crossing-change.}\label{AM}}   
\bigbreak 
\end{figure}

If two 1-dimensional knots $K$ and $K'$ in $\R^3$ differ only in the 3-ball $B^3$ as shown 
in Figure \ref{Alabama}, 
then we say that $K$ is obtained from $K'$ by one crossing-change. See Figure \ref{AM}.

Are there local moves on high dimensional knots 
like crossing-changes of 1-dimensional knots?
We answer this question in \S\ref{localmoves}.

\np 
The projection of any nontrivial 1-dimensional knot is 
that of an unknotted knot. 
The projection is the image of the map $\R^3\to \R^2,$\quad $(x,y,z)\mapsto(x,y)$. 
Probably you will feel it is true.  Yes. It is true. 
%We believe that you feel it is true.  
Indeed, the proof is easy. 

\begin{figure}[H]
\bigbreak
\includegraphics[width=12cm]{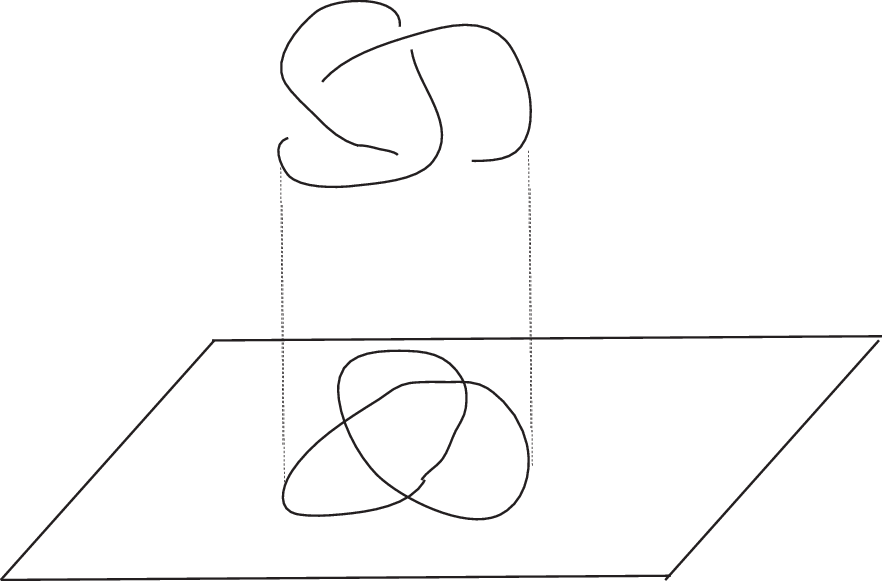}   
\bigbreak\caption{{\bf 
The projection of a 1-dimensional knot}
\label{0I1232pq}}   
\bigbreak 
\end{figure}

If there is a nontrivial high dimensional knot, 
then is 
the projection of it that of an unknotted knot? 
We answer this question in \S\ref{proj}.

\vskip1cm

We will show high dimensional figures to explain high dimensional knots and links.  
We draw them conceptually. 
We want the readers to captures them intuitively.  

{\Large Imagine!}

\np

We use  few literal expressions, equations, functions, etc.  
Of course the reason why we barely use them here 
is because this article is written for the beginners.  
However, note the following 
when you read papers and textbooks on mathematics and physics 
after finishing this article. 
Whether 
 the level of the topics in mathematics and physics is  high or low 
is not related to 
how complex are 
 the numbers and equations are used there.

\vskip1cm

We try to avoid using technical terms associated with manifolds, homology theory, or topics beyond them. 
Instead, the advanced readers may feel the explanation in this article to be quite loose. 
However, we put a high priority on explanations that beginners can understand  intuitively. 
(Of course, the advanced readers understand high dimensional figures  not only rigorously but also intuitively.)

\np
\section{High dimensional knots and links exist}\label{exist} 

\subsection{1-dimensional knots and links exist}\label{ichiichi}
We begin by explaining 1-dimensional knots. 

\begin{figure}[H]
\bigbreak
\includegraphics[width=9cm]{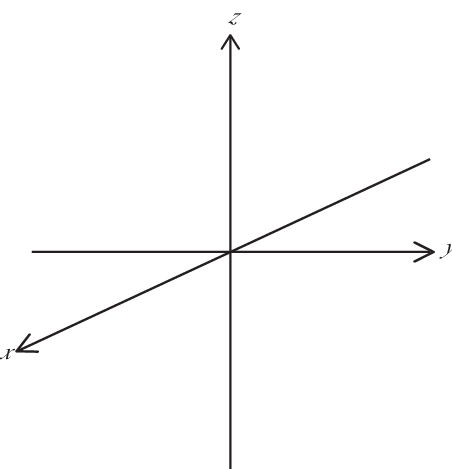}  
\bigbreak\caption{{\bf 
$\R^3$.}\label{xyz}}   
\bigbreak 
\end{figure}

Take the 3-dimensional space $\R^3$  as shown in Figure \ref{xyz}.
 
\noindent 
If a single copy of $S^1$ is embedded in $\R^3$, it is called 
a {\it 1-$($dimensional$)$ knot}.  Sometimes a knot means a 1-knot. 
When we say that $S^1$ is embedded in $\R^3$, 
$S^1$ is included in $\R^3$ and $S^1$ does not touch itself.  
%Embedding is a mathematical term. See textbooks on manifolds. 
%(Althogh, in fact, `embed' is a mathematical term, you may use the word as you use this word in daily life. \\
(In fact, `{\it embed}' is a mathematical term,  
which is seen textbooks on manifolds for embedding. 
However, if you currently associate the word 
`embed'  %`smooth' dashita ato naoshita
with 
how it is used in daily life, please continue to do so.  
In this article, we often hope that the readers take such attitude. 
This way is often good for the beginners.)  
\\

If a set of $m$ copies of $S^1$ is embedded in $\R^3$, 
it is called an 
{\it $m$-component 1$($-dimensional$)$ link}. 
Sometimes a link means a 1-link.

\noindent
{\bf To the advanced readers.}  
We replace $\R^3$ with the 3-dimensional sphere $S^3$ 
in a lot of literature. 
%We sometimes replace $\R^3$ with other manifolds than $S^3$.  

\np

Do you feel the two knots 
shown in Figure \ref{trefoil}    
seem to be different?

\begin{figure}[H]
\bigbreak
\input Y207.tex  
\bigbreak 
\bigbreak 
\bigbreak\caption{{\bf The trivial knot and the trefoil knot}\label{trefoil}}   
\bigbreak 
\end{figure}

%\np

Yes. It is true. 
Roughly speaking, when we say that they are not the same, it means the following. 
In $\R^3$ 
it is impossible to manipulate the trefoil knot to the trivial knot 
if it is never touching itself.

\begin{figure}[H]
\bigbreak
\includegraphics[width=18cm]{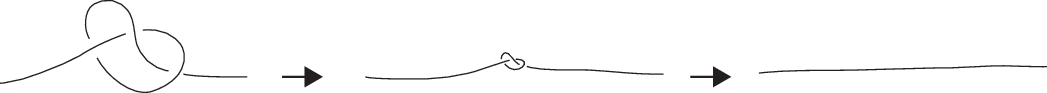}    
\bigbreak\caption{{\bf 
A funny move}\label{funny}}   
\bigbreak 
\end{figure}

\np
\noindent{\bf To advanced readers.}  
%we cannot move the trefoil knot to the trivial knot 
%such that it is always not touching itself, 
%%so that it does not touch itself,  
%%without touching itself, %self-touching, 
We show a funny move. 
See Figure \ref{funny}. 
The curve is not touching itself 
in the procedure 
where the left  figure is changed into the right one as shown in Figure \ref{funny}.

\bigbreak
\bigbreak

\begin{figure}[H]
\bigbreak
\includegraphics[width=8cm]{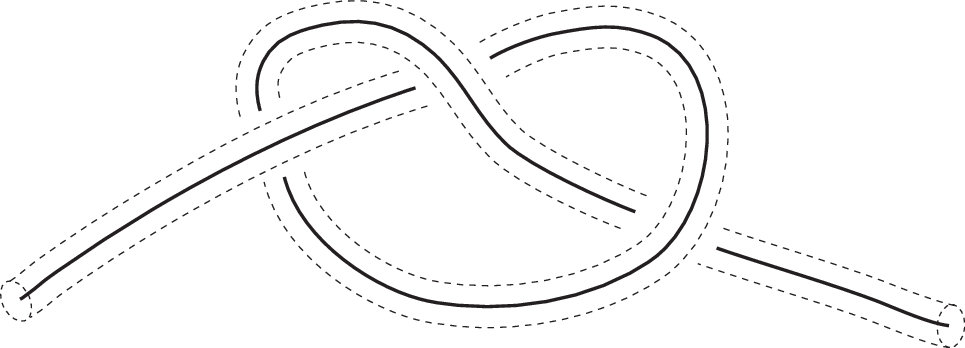}  
\bigbreak\caption{{\bf a thickened knot}\label{tubel}}   
\bigbreak 
\end{figure}

You feel that it is natural to prohibit this procedure, don't you?  
%Indeed the knot does not touch itself in this procedure.
%Indeed the knot is always not touching itself in the above procedure. 
%Indeed this procedure does not make self-touching.    
%Note that the knot is not thickened. 
So, when we say that we move a knot such that it is always not touching itself, 
%so that it does not touch itself, 
%without touching itself, %self-touching, 
we suppose that we move a knot with `tubular neighborhood' drawn in Figure \ref{tubel}
together 
and that the tubular neighborhood is not always touching itself. 
%does not touch itself.   %have self-touching. 
We can say that knots are `thickened' when we consider the procedure that knots are moved. 

\bigbreak
\bigbreak

In this article, if we move a knot as above,  
we say that we 
move it smoothly. 
(`Smooth' and `tubular neighborhood' are mathematical terms. 
%We explain it a little in the last Note of \S\ref{OH}. 
However, if you currently associate the word     
`smooth' and `tubular neighborhood' 
   %`smooth'   dashita ato naoshita
with 
how it is used in daily life, please continue to do so.)  
%Beginners may think this `smooth' is a word `smooth' in the daily life and read ahead.)

\np

If a 1-knot $K$ bounds an embedded 2-ball $B^2$ in $\R^3$, 
we say that $K$ is the {\it trivial knot} or the {\it unknot}. 

Suppose that an $m$-component  1-link $L$ satisfies the following condition:       
Each component bounds an embedded 2-ball $B^2$ in $\R^3$.  
Each $B^2$ does not touch any other $B^2$. 
Then we say that $L$ is the {\it trivial link}. 

%If  an $m$-component  1-link $L$ satisfies that        
%each component bounds an embedded 2-ball $B^2$ in $\R^3$ and that  
%each $B^2$ does not touch each other,  we say that $L$ is the trivial link. 

\vskip1cm
It is known that there are countably infinitely many different nontrivial 1-knots. 

In order to prove it, we use the Alexander polynomial, the Alexander module, 
the Jones polynomial, the fundamental groups, 
the signature, 
%the Arf invariant, 
the knot cobordism groups, etc. 
See 
\cite{Kauffman1994, Rolfsen} 
for detail.

\np

{\Large Can $S^2$ be `knotted' in $\R^4$?} 

\vskip1cm

 {\Large Can $S^n$ be `knotted' in $\R^{n+2}$?}  

\vskip1cm

Are there nontrivial high dimensional knots? 

\vskip1cm

These are among the main themes of this section. %Part III.  

\np

\subsection{2-dimensional knots and links exist}
\label{ichini}
Next we explain 2-dimensional knots. 

If a single copy of $S^2$ is embedded in $\R^4$, it is called 
a {\it 2-$($dimensional$)$ knot}. 

If a set of $m$ copies of $S^2$ is embedded in $\R^4$, 
it is called an {\it $m$-component 2-$($dimensional$)$ link}.

If a 2-knot $K$ bounds an embedded 3-ball $B^3$ in $\R^4$, 
we say that $K$ is the {\it trivial knot}. 

Suppose that  an $m$-component  2-link $L$ satisfies the following condition: 
Each component bounds an embedded 3-ball $B^3$ in $\R^4$.  
Each $B^3$ does not touch any other $B^3$. 
Then we say that $L$ is the {\it trivial link}. 

%If  an $m$-component  2-link $L$ satisfies that  
%each component bounds an embedded 3-ball $B^3$ in $\R^4$ and that 
%each $B^3$ does not touch each other, 
% we say that $L$ is the trivial link. 

We say that two 2-links are the same if we move one to the other smoothly 
in $\R^4$ so that it is not touching itself. 
%without touching itself, %self-touching, 
Here `moving smoothly' means that  
we move it with `tubular neighborhood' as in the 1-link case.  
Recall Figures \ref{funny} and  \ref{tubel}.

\vskip1cm
It is trivial that there is the 2-dimensional trivial knot (resp. 2-link).

\vskip1cm
There exist nontrivial 2-knots. We will explain it from here on.

\begin{figure}[H]
\bigbreak
\includegraphics[width=11cm]{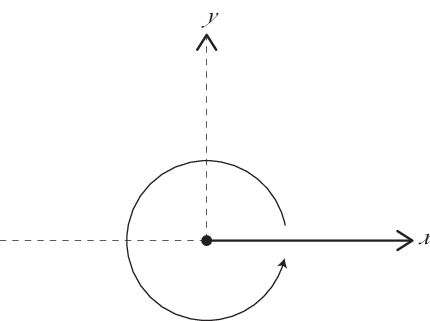}       
\bigbreak 
\bigbreak 
\bigbreak\caption{{\bf 
$\R^2$ is the result of rotating 
$\{(x,y) \vert x\geqq0, \quad y=0  \}$ 
around the point $(0,0)$}\label{roto1}}   
\bigbreak 
\end{figure}

\np
We regard $\R^2$ as the result of rotating $\{(x,y) \vert x\geqq0, \quad y=0  \}$ 
around the point $(0,0)$. 
See Figure \ref{roto1}.

\np

\begin{figure}[H]
\bigbreak
\includegraphics[width=9cm]{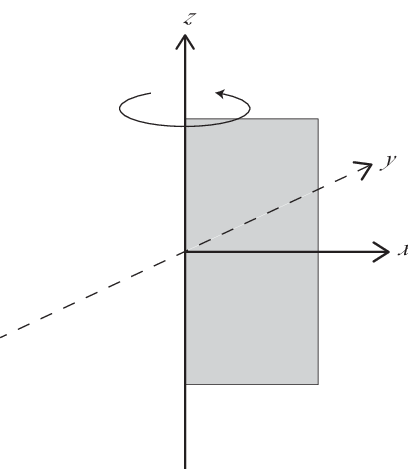}  
\bigbreak 
\bigbreak 
\bigbreak\caption{{\bf 
$\R^3$ is the result of rotating $\{(x, y,z) \vert x\geqq0, \quad  y=0 \}$ 
around the $z$-axis.}\label{rotz1}}   
\bigbreak 
\end{figure}

We regard $\R^3$ as the result of rotating $\{(x, y,z) \vert x\geqq0, \quad  y=0 \}$ 
around the $z$-axis. 
We draw the $x$-, $y-$, $z-$axis in a different fashion from what you are used to.  
See Figure \ref{rotz1}.

\np

\begin{figure}[H]
\bigbreak
\hskip41mm\includegraphics[width=5cm]{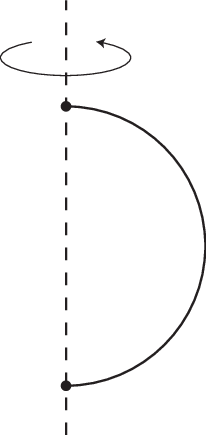}   
\bigbreak 
\bigbreak 
\bigbreak\caption{{\bf 
If we rotate the arc around the dotted line in $\R^3$, then the result is $S^2.$
}\label{arc1}}   
\bigbreak 
\end{figure}

See Figure \ref{arc1}. 
If we rotate the arc around the dotted line in $\R^3$, then the result is $S^2.$

\begin{figure}[H]
\bigbreak
\includegraphics[width=12cm]{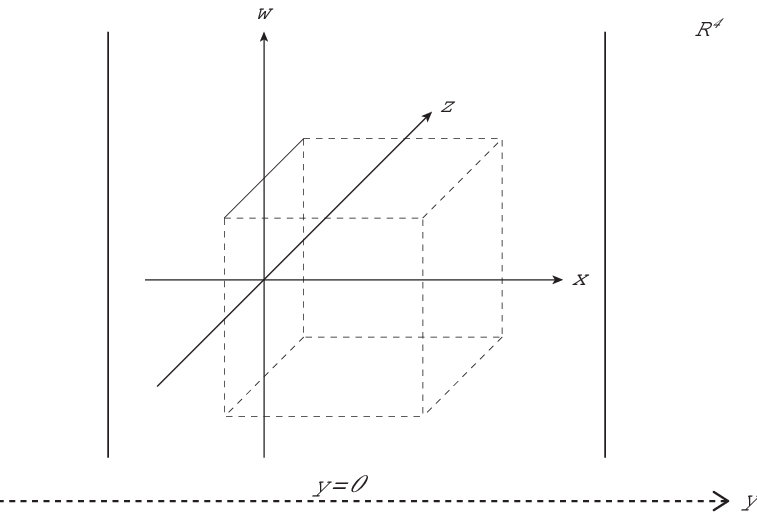}   
\bigbreak
\bigbreak
\caption{{\bf
We regard $\R^4$ as the result of rotating 
$\R^3_{\geq0}$
\newline
$=\{(x, y,z,w) \vert x\geqq0, \quad  y=0 \}$ 
around the $zw$-plane.  
}\label{rotzw1}}   
\bigbreak 
\end{figure}

\np
We regard $\R^4$ as the result of rotating \\
$\R^3_{\geq0}=\{(x, y,z,w) \vert x\geqq0, \quad  y=0 \}$ 
around the $zw$-plane. 
See Figure \ref{rotzw1}. Imagine!

\np

\begin{figure}[H]
\bigbreak
\hskip-10mm\includegraphics[width=170mm]{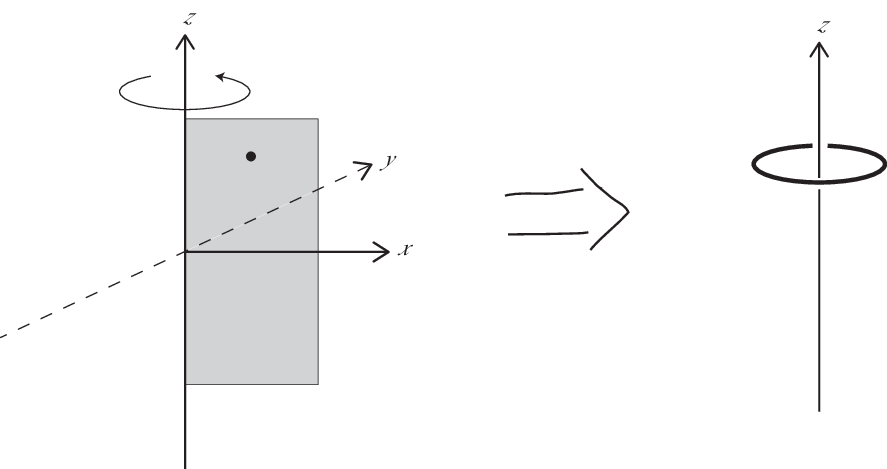}    
\bigbreak\caption{{\bf 
A circle is made by the rotation of a point around the $z$-axis. 
}\label{point1}}   
\bigbreak 
\end{figure}

Suppose that there is a point in $\{(x, y,z) \vert x\geqq0, \quad  y=0 \}$ as shown in Figure \ref{point1}.  
When we rotate $\{(x, y,z) \vert x\geqq0, \quad  y=0 \}$ 
around the $z$-axis,  
 rotate the point as well.  Then we obtain $S^1$ in $\R^3$.

%\np{\color{white}a}\vskip20cm

\np
Suppose that there is an arc in 

\begin{figure}[H]
\bigbreak
\hskip20mm\includegraphics[width=111mm]{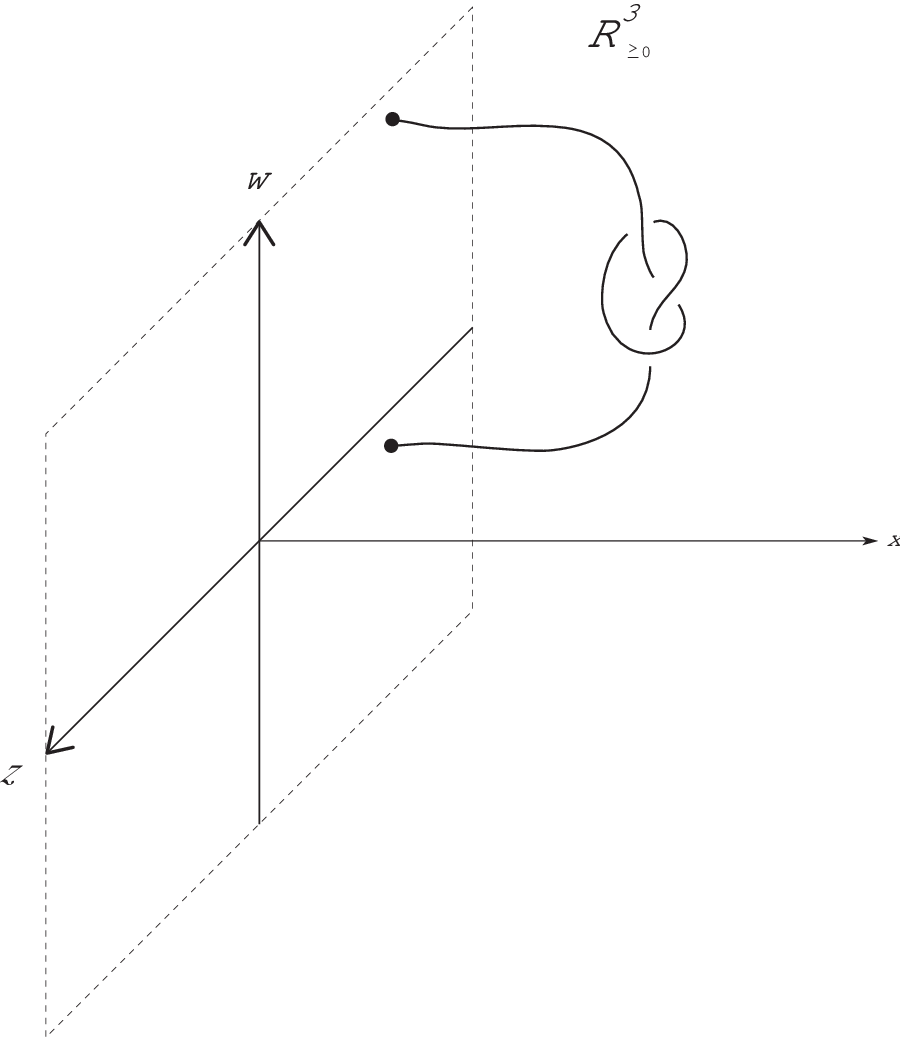}   
\bigbreak\caption{{\bf 
When we rotate $\{(x, y,z,w) \vert x\geqq0, \quad  y=0 \}$ 
around the $zw$-plane, 
 rotate the arc as well.  Then we obtain $S^2$ in $\R^4$.}\label{spun00}}   
\bigbreak 
\end{figure}

$\R^3_{\geq0}=\{(x, y,z,w) \vert x\geqq0, \quad  y=0 \}$ as shown in Figure \ref{spun00}.  

When we rotate $\{(x, y,z,w) \vert x\geqq0, \quad  y=0 \}$ 
around the $zw$-plane, 
 rotate the arc as well.  Then we obtain $S^2$ in $\R^4$.  Imagine!
Thus we obtain a 2-knot. 

%Note: We used the similar idea in Figure 078. 

\np

\begin{figure}[H]
\bigbreak
\includegraphics[width=12cm]{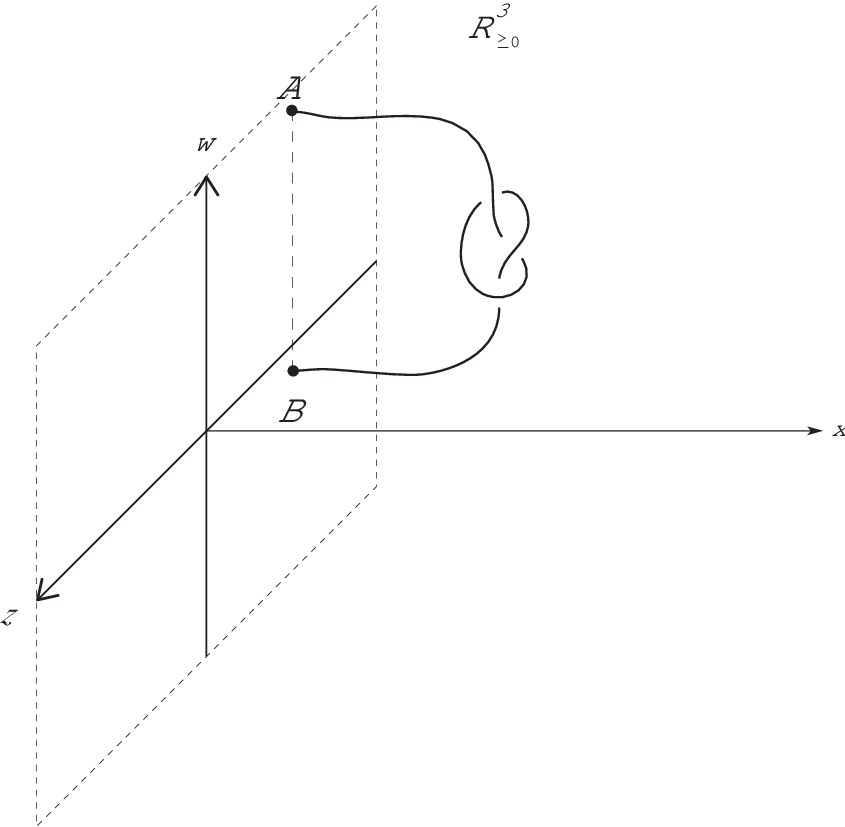}  
\bigbreak\caption{{\bf 
Making a spun knot from this.
}\label{spun2p}}   
\bigbreak 
\end{figure}

Let $A$ and $B$ be the endpoints of the arc. 
Suppose that $A$ and $B$ are in the $wz$-plane. 
Connect $A$ with $B$ by the dotted segment as shown in 
Figure \ref{spun2p}. 
Then a 1-knot is made from the arc and the dotted segment. Call this 1-knot $K$. 

The 2-knot made from this arc as shown in 
Figure \ref{spun2p} 
is called the {\it spun knot} of $K$. 

\bigbreak
\bigbreak

\np
It is known that if $K$ is a nontrivial 1-knot, then   
the spun-knot of $K$ is a nontrivial 2-knot. 

\vskip1cm
It is known that, there are countably infinitely many different 1-knots 
such that their spun knots are different each other.  
Thus there are countably infinitely many different 2-knots.  

\vskip1cm

We can define the spun-link of a 1-link as in the same manner. 
It is known that there are countably infinitely many different 1-links 
such that their spun links are different each other.  
Thus there are countably infinitely many different 2-links.  
\vskip1cm

In order to prove that they are nontrivial, 
we use the Alexander polynomial, the Alexander module, 
the fundamental groups, 
the $\mu$-invariant, the $\widetilde\eta$-invariant, etc.

\vskip1cm

See \cite{LevineOrr, Rolfsen, Ruberman, Zeeman}
 for detail.

\np

\begin{figure}[H]
\bigbreak
\includegraphics[width=120mm]{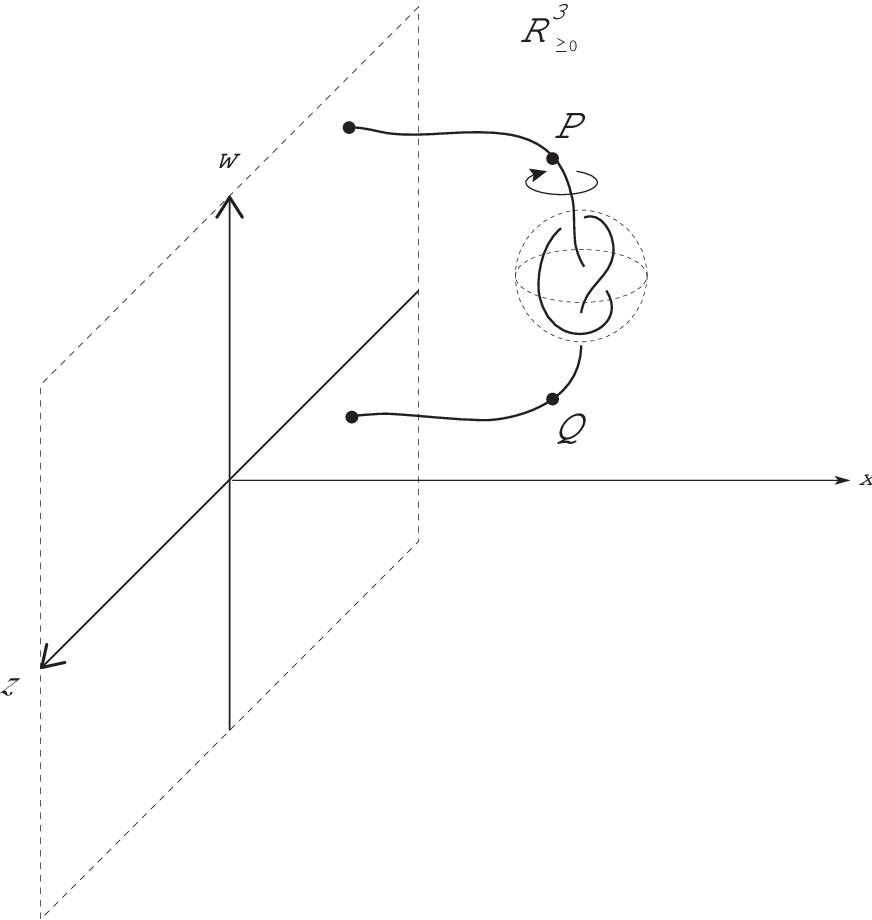}   
\bigbreak\caption{{\bf 
Making a twist-spun knot from this.
}\label{spun3p}}   
\bigbreak 
\end{figure}

When we rotate $R^3_{\geq0}=\{(x, y,z,w) \vert x\geqq0, \quad  y=0 \}$ 
and the arc as well  
around the $zw$-plane, 
rotate the part of the arc between $P$ and $Q$  $k$-times as shown in Figure \ref{spun3p}, 
where $k$ is an integer. 

Note that we rotate the dotted sphere as 
in Figure \ref{spun3p}. 
Associate this operation with `the revolution and the rotation of the planet'.  
%See Figure \ref{spun3p}

We obtain a 2-knot in $\R^4$. 

The result is called 
a {\it $k$-twist spun-knot} of $K$. 
Note that spun knots are 0-twist spun-knots.

\np
It is known that the 1-twist spun knot of any 1-knot is the trivial 2-knot. 

\vskip1cm
It is known that there are 
 countably infinitely many different 2-knots that are $k$-twist spun knots if $k\neq0,\pm1$.

In order to prove that they are nontrivial, 
we use the Alexander polynomial, the Alexander module, 
the fundamental groups, 
the $\mu$-invariant, the $\widetilde\eta$-invariant, etc.  

\vskip1cm
Furthermore it is known that there are countably infinitely many nontrivial 2-knots 
which are not $k$-twist spun knot for any integer $k$. 

\vskip1cm

See \cite{LevineOrr, Rolfsen, Ruberman, Zeeman}
 for detail. 

%\vskip1cm
%Furthermore it is known that there are countably infinitely many nontrivial 2-knots 
%which are not $k$-twist spun knot for any integer $k$. 

%\np

%\quad

%\vskip12cm\hskip8cm
%\includegraphics[width=6cm]{I11}    %zz-I11

%%%%%%%%%%%%%%%%%%%
\np

\subsection{$n$-dimsnional knots and links exist ($n\geqq3$)}
\label{ichisan}
We explain $n$-dimensional knots. % at the end of this Chapter \ref{exist}. 

\vskip1cm
We can define 
$n$-knots, $n$-links, 
the trivial $n$-knots, the trivial $n$-links,   
spun-knots of $(n-1)$-knots　(resp $(n-1)$-links), 
$k$-twist spun-knots of $(n-1)$-knots,  
as in a similar fashion to the 2-knot case.

\vskip1cm
If a single copy of $S^n$ is embedded in $\R^{n+2}$, it is called 
an {\it $n$-$($dimensional$)$ knot}. 

If a set of $m$ copies of $S^n$ is embedded in $\R^{n+2}$, 
it is called an {$m$-component $n$-$($dimensional$)$ link}.

If an $n$-knot $K$ bounds an embedded $n$-ball $B^n$ in $\R^{n+2}$, 
we say that $K$ is the {\it trivial knot}.

%If an $m$-component  $n$-link $L$ satisfies the condition that 
%each component bounds an embedded $n$-ball $B^n$ in $\R^{n+2}$ and that 
%each $B^n$ does not touch each other, we say that $L$ is the trivial link. 

Suppose that an $m$-component  $n$-link $L$ satisfies the following condition: 
Each component bounds an embedded $n$-ball $B^n$ in $\R^{n+2}$.   
Each $B^n$ does not touch any other $B^n$. 
Then we say that $L$ is the {\it trivial link}.

We say that two $n$-links are the same if we move one to the other smoothly 
in $\R^{n+2}$ 
such that it is always not touching itself, 
%so that it does not touch itself,  
%without touching itself, %self-touching, 
Here `moving smoothly' means that 
we move it with `tubular neighborhood' as in the 1-link case.

It is trivial that there is the $n$-dimensional trivial knot (resp. $n$-link).

\vskip1cm
\noindent{\bf To the advanced readers.}   
We usually define $n$-knots in $S^{n+2}$ rather than $\R^{n+2}$ in the smooth (resp. PL, topological) category  ($n\geqq1$).  
In the smooth category 
we usually define it to be a codimensision two smooth submanifold of $S^{n+2}$ 
which is PL homeomorphic to the standard $n$-sphere. 
It means that 
some $n$-knots for some natural numbers $n$
are exotic spheres. See 
\cite{Levineslice, Milnorexotic}.

\np

We regard $\R^{n+2}$ as the result of rotating \newline
$R^{n+1}_{\geq0}=\{(x_1,...,x_{n+2}) \vert x_1\geqq0, \quad  x_2=0 \}$ 
around \newline
the $x_3...x_{n+2}$-space
$=\{(x_1,...,x_{n+2}) \vert x_1=0, \quad  x_2=0 \}$.

Suppose that there is a $B^{n-1}$ in 

$R^{n+1}_{\geq0}=\{(x_1,...,x_{n+2}) \vert x_1\geqq0, \quad  x_2=0 \}$ 

\noindent 
and that 
the boundary of $B^{n-1}$ is in 

the $x_3...x_{n+2}$-space$=\{(x_1,...,x_{n+2}) \vert x_1=0, \quad  x_2=0 \}$.

\begin{figure}[H]
\bigbreak
\includegraphics[width=10cm]{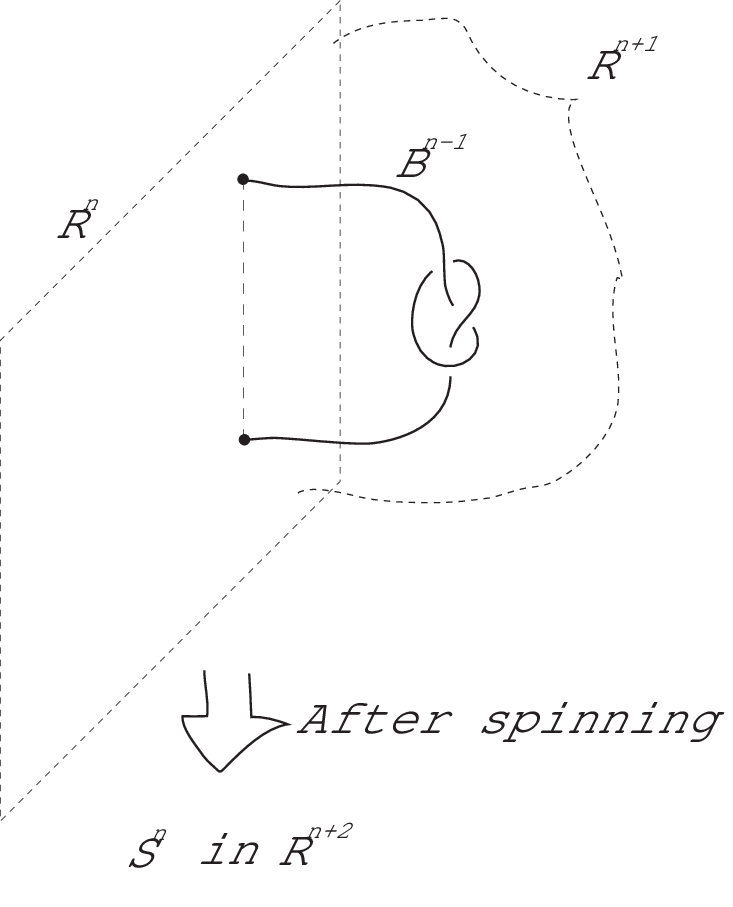}  
\bigbreak\caption{{\bf $n$-dimensional spun knots}\label{kyou}}   
\bigbreak 
\end{figure}

See Figure \ref{kyou}.
When we rotate 
$R^{n+1}_{\geq0}=\{(x_1,...,x_{n+2}) \vert x_1\geqq0, \quad  x_2=0 \}$ 
around \newline the $x_3...x_{n+2}$-space,  
rotate the $B^{n-1}$ as well.  
Then we obtain $S^n$ in $\R^{n+2}$.  Imagine!
Thus we obtain an $n$-knot.

We can make an $(n-1)$-knot $K$ 
from 
the $B^{n-1}$ in 

$R^{n+1}_{\geq0}=\{(x_1,...,x_{n+2}) \vert x_1\geqq0, \quad  x_2=0 \}$ 

\noindent 
as in a similar fashion to the case that we make spun-knots from 1-knots. 
See Figure \ref{kyou}.

The $n$-knot made from this $B^{n-1}$ is called 
the {\it spun-knot} of $K$.

%It is known that if $K$ is a nontrivial $(n-1)$-knot, 
%then the spun-knot of $K$ is a nontrivial $n$-knot. 

\np
It is known that there are countably infinitely many different $(n-1)$-knots 
such that their spun knots are different  each other.  
Thus there are countably infinitely many different $n$-knots ($n\geqq3$).  

\vskip1cm
We can define the spun-link of an $(n-1)$-link as in the same manner. 
It is known that there are countably infinitely many different $(n-1)$-links 
such that their spun links are different each other.  
Thus there are countably infinitely many different $n$-links.  

\vskip1cm
See \cite{LevineOrr, Rolfsen, Ruberman, Zeeman}
 for detail.

\np

\begin{figure}[H]
\bigbreak
\includegraphics[width=13cm]{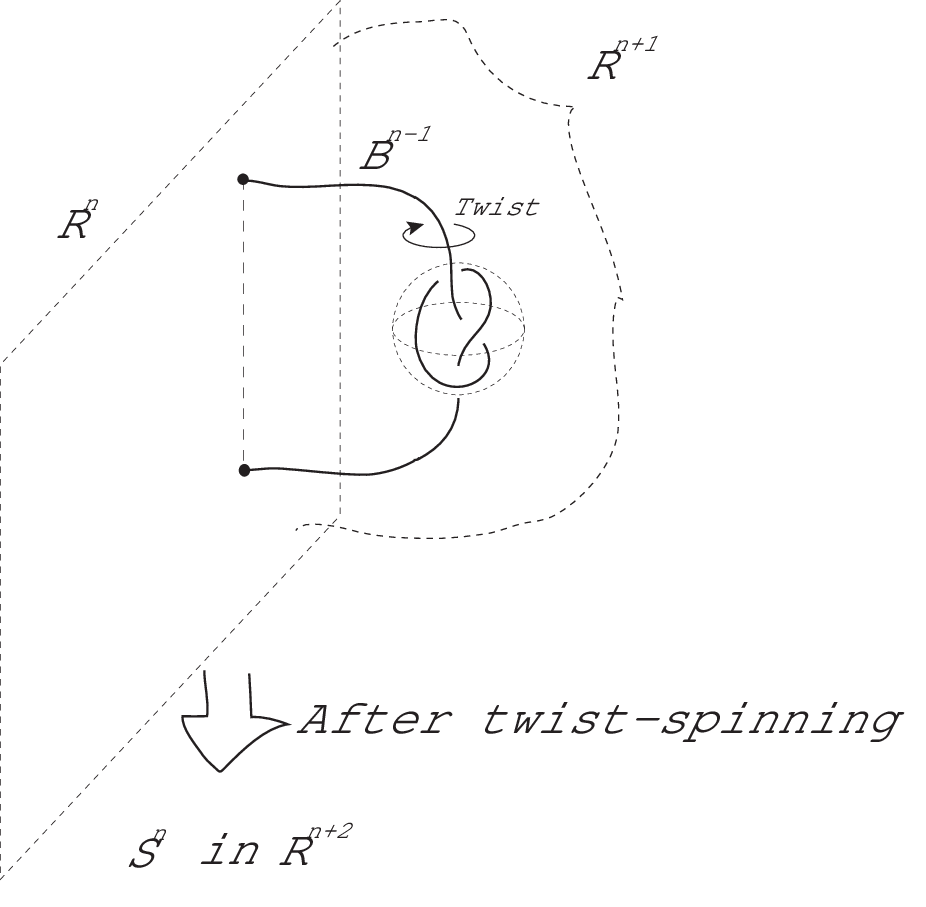}  
\bigbreak\caption{{\bf $n$-dimensional spun knots}\label{asu}}   
\bigbreak 
\end{figure}

We can define {\it $k$-twist spun-knot} of $K$ as follows. 
When we rotate $B^{n-1}$ around the $x_3...x_{n+2}$-space, 
we rotate a part of $B^{n-1}$ as in a similar fashion to the case 
that we make $k$-twist spun knots of 1-knots. 
See Figure \ref{asu}.

\np
It is known that 
the 1-twist spun knot of any $(n-1)$-knot is the trivial $n$-knot. 

\vskip1cm

It is known that there are 
 countably infinitely many different $n$-knots 
which are $k$-twist spun knots for each $k\neq0,\pm1$.

In order to prove that they are nontrivial, 
we use the Alexander polynomial, the Alexander module, 
the fundamental groups, 
%the $\mu$-invariant, 
the $\widetilde\eta$-invariant, 
the signature, %the Arf invariant, 
knot cobordism groups, etc.

\vskip1cm
Furthermore it is known the following: Let $n$ be an integer $\geqq3$. There are countably infinitely many nontrivial $n$-knots 
which are not $k$-twist spun knot for any integer $k$. 

\vskip1cm
See \cite{LevineOrr, Rolfsen, Ruberman, Zeeman}
 for detail. 

\vskip1cm
\noindent{\bf To the advanced readers.} 
If an $n$-knot $K$ is an exotic sphere, 
the $k$-twist spun knot of $K$ is also defined.

\vskip1cm
\noindent{\bf To the advanced readers.} 
We can also discuss $S^p$ in $\R^n$ even if $n>p+2$. 
See \cite{Haefliger}  for detail.

\np

\section{The projections of $n$-dimensional knots}\label{proj}
In this section, we explain `the projections of $n$-knots'. 
%Associated with this topic, we will introduce an application of the Boy surface.  
%%as previously mentioned %preannounced  in \S\ref{Boy}. 
%%Before we do so, we state the following fact.

\begin{figure}[H]
\bigbreak
\includegraphics[width=12cm]{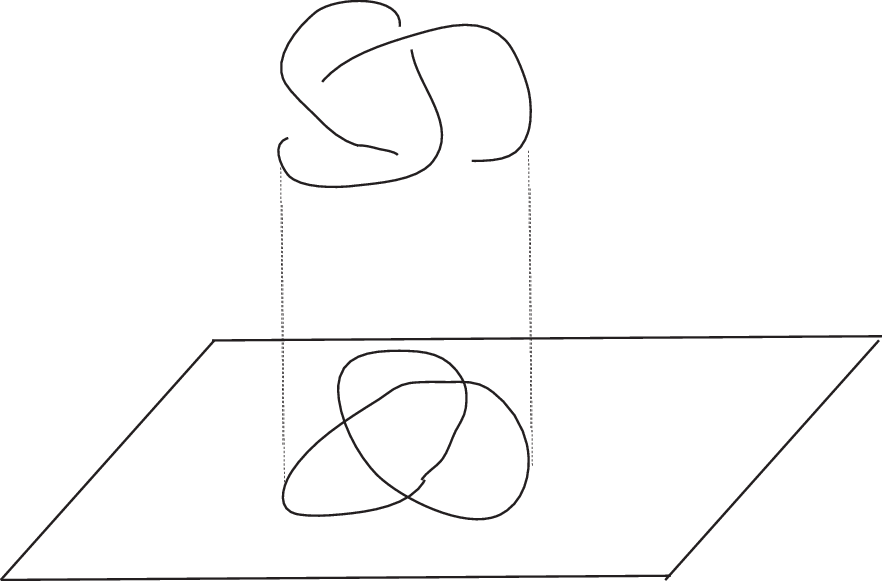}   
\bigbreak\caption{{\bf 
The projection of a 1-dimensional knot}
\label{I1232pq}}   
\bigbreak 
\end{figure}

The {\it projection} of a 1-knot is the image of the knot by a map 
from $\R^3$ to $\R^2$:\\ $(x,y,z)\mapsto(x,y)$.  %See \S\ref{linking}. 
See Figure \ref{I1232pq} for an example. 

\np

\begin{figure}[H]
\bigbreak
\includegraphics[width=11cm]{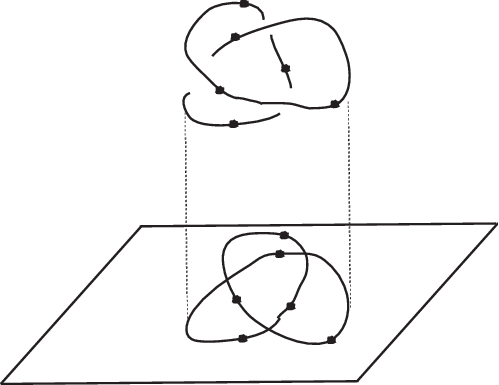}   
\bigbreak\caption{{\bf 
$S^1$ is divided into a finite numbed of 
embedded curved segments.}
\label{waketa}}   
\bigbreak 
\end{figure}

\begin{que}\label{mondai1}
Suppose that $S^1$ is in $\R^2$. 
It may touch itself. %have self-touching. 

Suppose the following conditions. 

\begin{enumerate}
\item 
$S^1$ can be divided a finite number of curved segments so that each curved segment is embedded 
(hence each curved segment does not touch itself). 
An example is drawn in Figure \ref{waketa}.

\item
If it touches itself, %If it has self-touching,
it touches itself like the intersection %the part of touching itself is the intersection like 
\includegraphics[width=7mm]{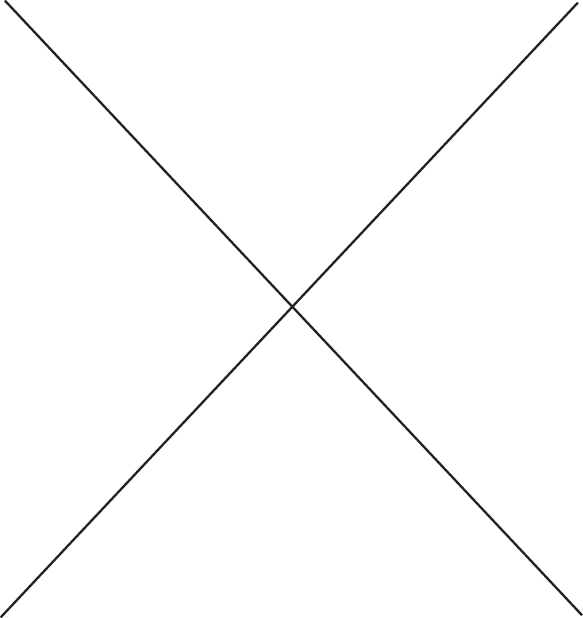}: 
Two curved segments intersect. 
\end{enumerate}

Then is this the projection of a 1-knot in $\R^3$?
\smallbreak 
\end{que}
Here, we say that $S^1$ is 
{\it immersed transversely}
into $\R^2$.

\vskip5mm\noindent 
Answer:  Yes.   
%
%\vskip5mm
%It is known.  
%The proof is easy. Please try. 
(The proof is easy. Please try.)

\np

\begin{figure}[H]
\bigbreak
\hskip30mm
\includegraphics[width=27mm]{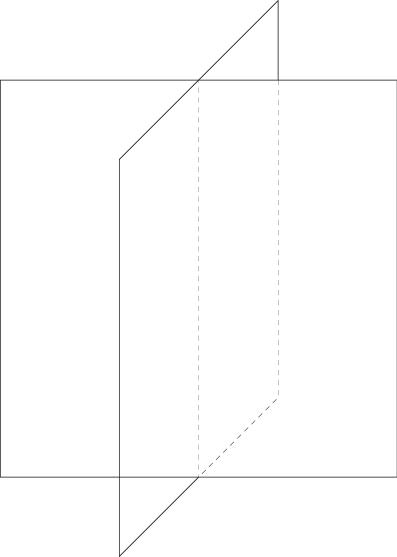}  
\hskip10mm
\includegraphics[width=36mm]{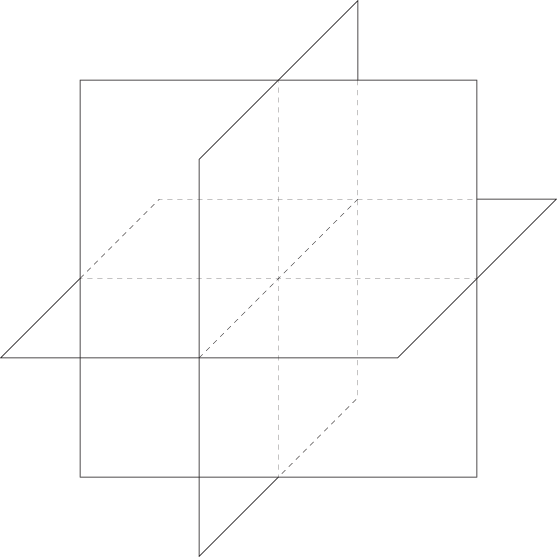}  
\bigbreak\caption{{\bf 
Two or three sheets intersect.
}\label{dotr}}   
\bigbreak 
\end{figure}

Take $S^2$ in $\R^3$ 
so that 
$S^2$ can be divided into a finite number of curved 2-dimensional discs each of which is embedded in $\R^3$.  
We suppose that $S^2$ may touch itself 
and that if so, it touches itself 
as shown 
in Figure \ref{dotr}.

Here, two or three sheets intersect. 
The sheets can be bended a little.

Then we say that $S^2$ is 
{\it immersed transversely} into $\R^3$  
and that 
the $S^2$ which touches itself 
is a transverse immersion of $S^2$. 
See textbooks on manifolds for transverse immersion.

\bigbreak
The {\it projection} of a 2-knot is the image of the knot by a map 
from $\R^4$ to $\R^3$: \\ $(x,y,z,w)\mapsto(x,y,z)$.  %See 

\bigbreak
\begin{que}\label{mondai2}
For an arbitrary transverse immersion of $S^2$ into $\R^3$, 
is there a 2-knot in $\R^4$ whose projection is the transverse immersion?  
\end{que}

\bigbreak

Question \ref{mondai2} is a one-dimensional-higher 
analogue of Question \ref{mondai1}.

\np
\noindent
Answer: No.  (\cite{Giller} proved.)   

\vskip1cm
Of course, there are countably infinitely many transverse immersions of $S^2$ into $\R^3$ 
which are the projections of 2-knots. 
Here, we ask whether any  transverse immersion is the projection of a 2-knot.

%\vskip1cm
%The idea of his proof is that 
%the Boy surface in $\R^3$ is not 
%the image of any embedded $\R P^2$ in $\R^4$
%by the projection map $\R^4\to\R^3$: $(x,y,z,w)\mapsto(x,y,z)$. 
%We will explain it 
%further %a little more 
%in the following pag.  

\bigbreak

\begin{figure}[H]
\bigbreak
\includegraphics[width=100mm]{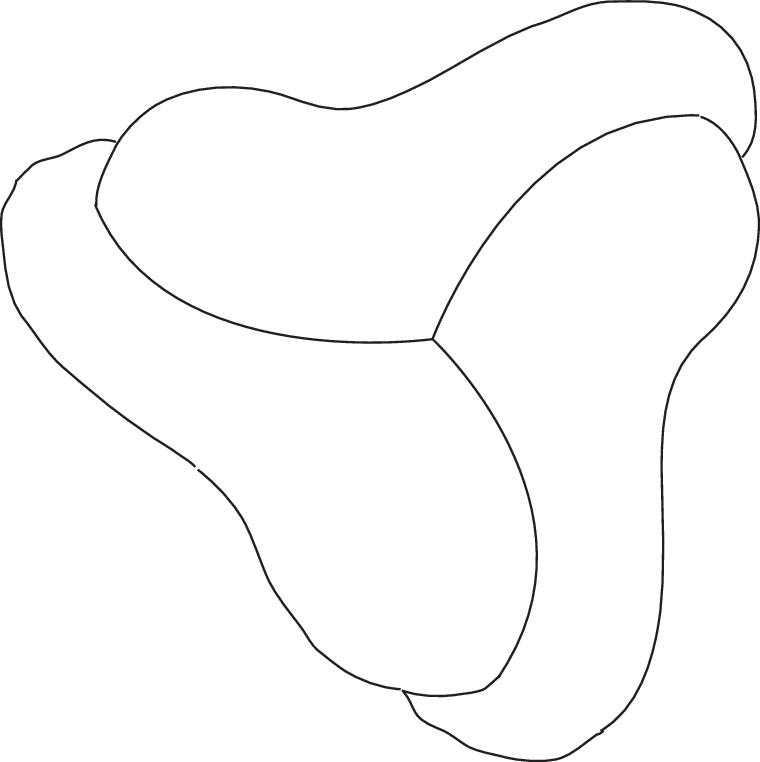}  
\bigbreak\caption{{\bf 
Boy surface
}\label{bosf}}   
\bigbreak 
\end{figure}

The idea of making an  transverse immersion of $S^2$ into $\R^3$ 
which is not the projection of any 2-knot in $\R^4$ is as follows. 
We use Boy surface. See Figure \ref{bosf}.

\bigbreak
In his paper \cite{Boy} 
Boy found a figure, which we call Boy surface.  
P121 of 
\cite{MilnorStasheff} 
quotes Boy's paper. 

\bigbreak 
In \cite{OgasaBoy} 
the author have created a way to construct the Boy surface 
by using a pair of scissors, a piece of paper, and a strip of scotch tape. 

We made a movie of our paper-craft. We put it in 

\bigbreak
\noindent
http://www.geocities.jp/n\_dimension\_n\_dimension/MakeyourBoysurface.html
\bigbreak

\noindent Don't forget the three \_ in the address  if you type it. 

It is connected with the author's website  

\bigbreak
\noindent
http://www.geocities.jp/n\_dimension\_n\_dimension/list.html

\bigbreak
\noindent
Click the indication: The movie `Make your Boy surface'  in the first few lines in his one.

You can find his website by typing in the author's name `Eiji Ogasa' in the search engine 
although you will not type the address of the website. 

\np
Take the Boy surface in $\R^3$ 

\begin{figure}[H]
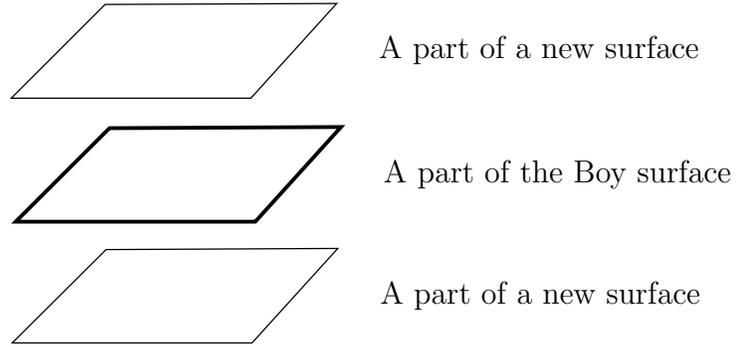

\input I234.tex
\vskip40mm
\bigbreak\caption{{\bf 
Take two sheets for a sheet which is a part of Boy surface.}\label{ninosan}}   
\bigbreak 
\end{figure}

Take two sheets at each point in the Boy surface 
as shown in Figure \ref{ninosan}.

For a point in the self-intersection of the Boy surface,  do the same procedure. 
%do as follows. 
%%There are two or three parts of the Boy surface which include the point. 
%The point is included in two or three parts of the Boy surface.   
%For the two or three parts, take two sheets as shown above, too.   
%Note that the new sheets intersect each other around the point. 
%
%
%
%The union of such all sheets makes an immersion of $S^2$ into $\R^3$. 
%Recall that a relation between $\R P^2$ and $S^2$.   %See \S\ref{ano}. 
%
%This immersion of $S^2$ into $\R^3$ is not the projection of any 2-knot. 
%
The union of such all sheets makes an  transverse immersion of $S^2$ into $\R^3$ 
which is not the projection of any 2-knot.

\np

It is natural to ask the following question, which is a high dimensional version of 
Question \ref{mondai2}.

\begin{que}\label{mondai3easy}
Take $S^n$ in $\R^{n+1}$.  
$S^n$ may touch itself.  

Is there an $n$-knot in $\R^{n+2}$ whose projection is the $S^n$ in $\R^{n+1}$?  
\end{que}

\bigbreak
The {\it projection} of an $n$-knot $K$ is the image of $K$ by the map  
from $\R^{n+2}$ to $\R^{n+1}$; $(x_1,...,x_{n+1},x_{n+2})\mapsto(x_1,...,x_{n+1}).$

\vskip10mm
\noindent{\bf To the advanced readers.} 
Consider the following question.
\begin{que}\label{mondai3}
Take a transverse immersion $S^n$ into $\R^{n+1}$. 
Is there an $n$-knot in $\R^{n+2}$ whose projection is the $S^n$ in $\R^{n+1}$?  
\end{que}

\bigbreak
\bigbreak
A {\it transverse immersion} $S^n$ into $\R^{n+1}$ is a natural generalization of the lower dimensional one: 
Let $P_i=\{(x_1,...,x_{n+1})|x_i=0\}$. 
$S^n$ can be divided into 
a finite number of embedded curved $n$-dimensional ball $D^n$. 
If a point is in 
the self-intersection of $S^n$, 
the part near $P$ is like 
`$P_1$ and $P_2$',    
`$P_1$ and $P_2$ and $P_3$', 
,...,or
`$P_1$,..., and $P_{n+1}$'. 

See textbooks on manifolds for detail.

\np

In his paper \cite{Ogasa99S}, 
the author   found that the answers 
to Question \ref{mondai3easy} and \ref{mondai3}
are negative for any integer $\geqq3$. 

\vskip1cm
The author made such a transverse immersion of $S^n$ into $\R^{n+1}$ in an explicit way.   
%(Note that a transverse  immersion of manifolds in $n$-dimensional case is a mathematical term. See textbooks on manifolds  for detail.) 

The idea of construction:  
%is a generalization of the earlier idea of the 2-sphere case where we use the Boy surface. 
We rotate the  transverse immersion  of $S^{n-1}$ as making spun-knots and 
obtain the  transverse immersion of $S^n$ for $n\geqq3$. 
The $n=2$ case is the example which is made earlier by using the Boy surface.

\vskip1cm
This is the sub-theme of the author's paper 
\cite{Ogasa99S}.

\vskip1cm
The author also found the following. 
It is the main theme of the above paper. 

\np

\noindent 

First, consider an easy question. 

\begin{que}\label{mondaitsuika}
Take a transverse immersion of $S^1$ into $\R^2$. 
%It may touch itself. Suppose the  conditions in Question \ref{mondai1}. 

Then is this the projection of a trivial 1-knot in $\R^3$?
\smallbreak 
\end{que}

Note that the difference between 
Questions \ref{mondai1} and \ref{mondaitsuika}.

\vskip1cm
It is almost trivial that the answer is positive. 
Please prove so. 

\vskip10mm

\begin{figure}[H]
\bigbreak
\includegraphics[width=15cm]{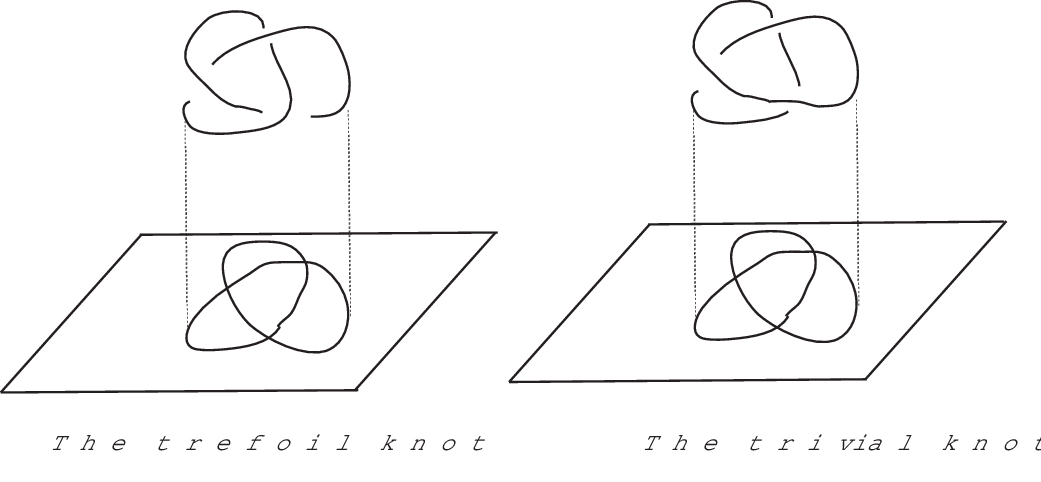} 
\bigbreak\caption{{\bf 
A nontrivial 1-knot whose projections is the projection of 
a trivial knot.
}\label{tanpuku}}   
\bigbreak 
\end{figure}

An example is drawn in Figure \ref{tanpuku}.

\np

Consider the $n$-dimensional version of Question \ref{mondaitsuika}.

\begin{que}\label{mondai4easy}
Take  $S^n$ in $\R^{n+1}$. 
$S^n$ may touch itself. 

Is there a trivial $n$-knot in $\R^{n+2}$ whose projection is 
the $S^n$ in $\R^{n+1}$?  
\end{que}

Note that the difference between 
Questions \ref{mondai3easy} and \ref{mondai4easy}.

\vskip10mm
\noindent{\bf To the advanced readers.}
Consider the following question. 

\begin{que}\label{mondai4}
Take  an arbitrary transverse immersion of $S^n$ into $\R^{n+1}$. 

Is there a trivial $n$-knot in $\R^{n+2}$ whose projection is 
the transverse immersion of $S^n$ into $\R^{n+1}$?  
\end{que}

\noindent 
Note that the difference between 
Questions \ref{mondai3} and \ref{mondai4}.

\np
The author found that 
the answers to 
Questions \ref{mondai4easy} and \ref{mondai4}
are negative for any integer $\geqq3$.

\vskip1cm
The author made a  transverse  immersion of $S^n$ into $\R^{n+1}$ 
which is the projection of a nontrivial $n$-knot but 
which is not the projection of any trivial $n$-knot 
in an explicit way. 
We used `the K3 surface', which you see in the textbooks on algebraic geometry. 

\vskip1cm
The $n=2$ case of Question \ref{mondai4} is open. Please solve.

\vskip1cm\noindent 
{\bf To the advanced readers:} 
The singularity of the transverse immersions  
in the author's papers \cite{Ogasa01P, Ogasa99S} 
consists of only double points. 
If $n\geqq5$, 
the connected components of the singularity  
of the author's example in 
\cite{Ogasa99S} 
are two.

In the $n=4,3$ (or $=2$),  
what is the least number of connected components of the singularity set? 
Please solve.

\vskip1cm   
See  the author's  papers
\cite{Ogasa01P, Ogasa99S}.

\np

\section{Local moves on $n$-dimensional knots}\label{localmoves}

\begin{figure}[H]
\bigbreak
\hskip15mm \includegraphics[width=10cm]{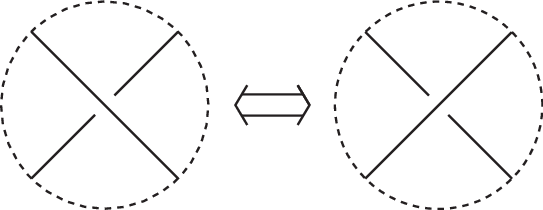}   
\bigbreak\caption{{\bf The crossing change of a 1-knot}\label{tcr1}}   
\bigbreak 
\end{figure}

\begin{figure}[H]
\bigbreak
\hskip15mm  \includegraphics[width=10cm]{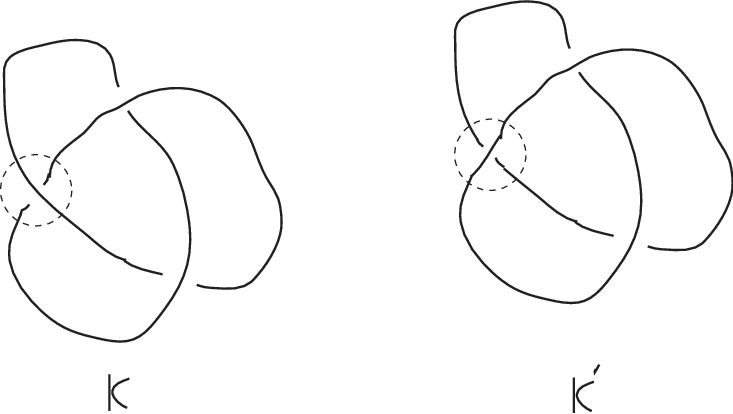}   
\bigbreak
\caption{{\bf $K$ is obtained from $K'$ by one crossing-change.}\label{AMPM}}   
\bigbreak 
\end{figure}

If two 1-dimensional links $K$ and $K'$ in $\R^3$ differ only in the 3-ball $B^3$ as shown 
in Figure \ref{tcr1},

\noindent 
then we say that $K$ is obtained from $K'$ by one  crossing-change. 
An example is shown in 
Figure \ref{AMPM}.

\np
Let $K_1,...,K_\mu$ be 1-links. 
If $K_\nu$ is obtained from $K_{\nu-1}$ by one  crossing-change\\ $(2\leqq\nu\leqq\mu)$, 
then we say that $K_\mu$ is obtained from $K_1$ 
by a sequence of a finite number of crossing-changes. 

It is known that 
any $m$-component 1-link is obtained from the trivial  $m$-component 1-link  
by a sequence of a finite number of crossing-changes. The proof is easy. Please try. 

\vskip1cm
It is important to consider such a procedure that 
we change a knot by `local moves' as shown 
in Figure \ref{tcr1}.
 
\vskip1cm
{\Large 
How about `local moves on high dimensional knots'? }

\bigbreak

It is the main theme of this section. %Chapter \ref{localmoves}. 

\np
From here on we will 
show an example that crossing-changes make a nontrivial link into 
the trivial link. 

\begin{figure}[H]
\bigbreak
\hskip30mm\includegraphics[width=9cm]{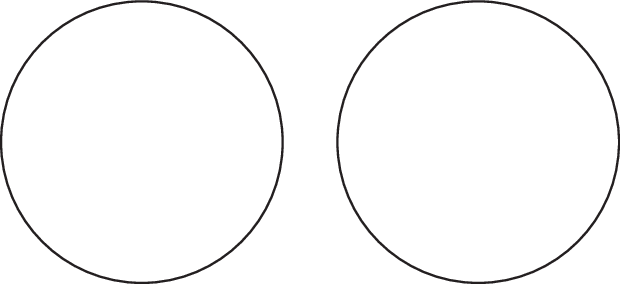} 
\bigbreak\caption{{\bf }\label{trilinzu}}   
\bigbreak 
\end{figure}

The 1(-dimensional) link  
in Figure \ref{trilinzu}

\noindent
is called the (1-dimensional 2-component) trivial link. 

\begin{figure}[H]
\bigbreak
\hskip30mm\includegraphics[width=8cm]{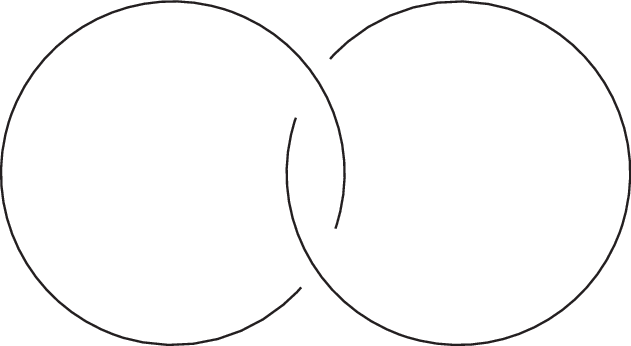} 
\bigbreak\caption{{\bf }\label{Hopfzu}}   
\bigbreak 
\end{figure}

The 1-link 
in Figure \ref{Hopfzu}

\noindent
is called the Hopf link. 
%Hopf is the name of a mathematician who studied this link. 

\np 
The two 1-links have a different property. 

\vskip1cm
Take a disc $D^2\subset\R^3$. 

We suppose that the disc does not touch itself. 
(Then we say that this disc is embedded in $\R^3$. 
See textbooks on manifolds for detail. 
We usually suppose the condition on relative topology furthermore.)

The boundary of $D^2$ is one of the two components of each 1-link. 
Here, we do not care whether $D^2$ touches the other component of the link or not. 

\vskip1cm
Can we take $D^2$ so that $D^2$ does not touch the other component of the link? 

\vskip1cm
In the case of the trivial link the answer is affirmative. 

In the case of the Hopf link the answer is negative. 

Do you feel they are true? Indeed, it is true. 

We can prove it 
by using the fundamental group or the homology groups of \\$\R^3-$(one of the two $S^1$). 
%See text books on knot theory and topology.

\np

\begin{figure}[H]
\bigbreak
\hskip20mm\includegraphics[width=8cm]{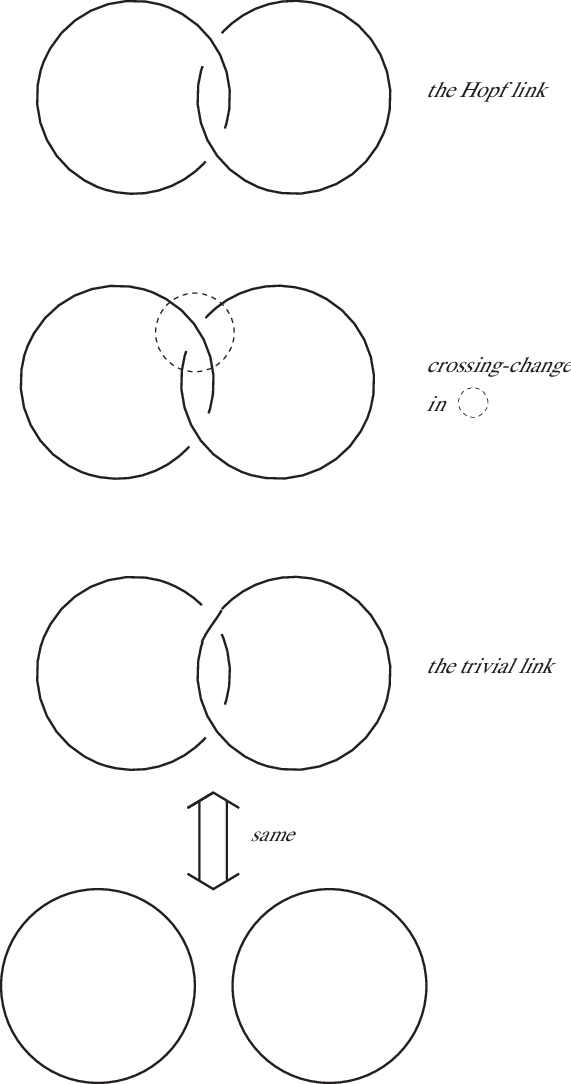}   
\bigbreak\caption{{\bf }\label{Hopftotrizu}}   
\bigbreak 
\end{figure}

A crossing-change makes the Hopf link into the trivial link as shown in Figure \ref{Hopftotrizu}.

%We will show a `one dimensional higher' analogue of this example from here on. 
We will generalize this example and consider a higher dimensional case from here on. 
%Recall that $L_0$ and $L_1$ which are 
%$S^1$ and $S^2$ in $\R^4$. See \S\ref{S12R4}, \S\ref{L0L1}. 

\np
%\begin{sub}\label{ribbonmoves}
%Ribbon moves on $n$-dimensional knots
%\end{sub}

%\part{tameshi}\label{tameshi}

\subsection{Ribbon moves on 2-dimensional knots}
\label{ribbonmoves}

\subsubsection{Ribbon moves on 
a circle and a sphere in the 4-space}
\label{ribb12}

We construct two types, $L_0$ and $L_1$, of  `$S^1$ and  $S^2$ in $\R^4$’ 
such that $S^1$ does not touch $S^2$ 
and that $S^1$ (resp. $S^2$) does not touch itself.

\vskip1cm
\noindent
$L_0$:  $S^1$ bounds an embedded 2-dimensional disc $D^2$ in $\R^4$. 
        
$S^2$  bounds an embedded 3-dimensional ball $B^3$ in $\R^4$. 
        
$B^3$ does not touch $D^2$. 

\begin{figure}[H]
\includegraphics[width=12cm]{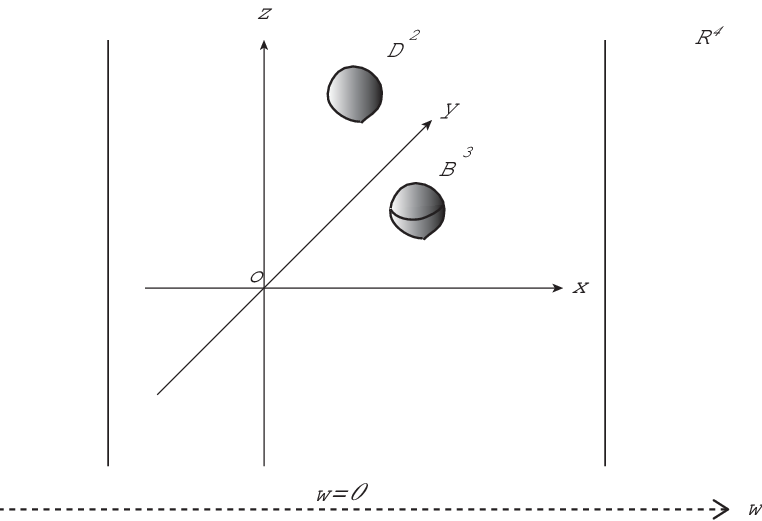}   
\bigbreak
\bigbreak
\caption{{\bf  $L_0$}\label{yojigen0}}   
\end{figure}

\vskip1cm
%We think that you feel that it is trivial that $L_0$ exists. 
Probably you will feel  $L_0$ exists. 
Indeed, it is true. 
Take 
$D^2$ and $B^3$ in $\R^3$ disjointly.  
Regard $\R^4$ as the the $xyzw$-space. 
Regard $\R^3$ as the $xyz$-space in the $xyzw$-space. 
See Figure \ref{yojigen0}.

\np 
\noindent
$L_1$:  $S^1$ bounds an embedded 2-dimensional disc $D^2$ in $\R^4$  
if  we do not 

care whether $D^2$ touches $S^2$ or not. 

$S^2$  bounds an embedded 3-dimensional ball $B^3$ in $\R^4$ 
if  we do not 

care whether $B^3$ touches $S^1$ or not. 

$D^2$ definitely touches $S^2$ even if we take $D^2$ in any way. 
        
We can take $D^2$ such that $D^2$ touches $S^2$ at a single point. 

\noindent{\bf To the advanced readers.} We suppose that they intersect transversely. 

\vskip1cm
In fact, such $L_1$ exists. 

Can you imagine $L_1$? 

Can you construct $L_1$? 

\vskip1cm
We construct $L_1$ from here on.

\vskip1cm
Note that $L_0$ and $L_1$ are the name used only in this article.

\np
 
Regard $\R^4$ as the result of rotating 

$R^3_{\geq0}={(x,y,z,t) \vert  x\geqq0, y=0 }$  around the $zw$-plane. 

(Recall the way when we define spun-knots 
in Figure \ref{spun2p}.)

\begin{figure}[H]
\bigbreak
\includegraphics[width=131mm]{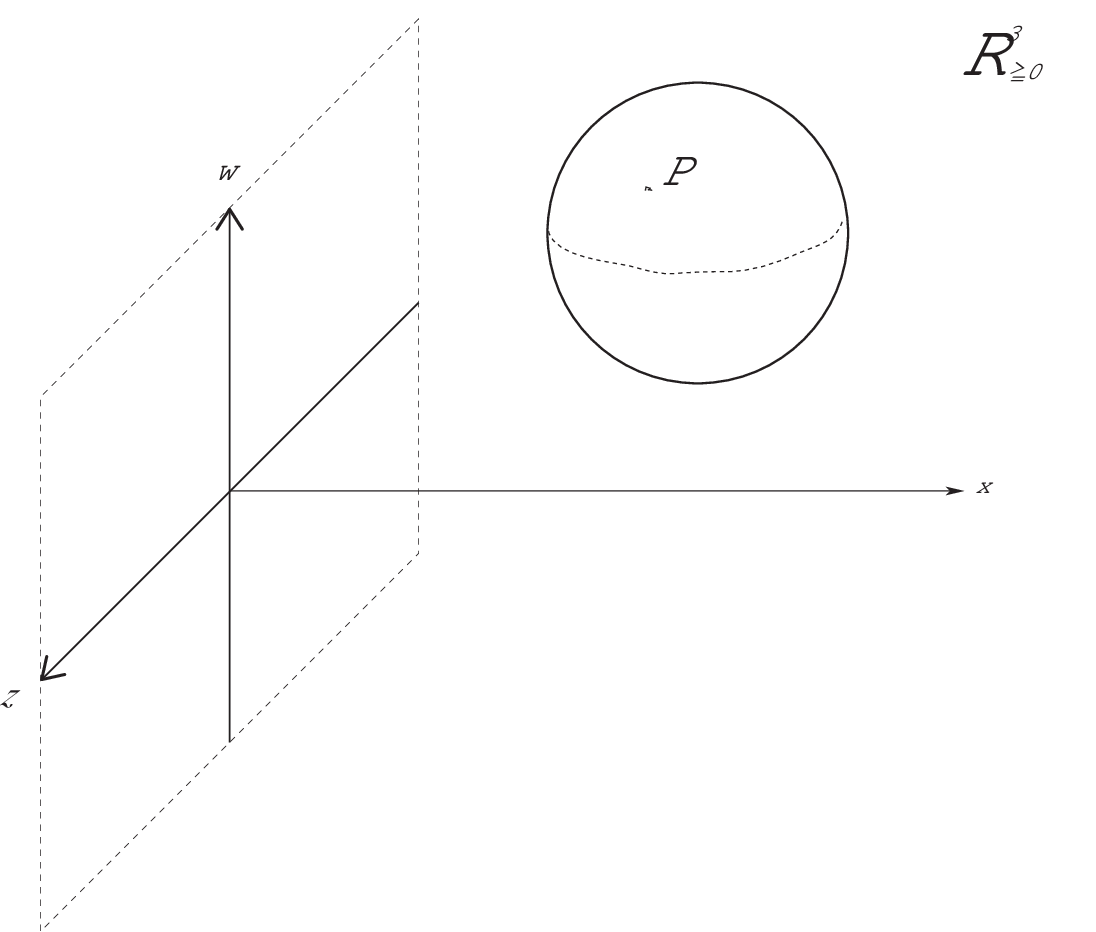}   
\bigbreak\caption{{\bf }\label{L1S1S2zu}}   
\bigbreak 
\end{figure}

Take $S^2$ in $R^3_{\geq0}$ 
as shown in Figure \ref{L1S1S2zu}. 
Take an embedded 3-dimensional ball $B^3$ in $R^3_{\geq0}$. 
Take a point $P$ in $B^3$ such that $P$ is not included in $S^2$. 
When we rotate $R^3_{\geq0}$ around the $zw$-plane, 
rotate $P$ as well but stop $S^2$ at the first place. 
Then $P$ becomes a circle $S^1$. 
This situation is drawn as shown 
in Figure \ref{conce1zu}. 

\np

\begin{figure}[H]
\bigbreak
\hskip10mm\input  conce1.tex
\bigbreak\caption{{\bf }\label{conce1zu}}   
\bigbreak 
\end{figure}

\vskip12mm
This $S^1$ and $S^2$ makes $L_1$. 

We can prove that this $L_1$ satisfies the required condition 
%by using  the homology groups of $R^4-S^2$ or $R^4-S^1$. 
by using the homology groups or the fundamental group of $\R^4-S^2$ 
or by using the homology groups of $\R^4-S^1$. 
%See text books on knot theory and topology.

\np
The figures drawn from here on 
will request the readers to have imagination  more than 
those until here. 

Imagine!

\bigbreak
\bigbreak

\begin{figure}[H]
\bigbreak
\includegraphics[width=123mm]{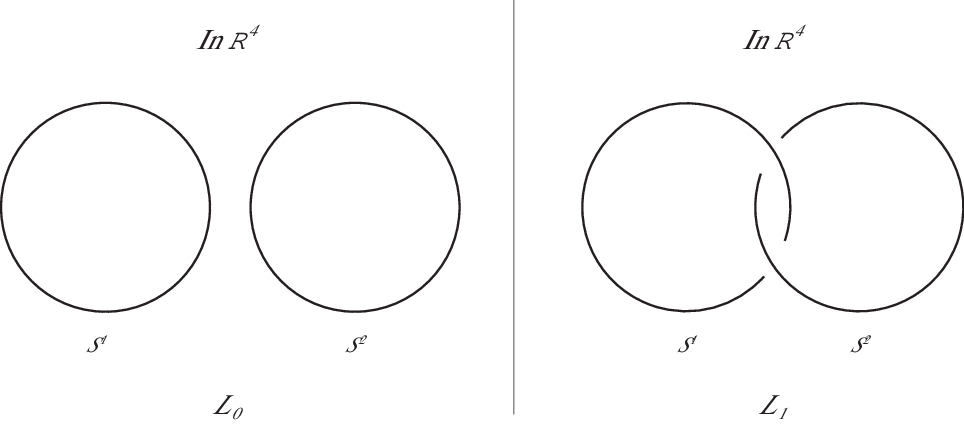}  
\bigbreak\caption{{\bf }\label{L1L0zu}}   
\bigbreak 
\end{figure}

\bigbreak
\bigbreak

\begin{figure}[H]
\bigbreak
\includegraphics[width=10cm]{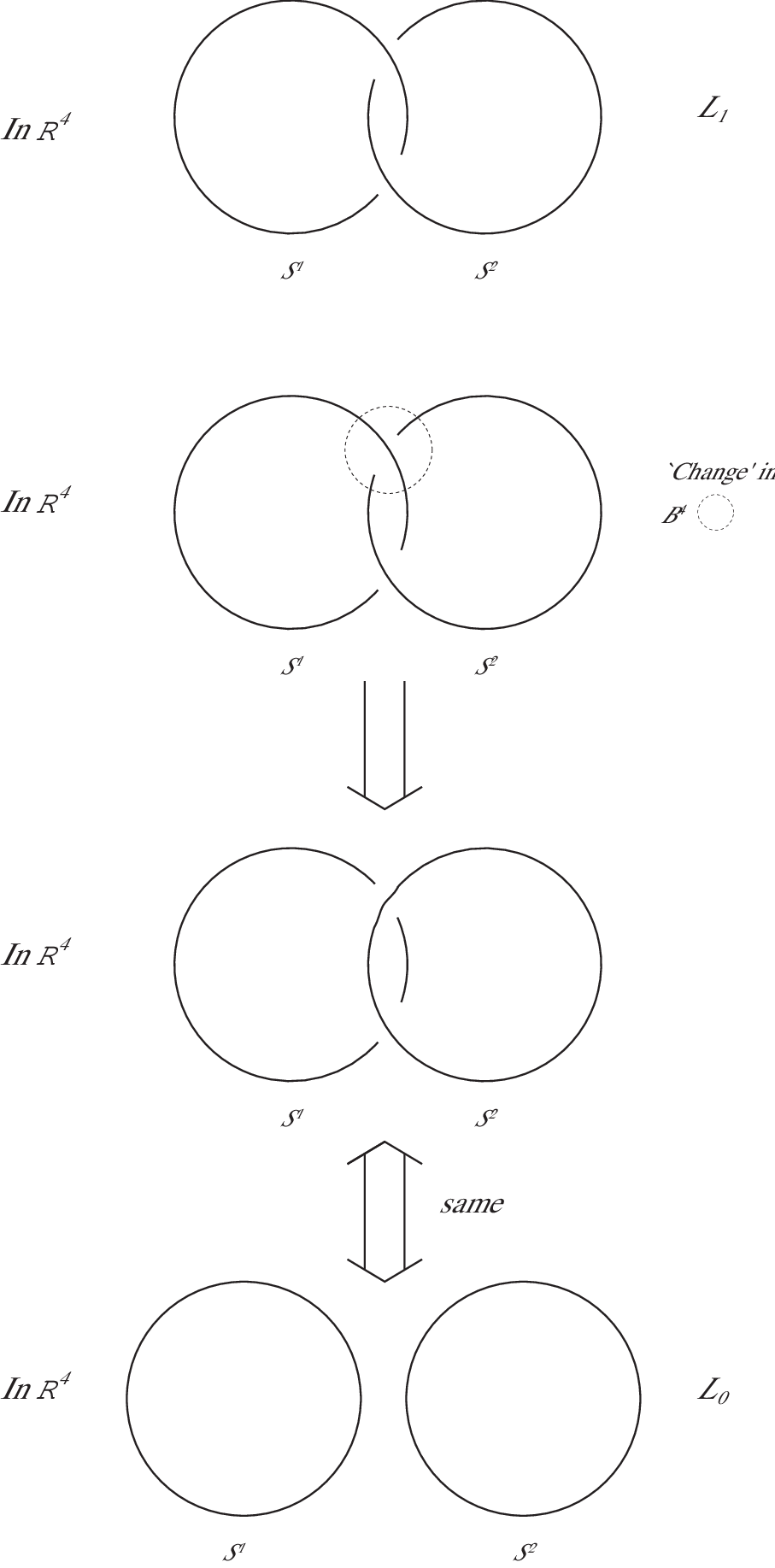}  
\bigbreak\caption{{\bf }\label{L1toL0zu}}   
\bigbreak 
\end{figure}

$L_1$ in Figure \ref{L1L0zu} 
is changed into $L_0$ in Figure \ref{L1L0zu}  
by a local move as shown 
in Figure \ref{L1toL0zu}.

\np

\np

\begin{figure}[H]
\bigbreak
\includegraphics[width=15cm]{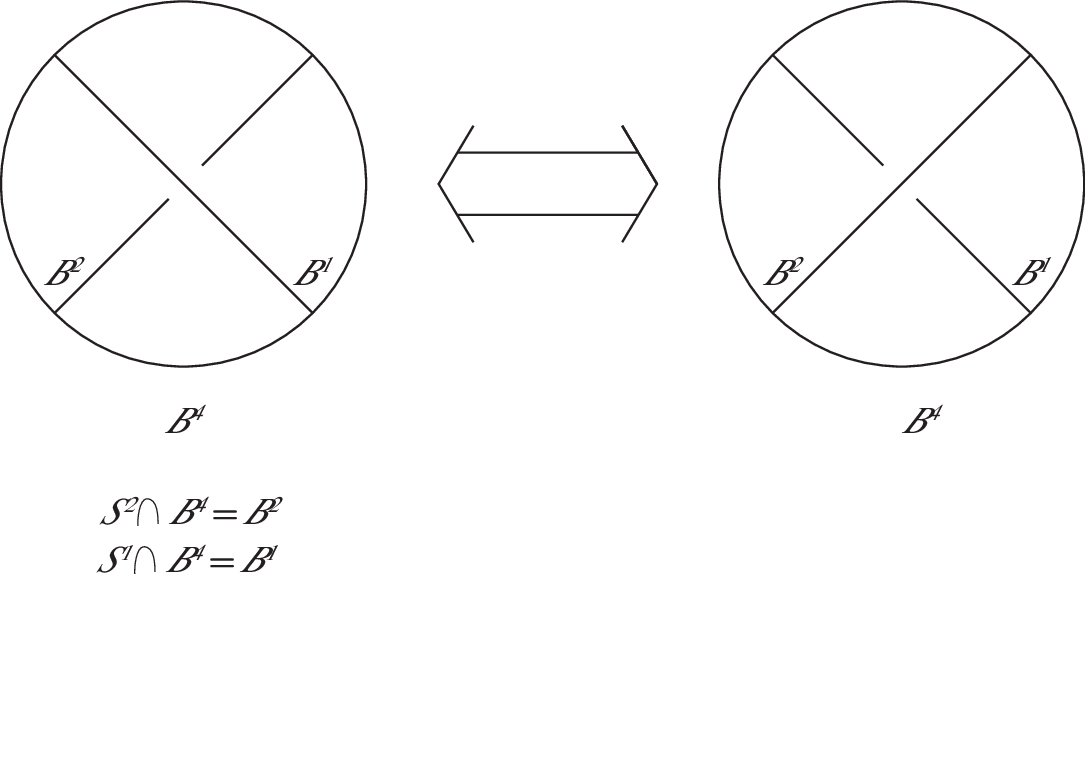}   
\bigbreak\caption{{\bf }\label{ribbonzu}}   
\bigbreak 
\end{figure}

`The local move' in $B^4$ 
%
%\begin{figure}[H]
%\bigbreak
\includegraphics[width=7mm]{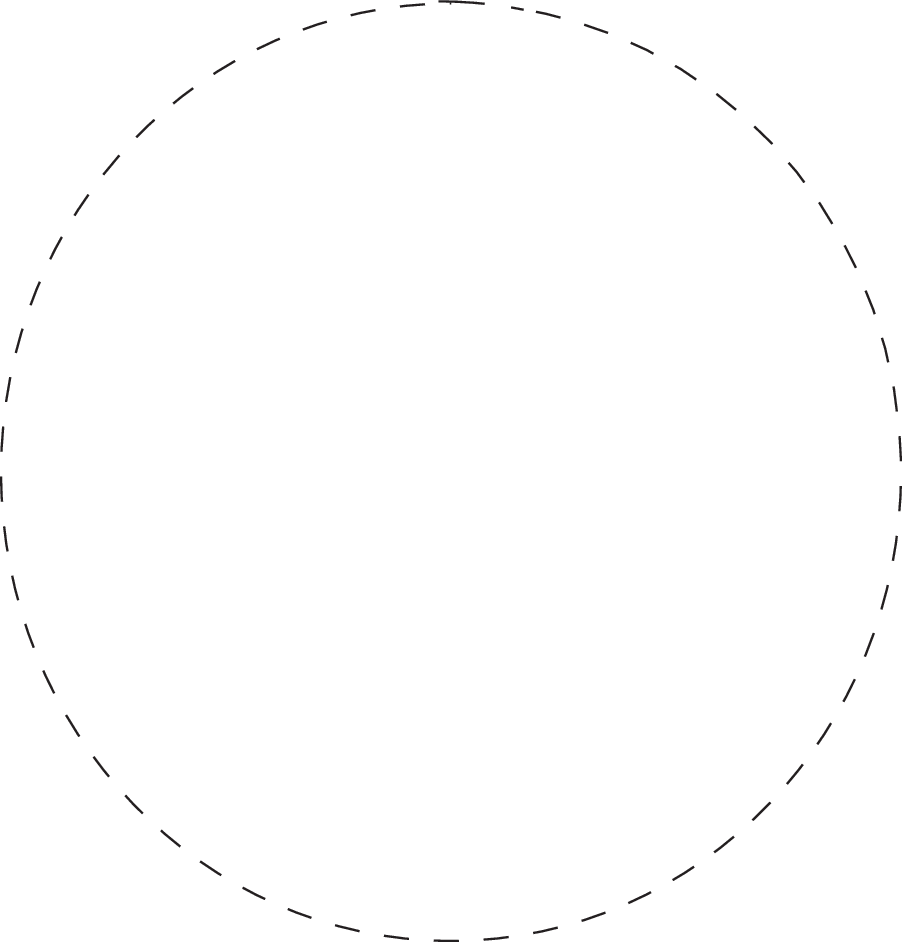}   
%\bigbreak\caption{{\bf }\label{}}   
%\bigbreak 
%\end{figure}
%
in Figure \ref{L1toL0zu}  
is drawn as shown 
in Figure \ref{ribbonzu}. 
Note $B^n$ is the $n$-dimensional ball ($n\in\N$).  
This procedure is called the {\it ribbon-move}. 

We say that 
$L_1$ (resp. $L_0$) is changed into $L_0$ (resp. $L_1$) 
by one ribbon-move.

%\vskip1cm   We will draw a little more rigorous figure of this procedure in Figures \ref{Figure1au} and \ref{Figure2au}.    

%\np

\vskip10mm

\np

We draw a little more explicit figure of the ribbon-move 
in Figures \ref{Figure1au} and \ref{Figure2au}.

\begin{figure}[H]
\bigbreak
\input Figure1.tex
\bigbreak 
\bigbreak 
\bigbreak\caption{{\bf }\label{Figure1au}}   
\end{figure}

\np

\begin{figure}[H]
\bigbreak
\input Figure2.tex
\bigbreak 
\bigbreak 
\bigbreak\caption{{\bf }\label{Figure2au}}   
\end{figure}

\np
\noindent  {\bf To  the advanced readers.}    
We will write the definition of  ribbon-moves rigorously 
in Definition \ref{ribodef1to2}  
by using Figures \ref{Figure1au} and \ref{Figure2au}.

We use some mathematical terms, so please check the literature. 

Although we will write the definition by using coordinates,  
it is better go in 
the way (i) for the beginners 
than in the way (ii).  

\smallbreak\noindent
(i) To perceive it intuitively the first time around.   
After that, understand the definition by using coordinates. 

\smallbreak\noindent
(ii)To try to understand it by coordinate representation the first time around.

\begin{defn}\label{ribodef1to2}
Let $K_1$ and $K_2$ be sets of $S^1$ and $S^2$ in $\R^4$. 
We say that $K_2$ is obtained from $K_1$ by one {\it ribbon-move } 
if there is a 4-ball $B$ embedded in $\R^4$ with the following properties.  

\vskip3mm
(1) 
$K_1$ coincides with $K_2$ in $\R^4-(\text{the interior of }B)$. 
This identity map from $K_1-(\text{the interior of }B)$ to 
$K_2-(\text{the interior of }B)$  
is orientation preserving.

(2) 
$B\cap K_1$ is drawn as shown in     
Figure \ref{Figure1au}. 

$B\cap K_2$ is drawn as shown in 
Figure \ref{Figure2au}.

\noindent 
Note: $\cap$ denotes the intersection.  
See textbooks on set theory for the intersection of sets. 
\vskip3mm

%%%%%%%%%%%%%%%%%%

We regard $B$ as 
(a closed 2-disc)$\times[0,1]\times\{t| -1\leqq t\leqq1\}$.
We suppose $B_t=$(a closed 2-disc)$\times[0,1]\times\{t \}$.  
Then $B=\cup B_t$. 
In 
Figures \ref{Figure1au} and \ref{Figure2au}  
we draw $B_{-0.5}, B_{0}, B_{0.5}$ in $B$. 
We draw $K_1$ and $K_2$ by the bold line. 
The fine line denotes $\partial B_t$. 
Here,  $\partial X$ means the boundary of $X$. 
  
$B\cap K_1$ (resp. $B\cap K_2$) is diffeomorphic to 
$D^2\amalg (\text{a segment} [0,1])$.  

\noindent Note: $\amalg$ denotes the disjoint union. 
See textbooks on set theory for the disjoint union.

Note that $[a,b]$ means $\{ x \vert a\leqq x\leqq b \}$.

\vskip1cm
$B\cap K_1$ has the following properties:  
$B_t\cap K_1$ is empty for $-1\leqq t<0$ and $0.5<t\leqq1$.
$B_0\cap K_1$ is diffeomorphic to 
$D^2\amalg( [0,0.3])\amalg( [0.7,1])$. 
$B_{0.5}\cap K_1$ is diffeomorphic to $( [0.3,0.7])$. 
$B_t\cap K_1$ is diffeomorphic to two points for $0<t<0.5$. 
(Here, we draw $[0,1]$ to have the corner 
in $B_0$ and in $B_{0.5}$. 
Strictly to say, $B\cap K_1$ in $B$ is a smooth embedding 
that is obtained by making the corner smooth.)

\vskip1cm
$B\cap K_2$ has the following properties:  
$B_t\cap  K_2$ is empty for $-1\leqq t<-0.5$ and $0<t\leqq1$.
$B_0\cap K_2$ is diffeomorphic to 
$D^2\amalg( [0, 0.3])\amalg( [0.7, 1])$. 
$B_{-0.5}\cap  K_2$ is diffeomorphic to $( [0.3, 0.7])$. 
$B_t\cap  K_2$ is diffeomorphic to two points for $-0.5<t<0$. 
(Here, we draw $[0,1]$ to have the corner 
in $B_0$ and in $B_{-0.5}$. 
Strictly to say, $B\cap K_1$ in $B$ is a smooth embedding 
that is obtained by making the corner smooth.)

\np

In Figure \ref{Figure1au} (resp. \ref{Figure2au}) 
there are an oriented segment $[0,1]$ 
and an oriented disc $D^2$ as we stated above. 
We do not make any assumption about 
the orientation of the segment and the disc. 
The orientation of $B\cap K_1$ (resp. $B\cap K_2$) 
coincides with that of the segment and that of the disc.

\vskip1cm

If $K_1$ is obtained from $K_2$ by one ribbon-move,  
then we also say that $K_2$ is obtained from $K_1$ by one {\it ribbon-move}.

\vskip1cm

Let $K_1$ and $K_2$  be two sets of $S^1$ and $S^2$ in $\R^4$. 
$K_1$ and $K_2$ are said to be {\it ribbon-move-equivalent} 
if there are sets of $S^1$ and $S^2$ in $\R^4$,  
$K_1=\bar{K}_1, \bar{K}_2,...,\bar{K}_{r-1},\bar{K}_r=K_2$  
 ($r,p$ are natural numbers and $p\geq2$) such that 
$\bar{K}_i$ is obtained from $\bar{K}_{i-1}$ $(1< i\leqq r)$ by one ribbon-move. 

\end{defn}

\np
\subsubsection
{Ribbon moves on 
a sphere and a torus in the 4-space}
\label{ribbTS}

Next we consider `$S^2$ and $T^2$ in $\R^4$'. 
Here,  $T^2$ denotes the torus. 
%See \S\ref{ano} and \S\ref{S2S1xS1} for tori. 

\vskip1cm

The figure in Figure \ref{torusau}   is called the {\it torus}. 
Note: 
the torus is like an inner tube. The interior of an inner tube is empty. 

\bigbreak 
\bigbreak 

\begin{figure}[H]
\bigbreak
\hskip20mm\input torus.tex
\bigbreak\caption{{\bf The torus}\label{torusau}}   
\bigbreak 
\end{figure}
\bigbreak 
\bigbreak

If the interior is filled, then it is called the {\it solid torus}.

\np
We construct two types, $L_0$ and $L_1$, of  
`$S^2$ and  $T^2$ in $\R^4$’ 
such that $T^2$ does not touch $S^2$ 
and that $T^2$ (resp. $S^2$) does not touch itself.

%It is known that there are countably infinitely many 
%different types of `$S^2$ and $T^2$ in $\R^4$'. 

\vskip10mm

\noindent
$L_0$:  
$S^2$ bounds a 3-dimensional ball $B^3$ embedded in $\R^4$. 

%so that  $B^3$ does not touch %avoids 
%$T^2$. 

$T^2$ bounds the solid torus  embedded in $\R^4$.  
%so that the solid torus does not touch %avoids 
%$S^2$. 

$B^3$ does not touch the solid torus. 

\vskip10mm
\noindent
$L_1$:  
$S^2$ bounds a 3-dimensional ball $B^3$ embedded in $\R^4$,  
where we do 

not care whether $B^3$ touches $T^2$ or not.

$T^2$ bounds the solid torus  embedded in $\R^4$ so that the solid torus 

does not touch %avoids 
$S^2$.

$B^3$ definitely touches the solid torus even if we take them in any way. 

We can take $B^3$ and the solid torus so that  
their intersection is a 

single 2-dimensional ball.

\vskip10mm
%We think that you feel that it is trivial that $L_0$ exists. 
Probably you will feel  $L_0$ exists. 
Indeed, it is true. 

You only have to 
replace 
$D^2$ with a solid torus in Figure \ref{yojigen0}. 

\vskip10mm

Can you imagine $L_1$? 

Can you construct $L_1$?

\vskip10mm
Note that $L_0$ and $L_1$ are the name used only in this article.

\np
We construct $L_1$ from here on.

\begin{figure}[H]
\bigbreak
\includegraphics[width=140mm]{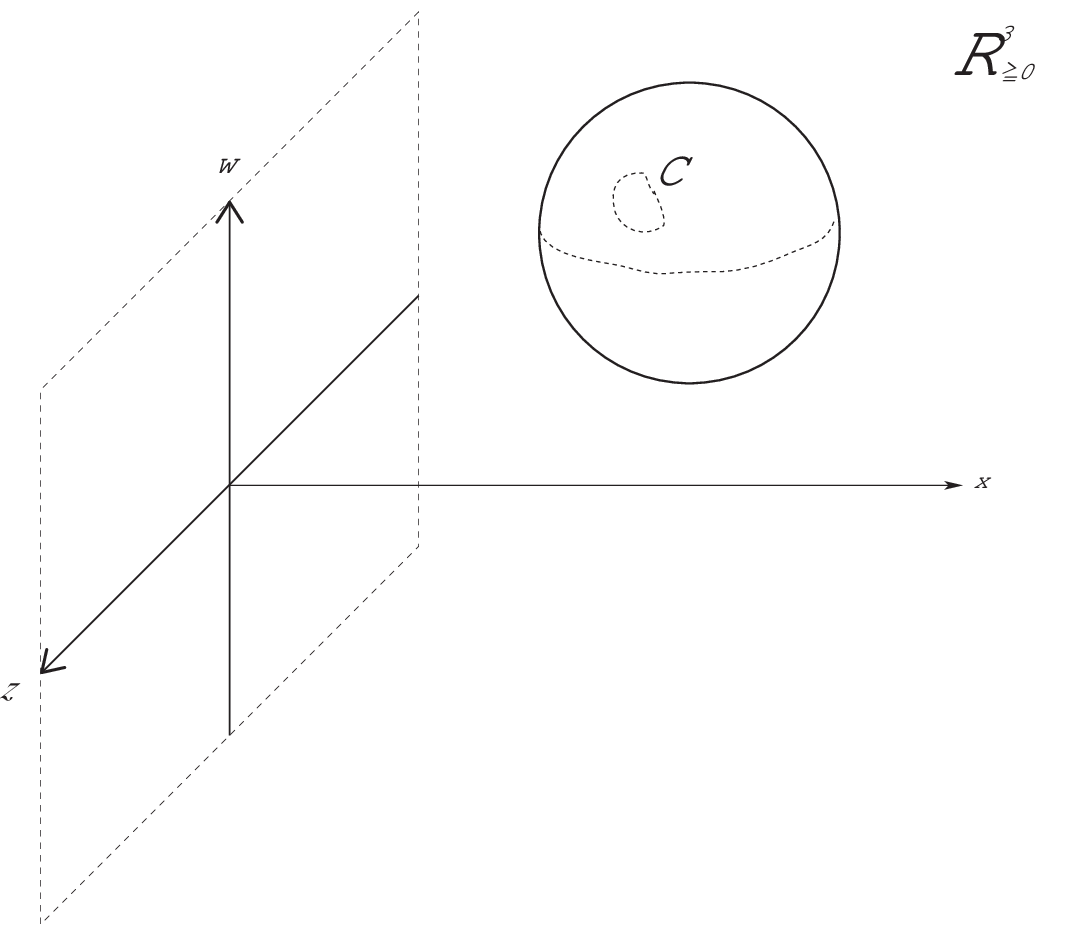}   
\bigbreak\caption{{\bf }\label{L1T2S2au}}   
\end{figure}

Recall the construction of $L_1$ of `$S^1$ and $S^2$ in $R^4$’ 
in Figures \ref{L1S1S2zu} and \ref{conce1zu}. 
Take a circle $C$ instead of the point $P$ in $B^3$ such that $C$ does not touch $S^2$. See Figure \ref{L1T2S2au}.

\np
When we rotate $R^3_{\geq0}$ around the $zw$-plane, 
rotate $C$ as well but stop $S^2$ at the first place. 
Then $C$ becomes a torus $T^2$. 
See Figure \ref{conce2au}.

\begin{figure}[H]
\bigbreak
\input  conce2.tex
\bigbreak\caption{{\bf }\label{conce2au}}   
\bigbreak 
\end{figure}

This $T^2$ and $S^2$ makes $L_1$. 

\np

\begin{figure}[H]
\bigbreak
\includegraphics[width=140mm]{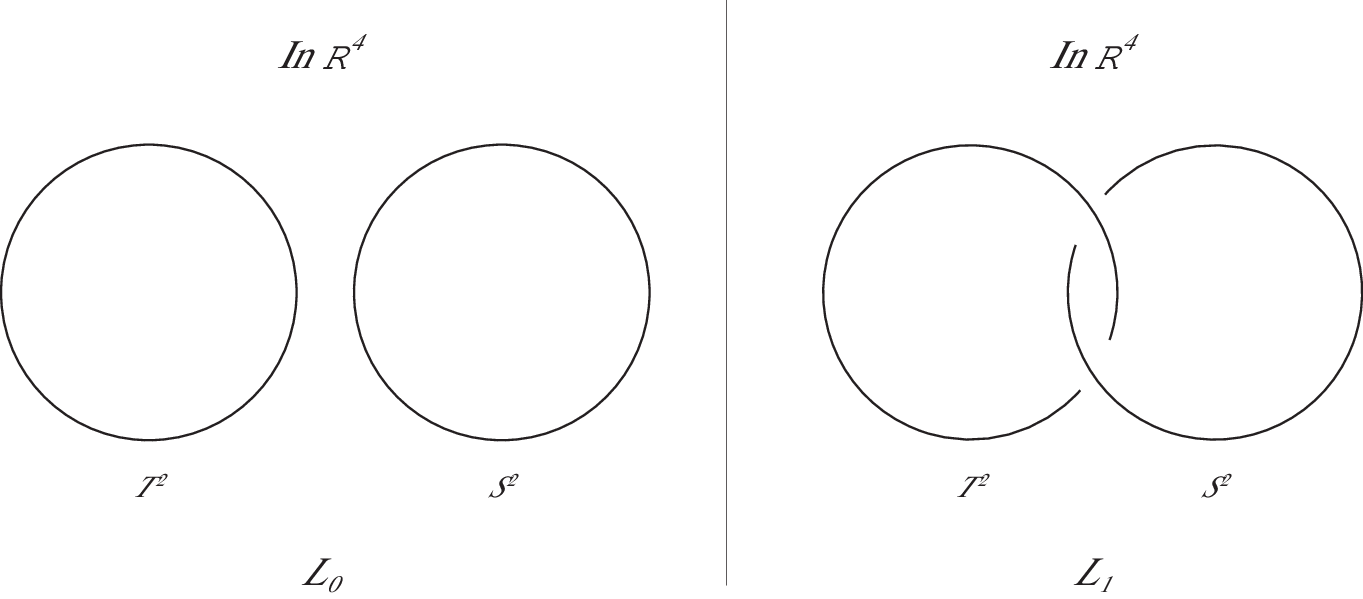}
\bigbreak\caption{{\bf }\label{T2a00}}   
\bigbreak 
\end{figure}

We draw $L_0$ and $L_1$ in Figure \ref{T2a00}. 

\bigbreak

\bigbreak

\vskip1cm

It is known that $L_0$ and $L_1$ are different. 
%You can prove it by using Chapter \ref{proof}. Please try. 
We can  prove it  
by using the homology groups or the fundamental group of $\R^4-S^2$ 
or by using the homology groups of $\R^4-T^2$. 
%See textbooks on knot theory and topology.

\np

\begin{figure}[H]
\bigbreak
\includegraphics[width=140mm]{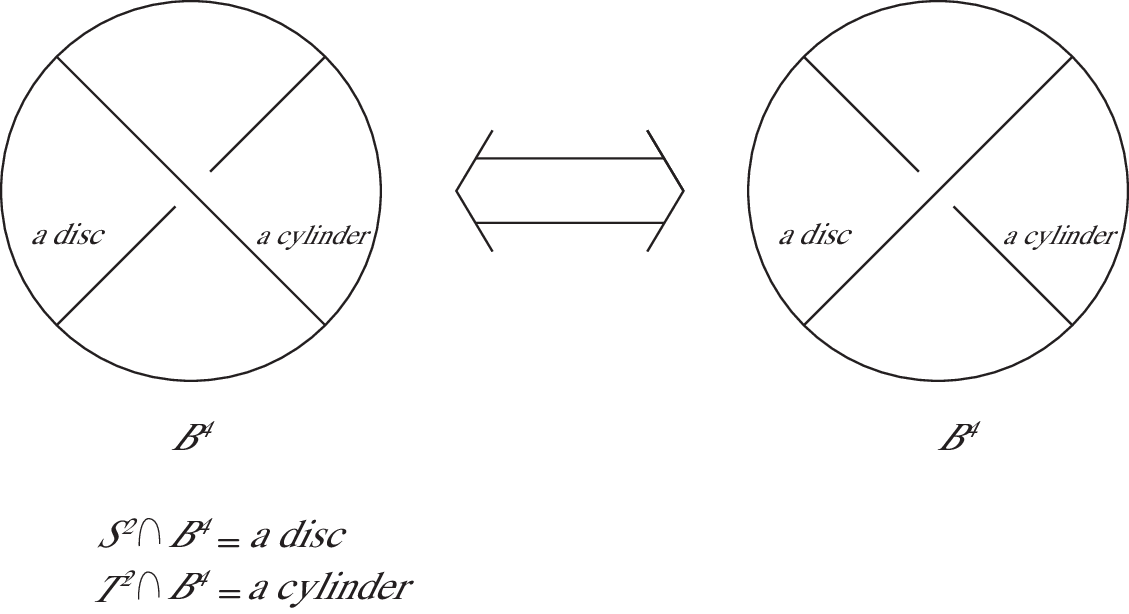}  
\bigbreak\caption{{\bf }\label{ribbonmoveau}}   
\bigbreak 
\end{figure}

\begin{figure}[H]
\bigbreak
\hskip18mm\includegraphics[width=10cm]{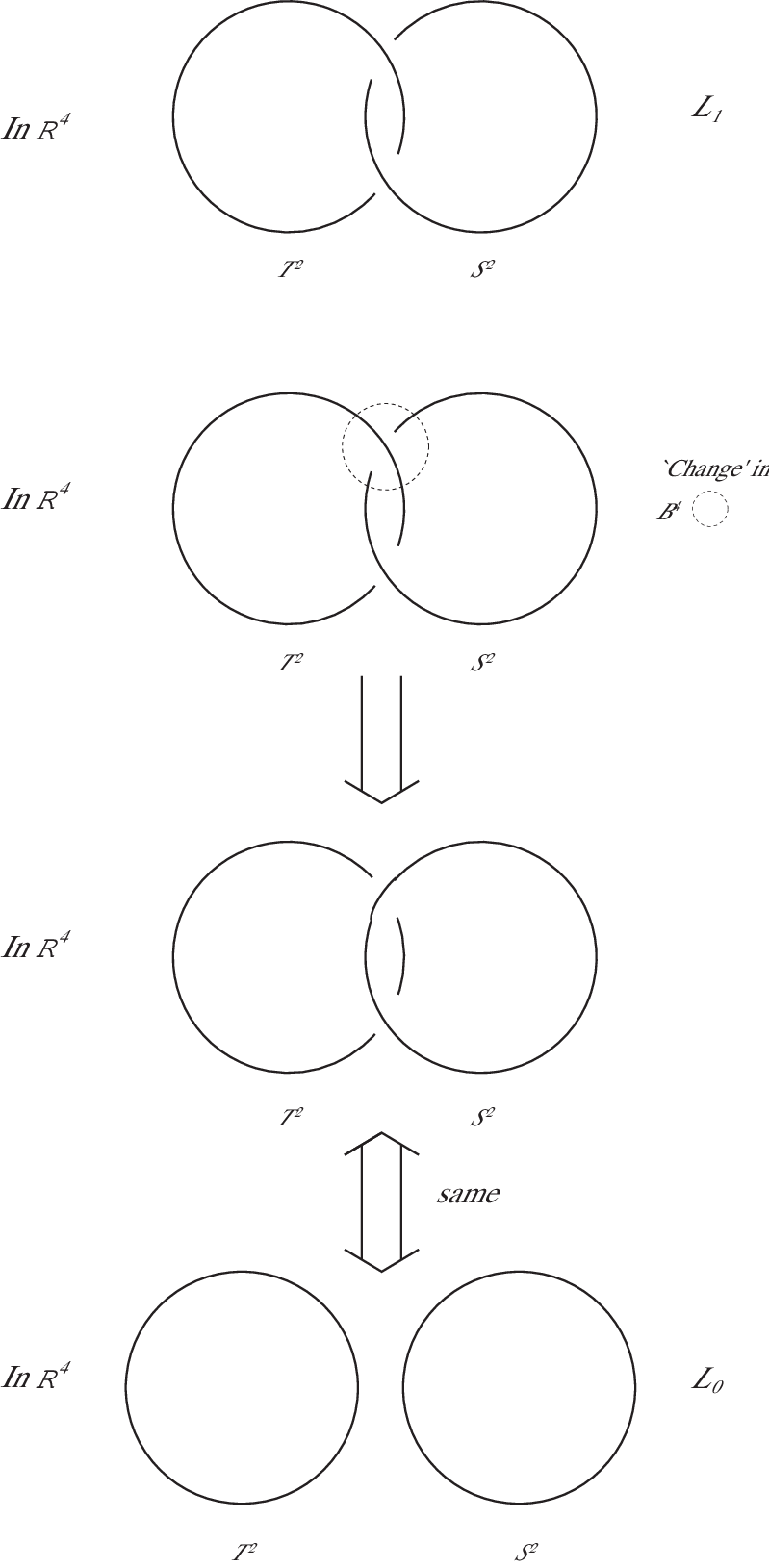} 
\bigbreak\caption{{\bf }\label{T2bau}}   
\bigbreak 
\end{figure}

\vskip1cm
By an analogy to the previous case `$S^1$ and $S^2$ in $\R^4$',  
we can guess 
$L_1$ is changed into $L_0$ by the operation 
as drawn in Figure \ref{ribbonmoveau}. 
Note that $B^n$ is the $n$-dimensional ball ($n\in\N$). 
Indeed it is true. 
The procedure is drawn in  Figure \ref{T2bau}.

\bigbreak

\bigbreak

This procedure is also called the {\it ribbon-move}. 
We say that 
$L_1$ (resp. $L_0$) is changed into $L_0$ (resp. $L_1$) by a ribbon-move. 
We will write a more rigorous figure of  ribbon-moves 
in Figures \ref{1.1au}  and \ref{1.2au}.

\np

\np

We draw a little more explicit figure of the ribbon-move 
in Figures 
 \ref{1.1au}  
and \ref{1.2au}.

\begin{figure}[H]
\bigbreak
\input 1.1.tex
\bigbreak 
\bigbreak 
\bigbreak\caption{{\bf }\label{1.1au}}   
\end{figure}

\np

\begin{figure}[H]
\bigbreak
\input 1.2.tex
\bigbreak 
\bigbreak 
\bigbreak\caption{{\bf }\label{1.2au}}   
\end{figure}

\noindent  {\bf To  the advanced readers.}    
We will write the definition of  ribbon-moves rigorously 
in 
Definition \ref{ribodef} 
by 
using Figures 
 \ref{1.1au}  
and \ref{1.2au}.

\np 
\subsubsection{Ribbon moves of 2-knots}

We have shown ribbon-moves on `$S^2$ and $T^2$ in $\R^4$'.  

Next we introduce ribbon-moves on 2-dimensional knots in $\R^4$. 

It can be defined as in the same manner in the case of `$S^2$ and $T^2$ in $\R^4$'.  

\bigbreak 
In  $B^4$ in $\R^4$ as shown in Figure \ref{ribbonmoveau}, 
the intersection of the 2-knot and $B^4$ is a set of a disc and a cylinder 
such that there is no common point between the disc and the cylinder. 

\bigbreak 

We will show 
a sequence of a finite number of ribbon-moves  
which makes the trivial 2-knot into a nontrivial 2-knots in Figure \ref{2reiau}. 

\begin{figure}[H]
\bigbreak
\includegraphics[width=120mm]{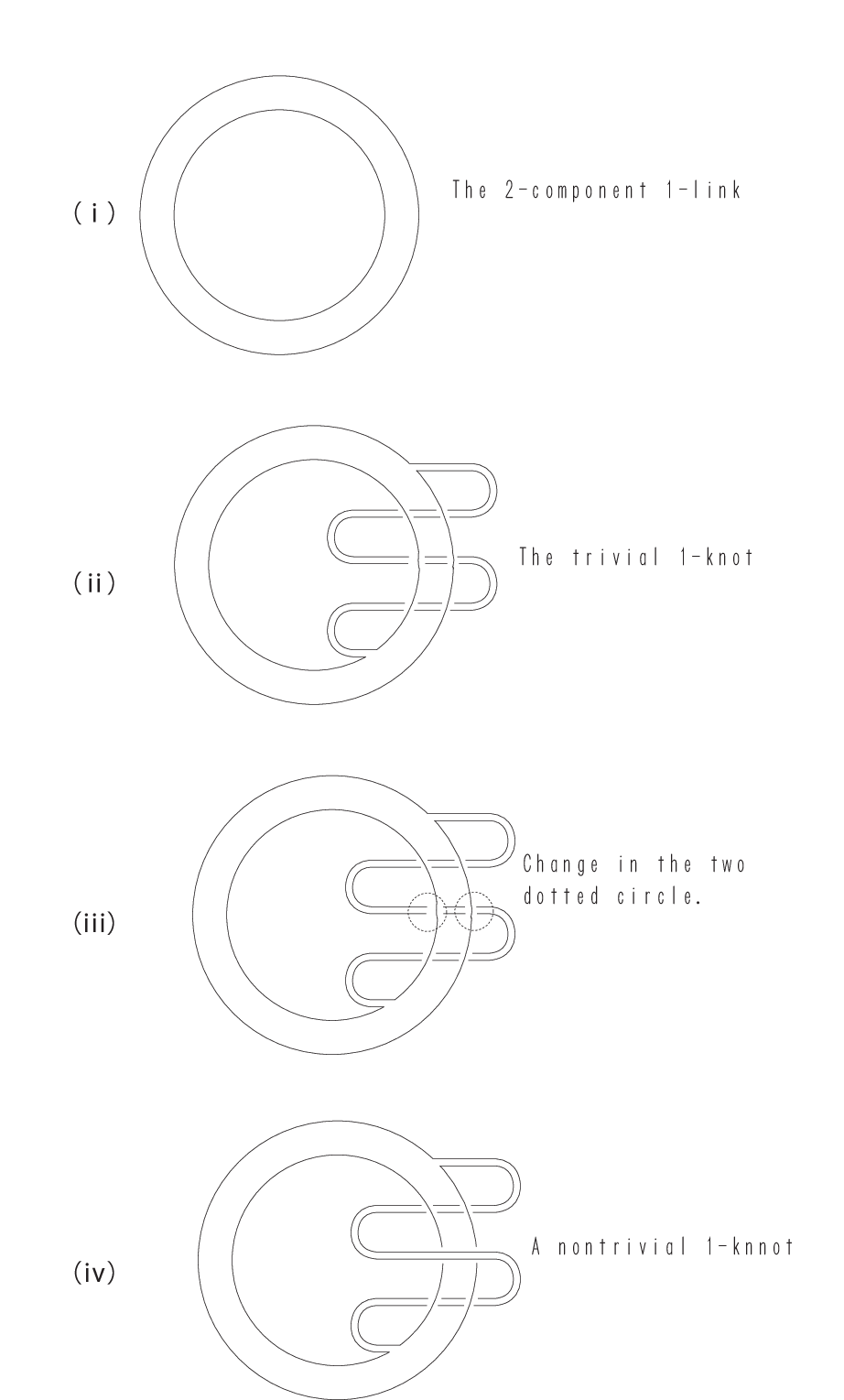}
\bigbreak\caption{{\bf 
(ii)-(iv) is a sequence of 1-knots}\label{1reiau}}   
\bigbreak 
\end{figure}

\begin{figure}[H]
\bigbreak
\includegraphics[width=120mm]{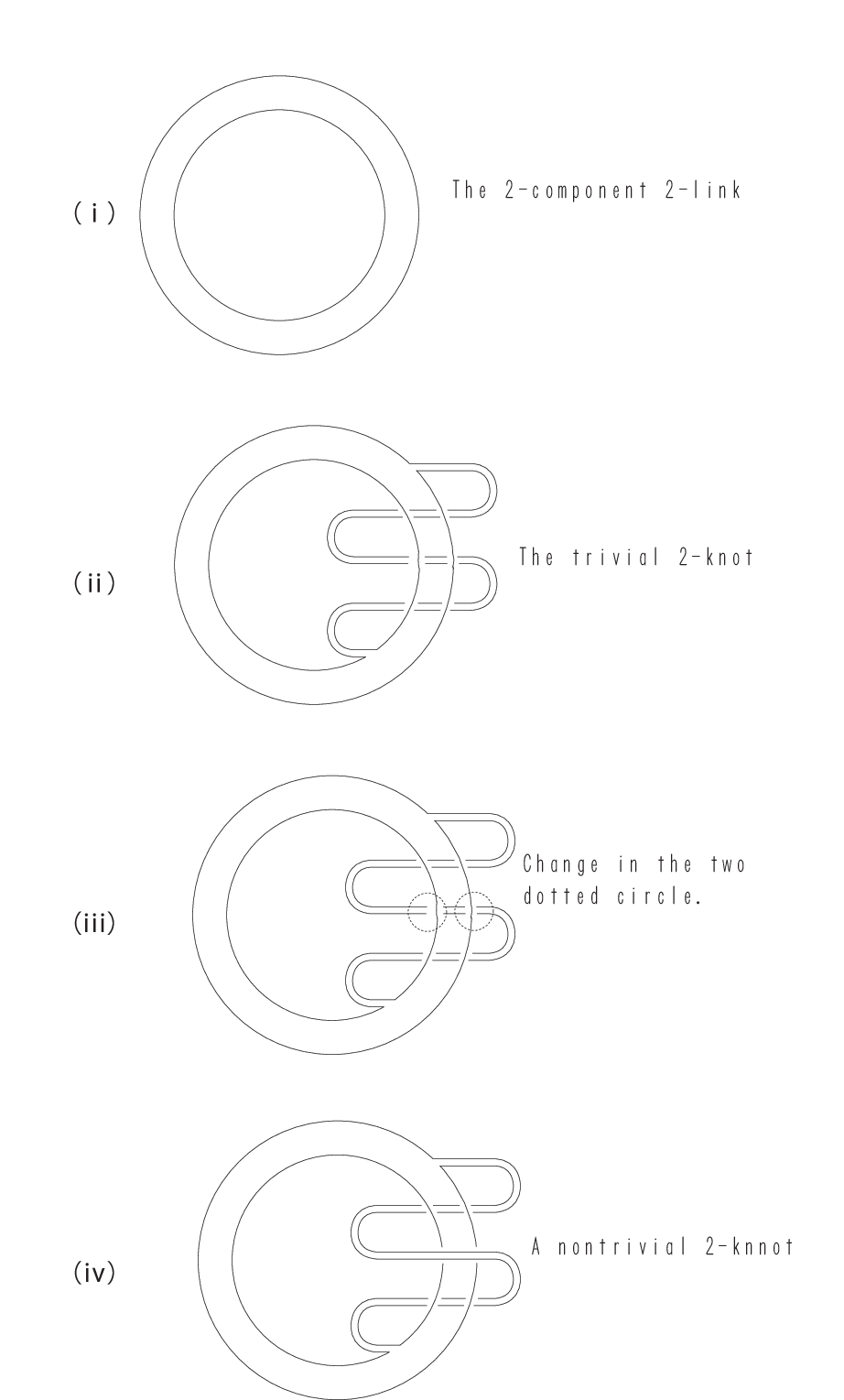}
\bigbreak\caption{{\bf 
We carry out the ribbon-move in (iii).}\label{2reiau}}   
\bigbreak 
\end{figure}

Before that,  
we will show 
a sequence of 
1-knots in Figure \ref{1reiau}.

The sequence of 2-knots in Figure \ref{2reiau}  
is a one-dimensional-higher analogue of 
the sequence of 1-knots in Figure \ref{1reiau}.

\np

\np

\np   

\noindent  {\bf To  the advanced readers.}     
We write the definition of the ribbon-move 
on 2-knots (resp. 2-links, $S^2\amalg T^2$ in $S^4$)
by using Figures \ref{1.1au} and \ref{1.2au}. 
We use some mathematical terms.    
Please check the literature. 
\vskip1cm

\begin{defn}\label{ribodef}
Let $K_1$ and $K_2$ be 2-links or 2-dimensional closed submanifolds in $\R^4$. 
We say that $K_2$ is obtained from $K_1$ by one {\it ribbon-move } 
if there is a 4-ball $B$ embedded in $\R^4$ with the following properties.  

\vskip3mm
(1) 
$K_1$ coincides with $K_2$ in $\R^4-(\text{the interior of }B)$. 
This identity map from $K_1-(\text{the interior of }B)$ to 
$K_2-(\text{the interior of }B)$  
is orientation preserving.

(2) 
$B\cap K_1$ is drawn as shown in Figure \ref{1.1au}.    
$B\cap K_2$ is drawn as shown in Figure \ref{1.2au}.    
\vskip3mm

\vskip1cm

%%%%%%%%%%%%%%%%%%

We regard $B$ as 
(a closed 2-disc)$\times[0,1]\times\{t| -1\leqq t\leqq1\}$.
We suppose $B_t=$(a closed 2-disc)$\times[0,1]\times\{t \}$.  
Then $B=\cup B_t$. 
In Figures \ref{1.1au} and \ref{1.2au}, 
we draw $B_{-0.5}, B_{0}, B_{0.5}$ in $B$. 
We draw $K_1$ and $K_2$ by the bold line. 
The fine line denotes $\partial B_t$. 
  
$B\cap K_1$ (resp. $B\cap K_2$) is diffeomorphic to 
$D^2\amalg (S^1\times [0,1])$, 
where $\amalg$ denotes the disjoint union. 

\vskip1cm

$B\cap K_1$ has the following properties:  
$B_t\cap K_1$ is empty for $-1\leqq t<0$ and $0.5<t\leqq1$.
$B_0\cap K_1$ is diffeomorphic to 
$D^2\amalg(S^1\times [0,0.3])\amalg(S^1\times [0.7,1])$. 
$B_{0.5}\cap K_1$ is diffeomorphic to $(S^1\times [0.3,0.7])$. 
$B_t\cap K_1$ is diffeomorphic to $S^1\amalg S^1$ for $0<t<0.5$. 
(Here, we draw $S^1\times [0,1]$ to have the corner 
in $B_0$ and in $B_{0.5}$. 
Strictly to say, $B\cap K_1$ in $B$ is a smooth embedding 
that is obtained by making the corner smooth.)

\vskip1cm

$B\cap K_2$ has the following properties:  
$B_t\cap  K_2$ is empty for $-1\leqq t<-0.5$ and $0<t\leqq1$.
$B_0\cap K_2$ is diffeomorphic to 
$D^2\amalg(S^1\times [0, 0.3])\amalg(S^1\times [0.7, 1])$. 
$B_{-0.5}\cap  K_2$ is diffeomorphic to $(S^1\times [0.3, 0.7])$. 
$B_t\cap  K_2$ is diffeomorphic to $S^1\amalg S^1$ for $-0.5<t<0$. 
(Here, we draw $S^1\times [0,1]$ to have the corner 
in $B_0$ and in $B_{-0.5}$. 
Strictly to say, $B\cap K_1$ in $B$ is a smooth embedding 
that is obtained by making the corner smooth.)

\np

In Figure \ref{1.1au} (resp. \ref{1.2au})   
there are an oriented cylinder $S^1\times [0,1]$ 
and an oriented disc $D^2$ as we stated above. 
We do not make any assumption about 
the orientation of the cylinder and the disc. 
The orientation of $B\cap K_1$ (resp. $B\cap K_2$ ) 
coincides with that of the cylinder and that of the disc. 

\vskip1cm

If $K_1$ is obtained from $K_2$ by one ribbon-move,  
then we also say that $K_2$ is obtained from $K_1$ by one {\it ribbon-move}.

\vskip1cm

%\vskip3mm
%\begin{defn}\label{ribbonmove2}    
Two 2-links $K_1$ and $K_2$ are said to be {\it ribbon-move-equivalent} 
if there are 2-links 
$K_1=\bar{K}_1, \bar{K}_2,...,\bar{K}_{r-1},\bar{K}_r=K_2$  
 ($r,p$ are natural numbers and $p\geq2$)such that 
$\bar{K}_i$ is obtained from $\bar{K}_{i-1}$ $(1< i\leqq r)$ by one ribbon-move. 
\end{defn}

\vskip1cm
\noindent  {\bf To  the advanced readers.}     
The ribbon-moves on $T^2\amalg S^2$ in $S^4$ can change 
the alinking number of the  $T^2\amalg S^2$ in $S^4$. 
See 
\cite{Ogasa04, OgasaZ, OgasaOCSL, Sato}  
for the alinking number.

\np

It is known that we have the following: 

%One time of crossing-change of makes the Hopf link into the 2-component 1-dimensional link. 
 
%The Hopf link is made into the 2-component 1-dimensional trivial link by one time of crossing-change. 

%A sequence of finite number of crossing-changes makes any $m$-component 1-dimensional link 
%into the trivial $m$-component 1-dimensional link, where $m$ is any natural number. 

Any $m$-component 1-dimensional link  is made into 
the trivial $m$-component 1-dimensional link 
by a sequence of a finite number of crossing-changes, where $m$ is any natural number.

In particular, 
any 1-dimensional knot  is made into 
the trivial 1-dimensional knot  
by a sequence of a finite number of crossing-changes.

\vskip1cm
By the way, we have the following.

Let  $L_0$ and $L_1$ be the sets of $T^2$ and $S^2$ in $\R^4$ in Figure \ref{T2bau} of  \S\ref{ribbTS}.

$L_0$ is made into $L_1$ by one ribbon-move.

\vskip1cm
We can make countably infinitely many nontrivial 2-knots in $\R^4$  
from the trivial 2-knot in $\R^4$ 
by a sequence of a finite number of ribbon-moves.  
An example is drawn in Figure \ref{2reiau}. 

It is well known that the spun knot of any 1-knot is ribbon-move-equivalent to the trivial 2-knot. 
Recall that a 1-knot is in $\R^3$ and that 
the spun knot of a 1-knot is a 2-knot in $\R^4$.

\vskip1cm
Well, it is natural to ask the following question. 

\bigbreak\noindent
{\Large{\begin{que}\label{kamakura}
Is any 2-dimensional knot in $\R^4$  ribbon-move-equivalent 
to the trivial 2-dimensional knot in $\R^4$?
\end{que}}}  
\bigbreak

\noindent
Note:  If $K$ is obtained from $K'$ 
by a sequence of a finite number of ribbon-moves, then we say that 
$K$ is {\it ribbon-move-equivalent} to $K'$.

\np

In his paper \cite{Ogasa04} 
the author proved that the answer to 
Question \ref{kamakura} is negative. 
He proved that 
the 5-twist spun knot of the trefoil knot is not ribbon-move-equivalent to the trivial 2-knot. 
Recall that the 5-twist spun knot of the trefoil knot is a 2-dimensional knot in $\R^4$. 
Furthermore he proved that 
there are countably infinitely many nontrivial 2-knots 
which are not ribbon-move-equivalent to the trivial knot. \\

\cite{Ogasa07, OgasaT3, OgasaZ, OgasaOCSL, Ogasanikai}   are sequels of the above paper. \\

The idea of the proof:  
Ribbon moves on 2-knots preserve the $\mu$-invariant of 2-knots,  
the $\Q/\Z$-valued $\widetilde\eta$-invariants of 2-knots, 
Farber-Levine tosion linking pairing. 
Furthermore they preserve 
some partial information of the homology groups 
with intersection products, 
a kind of intersections of submanifolds,  
and torsion linking pairings 
of the infinite cyclic covering spaces 
of the complements of 2-knots.  
%See the author's papers \cite{Ogasa04, Ogasa07, OgasaT3}   for detail. 

\vskip10mm
See the author's paper \cite{Ogasanikai}. 
It has been known that 
there is a 1-knot 
which is not changed into the trivial 1-knot by one crossing-change 
but which is made into the trivial one by two. 
The author proved that 
there is a 2-knot 
which is not changed into the trivial 2-knot by one ribbon move  
but which is made into the trivial one by two.

\np

Let $K$ and $K'$ be arbitrary 2-knots. 
Then is $K$ ribbon-move-equivalent to $K'$?

\vskip1cm
That is, classify 2-knots by ribbon-move-equivalence. 

\vskip1cm
It is an open problem. 

\vskip1cm
Please solve this.

\np 
%\begin{sub}\label{passmoves}
%$(p,q)$-pass-moves on $n$-dimensional knots ($p+q=n+1$)\end{sub}
\subsection
{$(p,q)$-pass-moves on $n$-dimensional knots ($p+q=n+1$)}
\label{passmoves}

You can imagine $\R^4$ so you can see $\R^n (n\geqq5)$!
Let's go into high dimmensional space! 

\vskip10mm

We introduce another local move on high dimensional knots, 
the $(p,q)$-pass-move on $n$-dimensional knots ($p+q=n+1$ and $n\geqq1$).   
%after we introduce products of figures. 

\vskip10mm
Before that, 
we introduce a local move on 1-links, 
the pass-move on 1-knots. 
We generalize it into the high dimensional case 
and 
define 
the $(p,q)$-pass-move on $n$-dimensional knots.  

\np
See Figure  \ref{passau}   
 for an illustration of the pass-move on 1-links. 
Here, we consider not only 1-knot case but also 1-link case.  
%We often abbreviate (1,1)-pass-move to pass-move. 

The 1-link is oriented and is represented by the arrow 
in Figure \ref{passau}.    

Each of four arcs in the 3-ball 
may belong to different components of the 1-link. 

\begin{figure}[H]
\bigbreak
\input pass.tex
\bigbreak\caption{{\bf }\label{passau}}   
\bigbreak 
\end{figure}

First we change the trivial knot into a nontrivial knot 
by a pass-move shown in Figures \ref{F2au}-\ref{F4au}.  

We will generalize this way to high dimensional case. 

Note that we sometimes abbreviate 1-dimensional knot to 1-knot or knot.  

\np

Regard $S^1$ as the union of two 1-ball $B^1_u$ and $B^1_d$.  
Note that the 1-ball is the segment.   
  
Then $S^1\x S^1$ is regarded as 
the union of four parts, 
$B^1_u\x B^1_u$, 
$B^1_u\x B^1_d$, 
$B^1_d\x B^1_u$, and 
$B^1_d\x B^1_d$. 

Remove the interior of $B^1_u\x B^1_u$ from $S^1\x S^1$, call it $F$. 

\begin{figure}[H]
\bigbreak
\includegraphics[width=160mm]{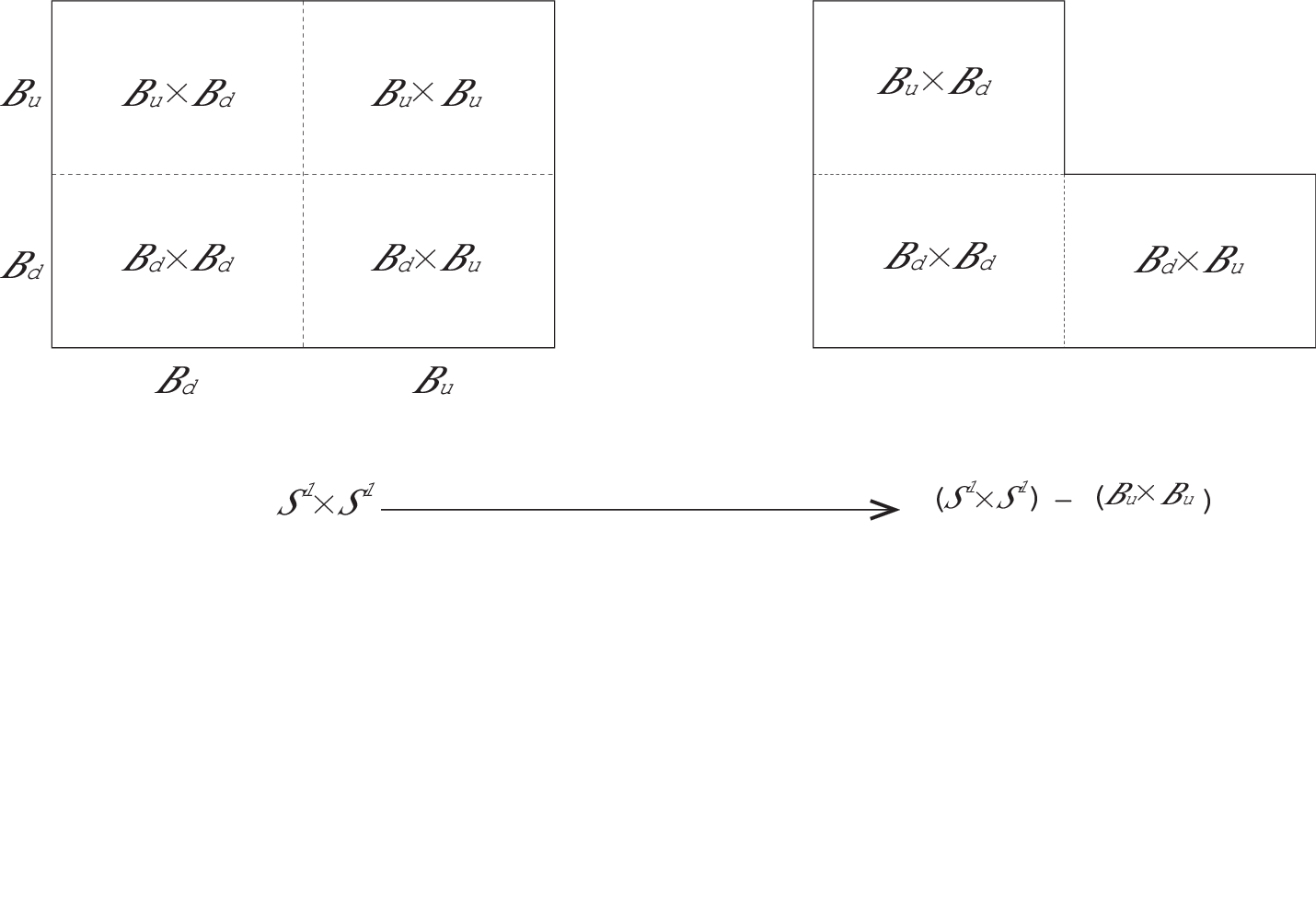}  
\vskip-50mm\bigbreak\caption{{\bf }\label{BuBdau}}   
\end{figure}

$F$ is drawn conceptually as shown in Figure \ref{BuBdau}. 
We abbreviate $B^\sharp_\star$ to  $B_\star$.

\vskip10mm

\np

\begin{figure}[H]
\bigbreak
\includegraphics[width=160mm]{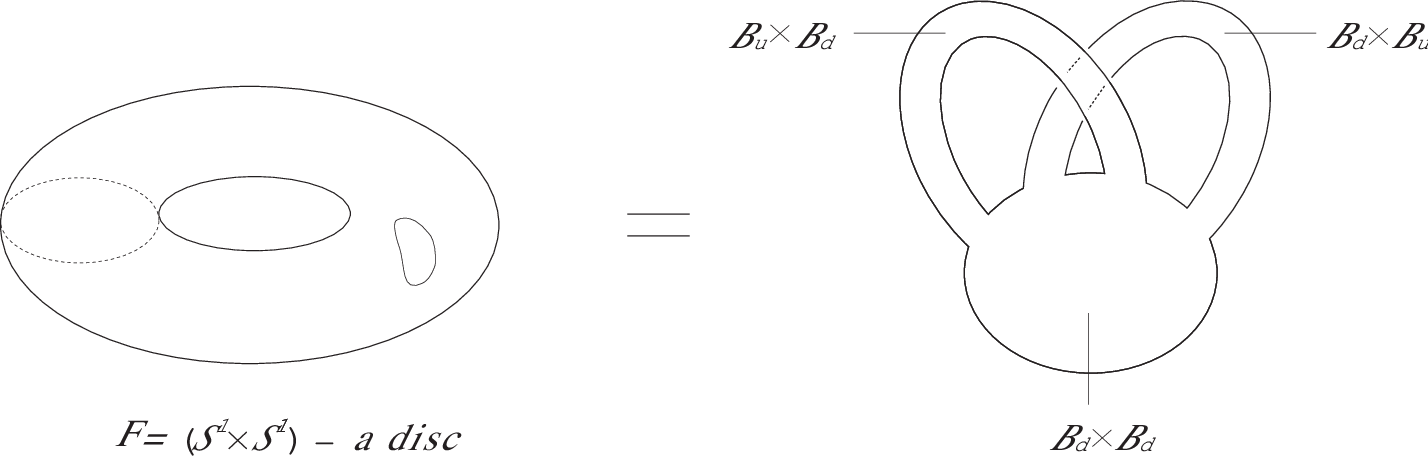}  
\bigbreak
\bigbreak
\bigbreak
\caption{{\bf }\label{Fau}}   
\bigbreak
\end{figure}

\begin{figure}[H]
\bigbreak
\includegraphics[width=70mm]{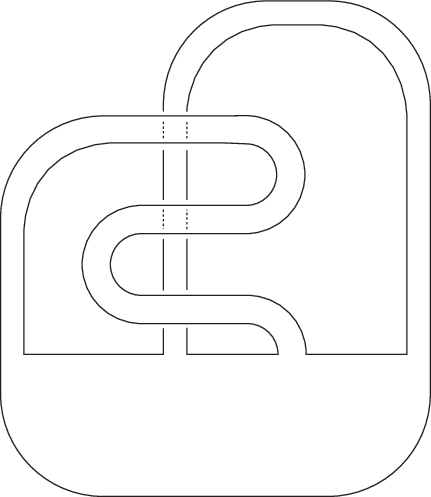} 
\bigbreak\caption{{\bf }\label{F2au}}   
\bigbreak 
\end{figure}

\begin{figure}[H]
\bigbreak
\includegraphics[width=70mm]{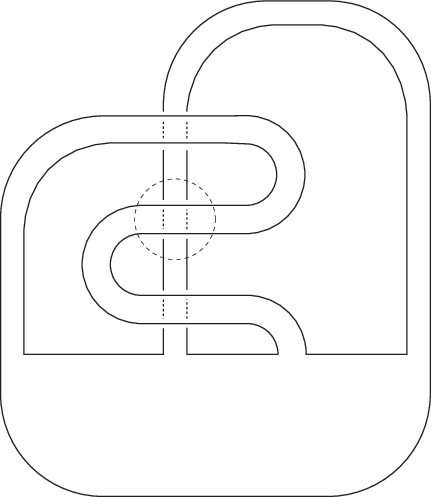}   
\bigbreak\caption{{\bf }\label{F3au}}   
\bigbreak 
\end{figure}

\begin{figure}[H]
\bigbreak
\includegraphics[width=70mm]{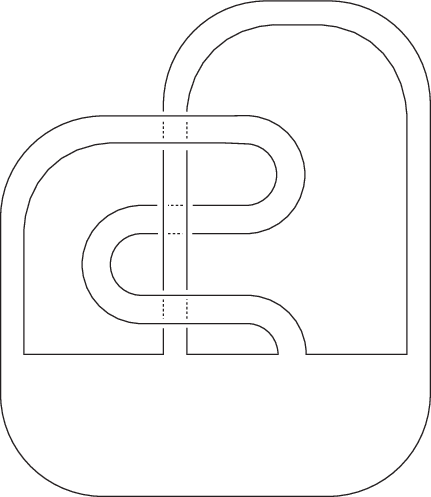}   
\bigbreak\caption{{\bf }\label{F4au}}   
\bigbreak 
\end{figure}

$F$ is drawn explicitly as shown in Figure \ref{Fau}. 
Note that we can bend the corner of $B^1_u\x B^1_u$ and change it into the 2-dimensional ball. 
Hence the boundary of $F$ is a single circle.  
\vskip10mm

\np

Take $F$ in $\R^3$ as 
drawn in Figure \ref{F2au}. 

\vskip1cm.

The boundary of $F$ in $\R^3$ is a 1-knot. 
Furthermore it is the trivial knot.

\np

Carry out a pass move on this knot 
in 
%\includegraphics[width=7mm]{B4}   %zz-dc
%\includegraphics[width=1cm]{dc}   %zz-dc
%
%
%\begin{figure}[H]
%\bigbreak
\includegraphics[width=7mm]{B4}   
%\bigbreak\caption{{\bf }\label{B4au}}   
%\bigbreak 
%\end{figure}
in Figure \ref{F3au}.

\np
%We think that you feel that the resulting knot 
Probably you will feel the resulting knot 
in Figure \ref{F4au}

\noindent 
is a nontrivial knot. 
Yes.  It is true. 
%A way to prove it is to use Seifert matrices and the Alexander polynomial. 
We can prove it by using Seifert matrices and the Alexander polynomial.

\np

It is known that there are countably infinitely many nontrivial knots 
which are pass-move-equivalent to the trivial knot 
and 
which bound  $F$ 
(resp. which do not bound $F$).  
 
Any 1-knot is pass-move-equivalent to the trivial knot or the trefoil knot. 
Furthermore the trivial knot is not pass-move-equivalent to the trefoil knot.

\vskip1cm
\noindent  {\bf To  the advanced readers.}    
We have the following. 
See \cite{Kauffman1987}.

If the Arf invariant of a 1-knot $K$ is zero (resp. one), then $K$ is pass-move-equivalent to 
the trivial knot (resp. the trefoil knot).

\smallbreak
Furthermore we have the following. 
Let $L_1$ and $L_2$ be 1-links. 
Then  $L_1$ and $L_2$ are pass-move-equivalent if and only if 
$L_1$ and $L_2$   satisfy one of the  following conditions (1) and (2).

\noindent 
(1)
 Both $L_1$ and $L_2$ are proper links, and 
$$\mathrm{Arf}(L_1) = \mathrm{Arf}(L_2).$$

\noindent 
(2) 
   Neither $L_1$ nor $L_2$ is a proper link, and 

$$\mathrm{lk}(K_{1j}, L_1-K_{1j})\equiv   
  \mathrm{lk}(K_{2j}, L_2-K_{2j}) 
  \hskip2mm\mathrm{mod 2}\hskip2mm for\hskip2mm  all\hskip2mm j. $$ 

\bigbreak
Here,  $\mathrm{lk}(K_{aj}, L_a-K_{aj})$ is the sum of all linking number 
lk$(K_{aj}, K_{ai})$, where $i\neq j$. 
%We introduced the linking number in \S\ref{linking}. 
Here, lk is the linking number. 
See \cite{Kauffman1987}  
for 
the linking number.

\np

Let's start the high dimensional case.

\vskip1cm
Let $A$ be in $\R^n$. 
Let $B$ be in $\R^m$.

Define $A\x B$ to be  

\noindent
$\{(x_1,...,x_n,y_1,...,y_m)\vert (x_1,...,x_n) \text{is in $A$.}  (y_1,...,y_m) \text{is in $B$}\}$
in $\R^{n+m}$. 

\begin{figure}[H]
\bigbreak
\includegraphics[width=133mm]{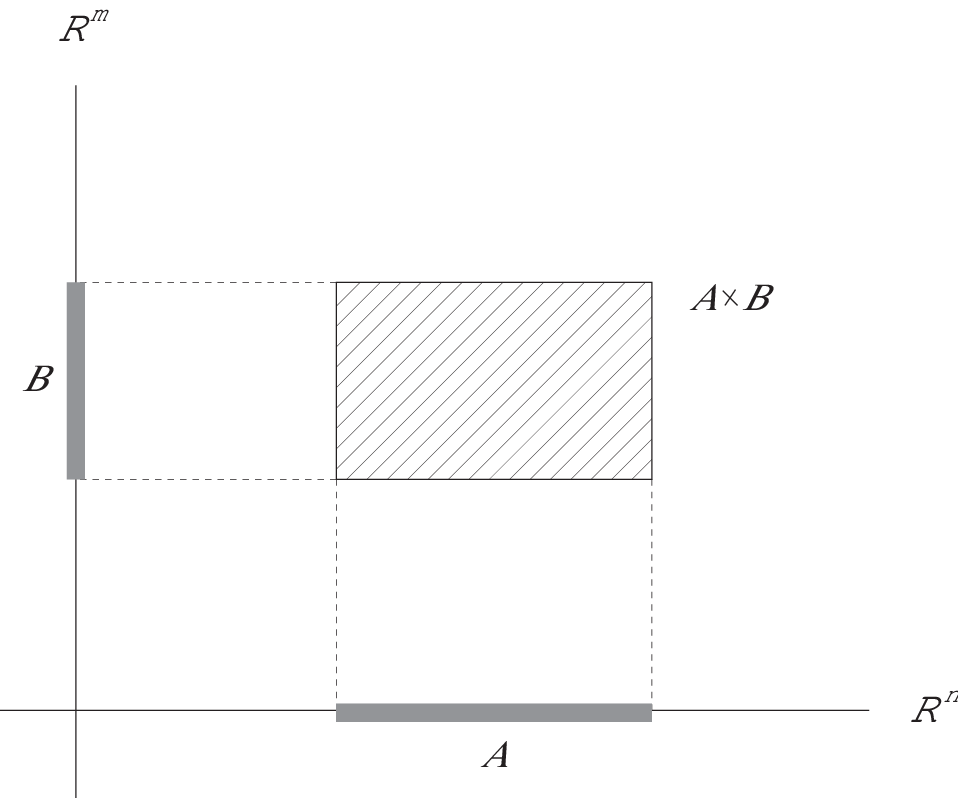}   
\bigbreak\caption{{\bf $A\x B$}\label{xau}}   
\bigbreak 
\end{figure}

%It is conceptually drawn as follows.
See Figure \ref{xau}.

\np

$S^1\x S^1$ is an example. 

$S^1$ is represented by 

$\{ (x_1,x_2)  \vert  (x_1)^2+(x_2)^2=1 \}$ in $\R^2$

$S^1\x S^1$ is represented by 

$\{ (x_1,x_2,y_1,y_2)  \vert  (x_1)^2+(x_2)^2=1,\quad  (y_1)^2+(y_2)^2=1 \}$ in $\R^4$. 

\noindent
Note:  
In this way, $S^1\x S^1$ is defined in  $\R^4$ 
but it can be embedded in $\R^3$. 
%We have this situation in some other cases. 
$S^1\x S^1$ is the torus.

%\vskip1cm  $S^1\x S^{n-1}$ ($n\geqq3$) is defined  in the similar fashion.  

\vskip1cm
\noindent  {\bf To  the advanced readers.}    
We want to talk about product sets and product manifolds. 
In fact they can be defined without using $\R^l$ which includes figures.

\np

We define the $(p,q)$-pass-move on $n$-knots ($p+q=n+1$). 

Take an $(n+2)$-ball $B^{n+2}$ trivially embedded in $\R^{n+2}$.

Regard $B^{n+2}$ as $B^1\x B^p\x B^q$. 
Note that  $B^l$ is 

\noindent 
$\{ (x_1,...,x_l) \vert (x_1)^2+...+(x_l)^2\leqq1  \}$ for any natural number $l$.  
In particular, that $B^1$ is $[-1,1]=\{ x \vert -1\leqq x\leqq1  \}$. 

Take a smaller $p$-ball $D^p$ in $B^p$ and 
suppose that  $D^p$ is 

\noindent
$\{ (x_1,...,x_p) \vert (x_1)^2+...+(x_p)^2\leqq\frac{1}{4}  \}$. 

Take a smaller $q$-ball $D^q$ in $B^q$ as well.

\begin{figure}[H]
\bigbreak
\hskip-43mm\includegraphics[width=190mm]{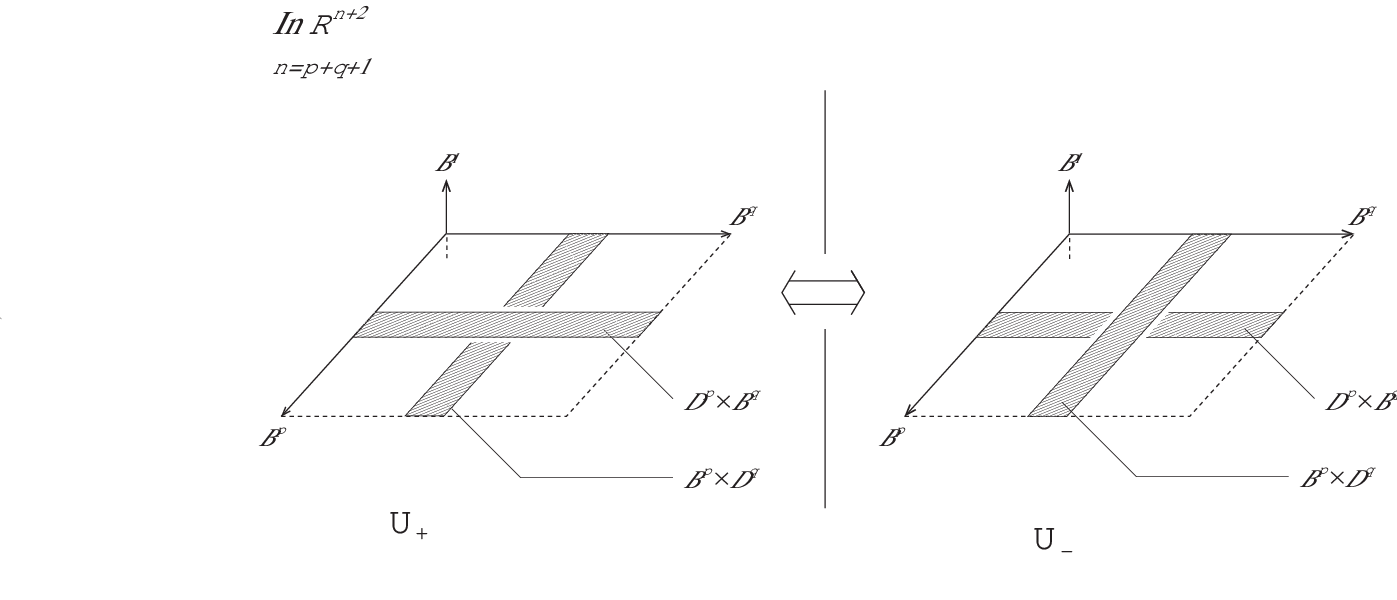}  
\bigbreak\caption{{\bf }\label{bpqau}}   
\end{figure}

\bigbreak 
Take  
 $D^p\x B^q$ and $B^p\x D^q$ 
in  the left (resp. right) figure of $B^{n+2}$  
as shown in Figure \ref{bpqau}, and call 
the union of  $D^p\x B^q$ and $B^p\x D^q$,  
$U_+$ (resp. $U_-$).     
Imagine!  
Note that  $D^p\x B^q$ and  $B^p\x D^q$ do not touch each other. 
Imagine!

\bigbreak
\bigbreak

\np

\bigbreak
$B^p\x D^q$ of $U_+$ (resp. $U_-$) is embedded in $\{0\}\x B^p\x B^q$ trivially.

$D^p\x B^q$ of  $U_+$ is embedded in $\{x\geqq0\}\x B^p\x B^q$ trivially. 

$D^p\x B^q$ of  $U_-$ is embedded in $\{x\leqq0\}\x B^p\x B^q$ trivially.

\bigbreak
\bigbreak
\noindent
{\bf Note to the advanced readers:} 
`Embedding trivially' is mathematically defined. 
However, as we told before, 
if you currently associate the word `embedding trivially' with 
how it is used in daily life, please continue to do so. 
%but 
%they can regard it as 
%what they associate `the word in the daily life' with for now.

\np

Let $S^{p-1}$ be the boundary of  $D^p$. 

Hence $S^{p-1}\x B^q$ is included in the boundary of  $D^p\x B^q$. 

 $B^p\x S^{q-1}$ is included in the boundary of $B^p\x D^q$ as well.  

\begin{figure}[H]
\bigbreak

\includegraphics[width=150mm]{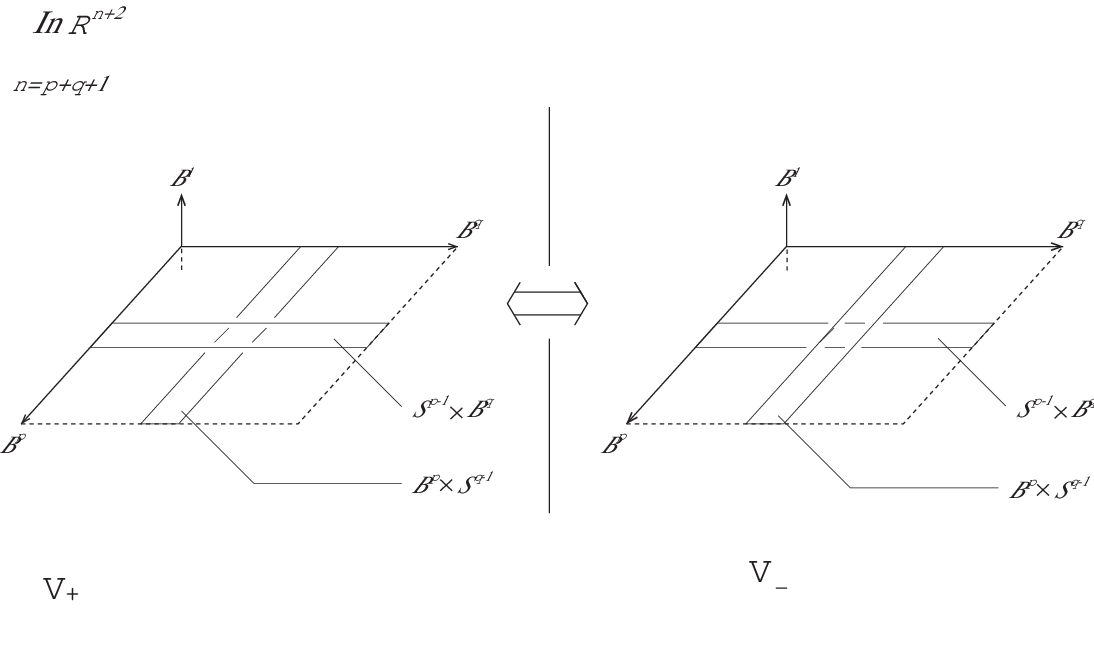} 
\bigbreak\caption{{\bf }\label{pqau}}   
\bigbreak 
\end{figure}

Take $S^{p-1}\x B^q$ and  $B^p\x S^{q-1}$ 
in the left (resp. right) figure of $B^{n+2}$ 
as shown in Figure \ref{pqau}, 
and call the union of $S^{p-1}\x B^q$ and  $B^p\x S^{q-1}$, 
$V_+$ (resp. $V_-$).    
Imagine again!

\np

Let $K_+$ and $K_-$ be $n$-dimensional knots in $\R^{n+2}$. 
Suppose that 
$K_+$ and $K_-$ differ only in the above $B^{n+2}$,  
and that 
the intersection of $K_+$ (resp. $K_-$) and $B^{n+2}$ is 
$V_+$ (resp. $V_-$).    
Then we say that  $K_-$ (resp. $K_+$) is obtained from 
$K_+$ (resp. $K_-$) by one {\it $(p,q)$-pass-move}. 

\vskip1cm
In order to imagine the high dimensional case, 
it is very often useful to consider a low dimensional case.  
Consider the case where $p=q=1$ and hence $n=3$. 
You can easily understand 
that 
this case, that is, 
the (1,1)-pass-move on 1-knots, 
is 
the pass-move on 1-knots.

\vskip1cm
Let $p$ (resp. $q$) be an fixed natural number. Let $p+q=n+1$. 

Could you imagine that a nontrivial $n$-dimensional knot is obtained from 
the trivial $n$-dimensional knot by one $(p,q)$-pass-move 
(resp. a sequence of a finite number of $(p,q)$-pass-moves)? 

We explain it from here on. 
See Figures \ref{F2au>}-\ref{F4zu}.

\vskip1cm
\noindent  {\bf To  the advanced readers.}    
See  the author's paper  
\cite{Ogasa98n, Ogasa04, Ogasa09, Ogasa07, OgasaT3, OgasaZ, Ogasanikai} 
for detail.

\begin{figure}[H]
\bigbreak
\includegraphics[width=170mm]{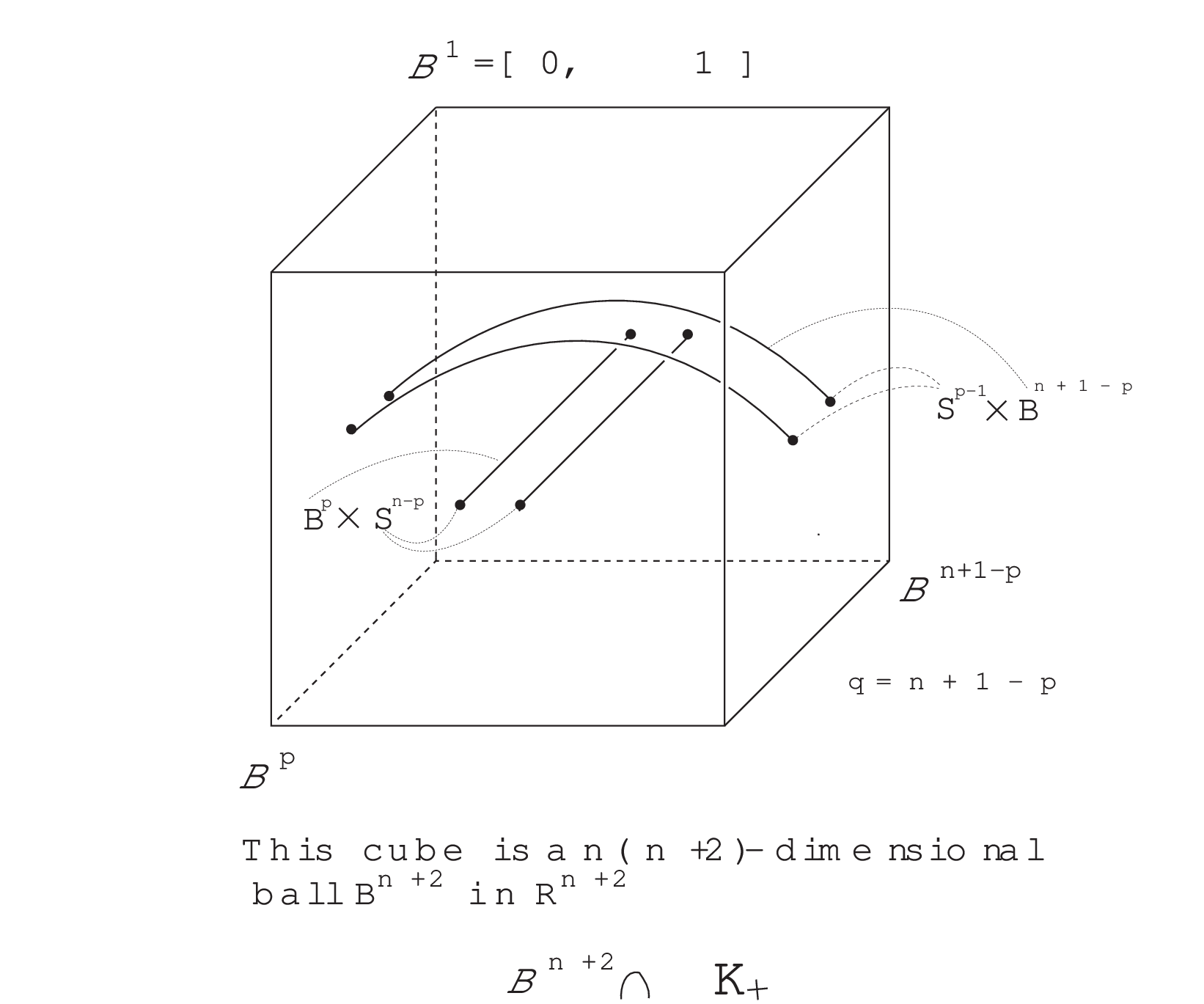}   
\bigbreak\caption{{\bf }\label{pqdaiau}}   
\bigbreak 
\end{figure}

\begin{figure}[H]
\bigbreak
\includegraphics[width=14cm]{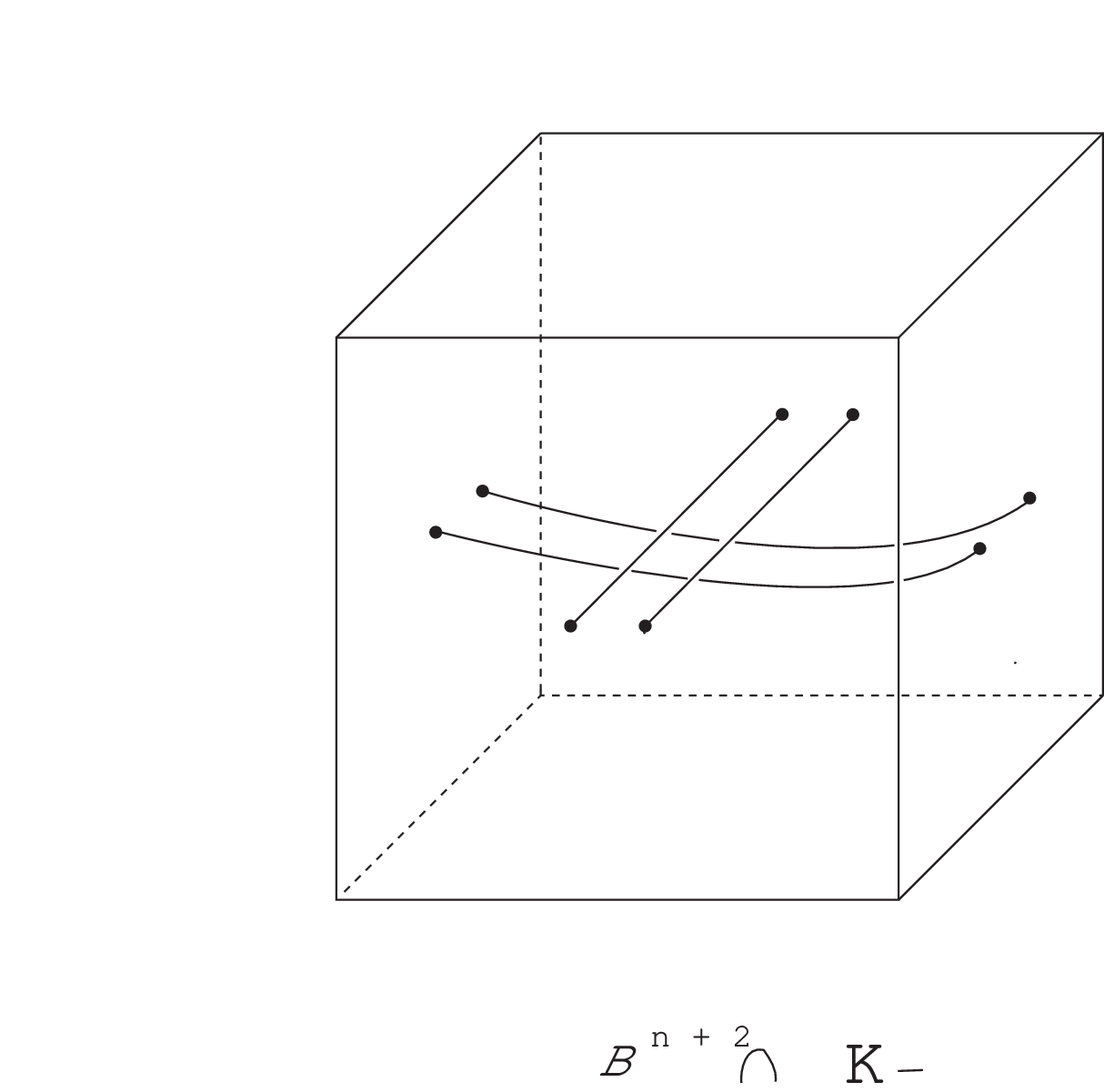}  
\vskip-20mm\caption{{\bf }\label{pqdai2au}}   
\end{figure}

\np  
A pair of 
Figures 
\ref{pqdaiau} and \ref{pqdai2au}
is
a little less conceptual figure of the $(p,q)$-pass-move.

\np

\np

Next consider the case where $p=1$ $q=2$, and $n=4$. 

\vskip1cm

One (1,2)-pass-move is two ribbon-moves. 

Can you imagine this?

\vskip1cm

Furthermore in fact if a 2-knot $K$ is obtained from a 2-knot $K'$ by one ribbon-move, 
then  
$K'$ is obtained from $K$ by one (1,2)-pass-move. 
The proof is written in the author's paper 
\cite{Ogasa04}.

\np

Could you imagine that a nontrivial $n$-dimensional knot is obtained from 
the trivial $n$-dimensional knot by one $(p,q)$-pass-move 
(resp. a sequence of a finite number of $(p,q)$-pass-moves)?

\vskip10mm
We can do it. We will show it from here on.

\np

Next we  change the trivial $n$-dimensional knot into a nontrivial $n$-dimensional knot 
by a $(p,q)$-pass-move ($p+q=n+1$).  Recall that  $n$-knots are in $\R^{n+2}$

%we show a $(p,q)$-pass-move on $n$-knot ($p+q=n+1$) in $\R^{n+2}$
% changes a trivial $n$-dimensional knotinto a nontrivial $n$-dimensional knot. 

Regard $S^p$ as the union of $p$-ball $B^p_u$ and $B^p_d$.  
Do $S^q$ as well.

Then $S^p\x S^q$ is regarded as 
the union of four parts, 
$B^p_u\x B^q_u$, 
$B^p_u\x B^q_d$, 
$B^p_d\x B^q_u$, and 
$B^p_d\x B^q_d$. 

Remove the interior of $B^p_u\x B^q_u$ from $S^p\x S^q$, call it $F$. 

\begin{figure}[H]
\bigbreak
\hskip-10mm\includegraphics[width=160mm]{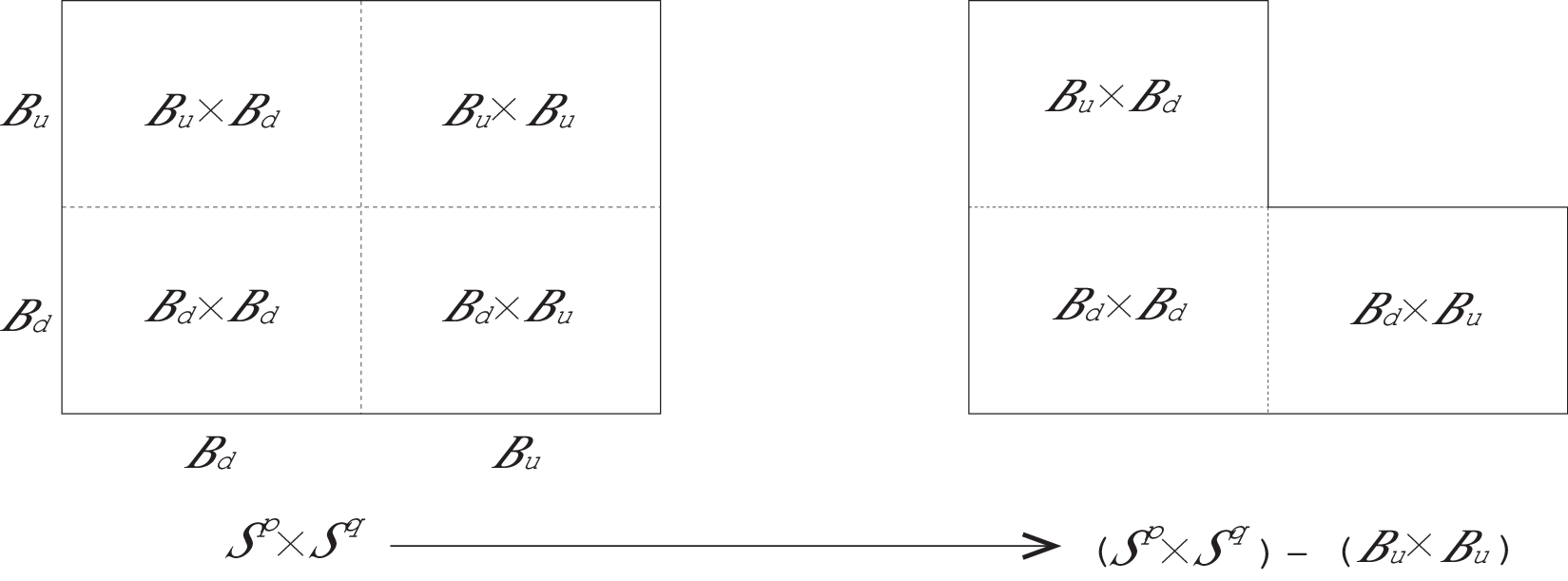} 
\bigbreak
\bigbreak
\caption{{\bf }\label{pqBuBdau}}   
\bigbreak 
\end{figure}

$F$ is drawn in Figure \ref{pqBuBdau}. 
We abbreviate $B^\sharp_\star$ to  $B_\star$.

\vskip10mm

\np

\begin{figure}[H]
\bigbreak
\hskip-10mm\includegraphics[width=160mm]{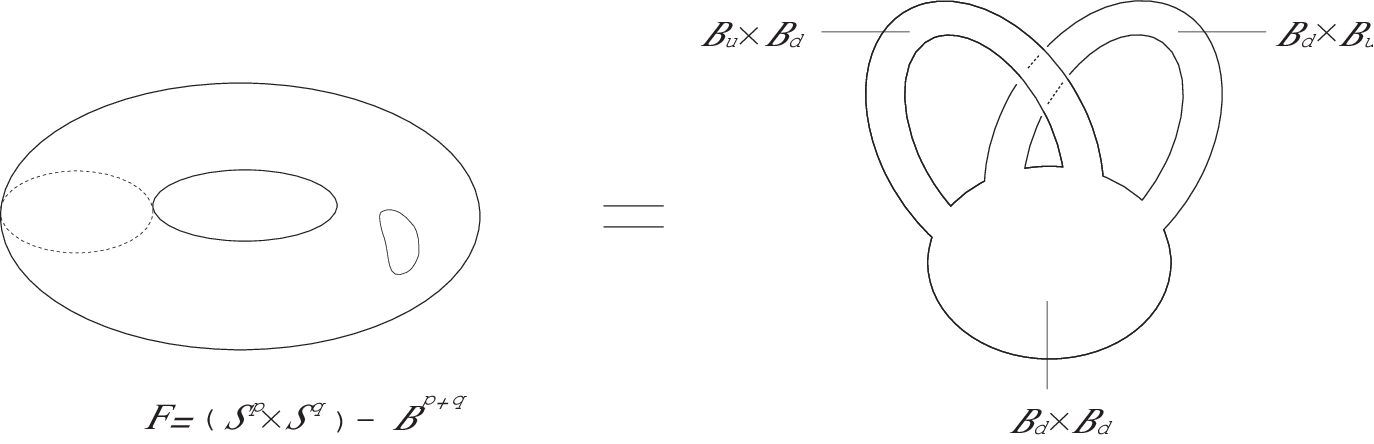} 
\bigbreak\caption{{\bf }\label{pqFau}}   
\bigbreak 
\end{figure}

\begin{figure}[H]
\bigbreak
\bigbreak
\bigbreak
\hskip20mm\includegraphics[width=80mm]{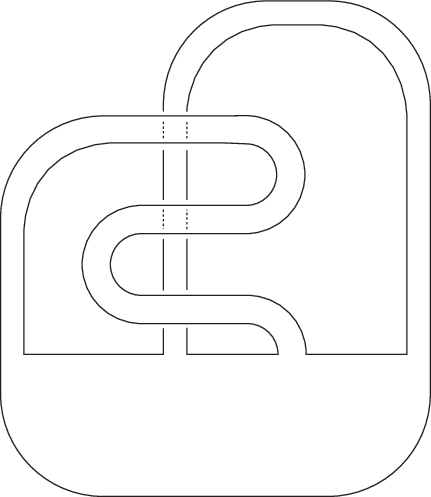}   
\bigbreak\caption{{\bf }\label{F2au>}}   
\bigbreak 
\bigbreak
\bigbreak
\end{figure}

\begin{figure}[H]
\bigbreak
\hskip20mm\includegraphics[width=80mm]{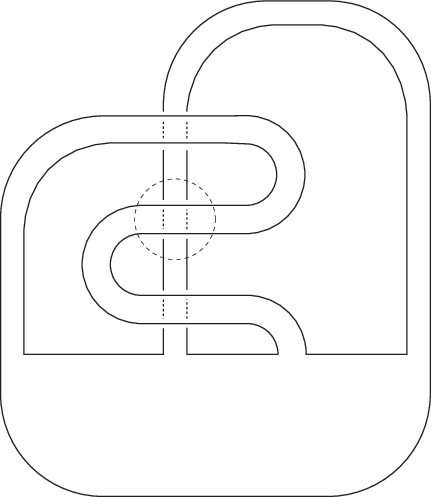}  
\bigbreak\caption{{\bf }\label{F3zu}}   
\bigbreak 
\end{figure}

\begin{figure}[H]
\bigbreak
\hskip20mm\includegraphics[width=80mm]{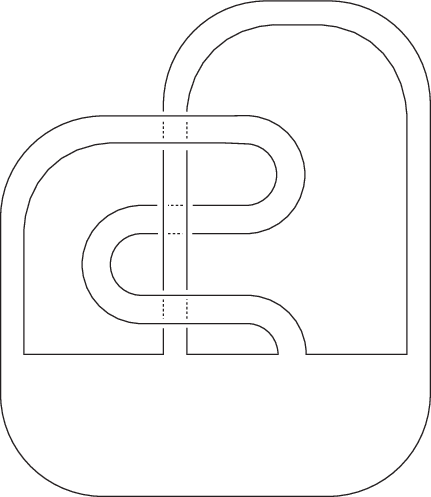}  
\bigbreak\caption{{\bf }\label{F4zu}}   
\bigbreak 
\end{figure}

$F$ is drawn in another way as shown in Figure \ref{pqFau}. 
Note that we can bend the corner of $B^p_u\x B^q_u$ 
and change it into the $(p+q)$-dimensional ball. 
Hence the boundary of $F$ is $S^n$.

\vskip10mm
\vskip10mm

\noindent  {\bf To  the advanced readers.}    
This figure is an example of handle decomposition of manifolds. 
See textbooks on differential topology and, especially, surgery theory for detail.

\np

Take $F$ in $\R^{n+2}$ as shown in Figure \ref{F2au>}.

The boundary of $F$ in  $\R^{n+2}$ is an $n$-knot.  
Furthermore it is the trivial $n$-knot. 
%It is drawn very conceptually. 

\np

Carry out a $(p,q)$-pass-move on this $n$-knot 
in 
%
%\begin{figure}[H]
%\bigbreak
\includegraphics[width=7mm]{B4}  
%\bigbreak\caption{{\bf }\label{}}   
%\bigbreak 
%\end{figure}
%
%
%\includegraphics[width=1cm]{dc}   %%zz-dc
in Figure \ref{F3zu}.

\np
It is known that the resulting $n$-knot 
in Figure \ref{F4zu}

\noindent 
is a nontrivial $n$-knot. 

We can prove this fact by using Seifert matrices and the Alexander polynomial. 
In the proof,  we use the fact that $S^p$ and $S^q$ can be `linked' in $\R^{p+q+1}$.    
%as implicitly explained in the last part of Chapter \ref{proof}. 

\np

\begin{figure}[H]
\bigbreak
\hskip20mm\includegraphics[width=80mm]{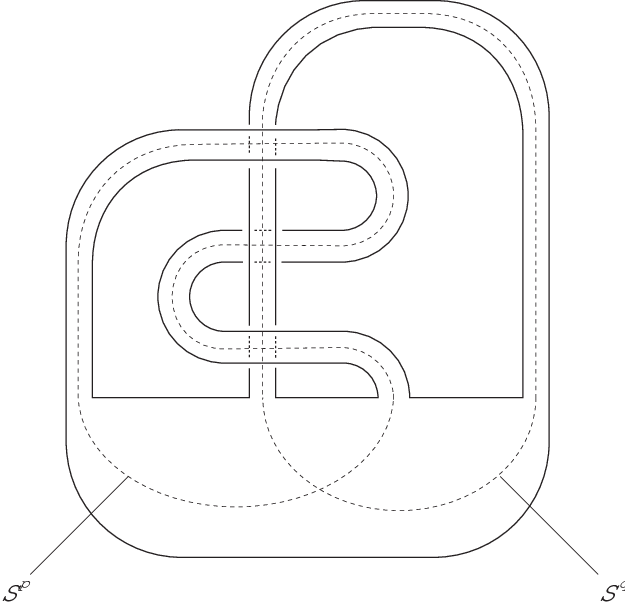}  
\bigbreak\caption{{\bf }\label{F5au}}   
\bigbreak 
\end{figure}

We use the following facts to prove it: 
$S^q$ and $S^p$ are included in $F$ as shown 
in Figure \ref{F5au}.    
$S^q$ and $S^p$ are `linked' in $S^{p+q+1}$.  
%(The problem associated with the last figur in \S\ref{S1Sn-1kotae} 
%is related to this situation.)  
Imagine!

\vskip10mm

It is known that there are countably infinitely many nontrivial $n$-dimensional knots 
which are pass-move-equivalent to the trivial $n$-knot 
and 
which bound  $F$ 
(resp. which do not bound $F$).

\vskip1cm
It is not solved completely that 
what kind of high dimensional knots are $(p,q)$-pass-move-equivalent to the trivial knot 
in both cases where $p$ and $q$ are fixed and where $p$ and $q$ run. 
The $p=q$ case is solved by the author. It is proved by using  
\cite[Theorem 4.1]{Ogasa98n} 
and 
\cite[Proposition 12.1]{KauffmanOgasa}. 
It is characterized as follows. Let $K$ be a $(2p+1)$ knot.  
$K$ is $(p+1,p+1)$-pass-move equivalent to the unknot 
if and only if $K$ is a simple knot and  
the Arf invariant (resp. the signature) of $K$ is zero in the case where 
$p$ is even (resp. odd).

\vskip1cm
It is also not solved to classify $n$-knots by $(p,q)$-pass-move-equivalence.

\np

%We found some examples of nontrivial $n$-knots which are not $(p,q)$-pass-move-equivalent to the trivial knot. 

We found some relations between some invariants of high dimensional knots 
and $(p,q)$-pass-moves. 

\bigbreak
The author  found an application of $(p,q)$-pass-moves to the intersection of submanifolds, which is a kind of figure, in $S^n$. 
See \cite{Ogasa98n}. \\

The author found  a high dimensional version of 
a well-known identity 

Alex($K_+)-$Alex($K_-$)=$(t-1)\cdot$Alex($K_0$). 

\noindent 
It is associated with the $(p,q)$-pass-move 
while the 1-dimensional case is with the crossing change. 
See 
\cite{Ogasa09, 
Ogasanu, 
OgasaZ}.  
\\

Kauffman and the author found some relation between  $(p,q)$-pass-moves  and knot products. See 
\cite{KauffmanOgasa, KauffmanOgasa2, KauffmanOgasa3}

\vskip1cm
There are many other exciting open problems on these local moves. 
Please try.

%\vskip1cm  See  the author's paper \cite{Ogasa98n, Ogasa09} and 

\np  
%\subsection*
\subsection 
{{\bf To the advanced readers.} 
Twist-moves on $(2m+1)$-dimensional knots}
\label{Appendix}   
%Note to the advanced readers:  
We introduce another  local move on high dimensional knots. 
It is called the  twist-move.  
We define 
the {\it twist-move} on $(2m+1)$-dimensional knots in $\R^{2m+3}$.

\vskip10mm
Before that, 
we review the crossing-change on 1-links, 
 generalize it, 
and define the twist-move on high dimensional knots.

\begin{figure}[H]
\bigbreak
\hskip20mm\includegraphics[width=90mm]{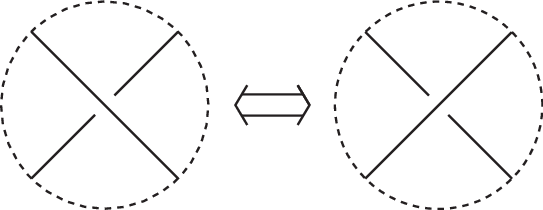} 
\bigbreak\caption{{\bf }\label{twistcross4au}}   
\bigbreak 
\bigbreak 
\end{figure}

\np
 Figure \ref{twistcross4au} draws  
 the crossing-change of 1-links.

\begin{figure}[H]
\bigbreak 
\bigbreak
\hskip20mm\includegraphics[width=90mm]{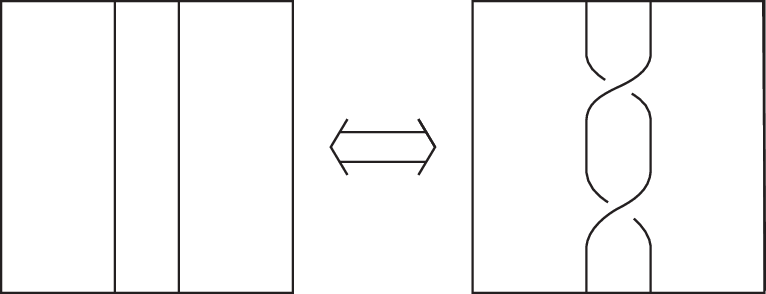} 
\bigbreak\caption{{\bf 
A local move on 1-links.
}\label{twistcross1au}}   
\bigbreak 
\bigbreak 
\end{figure}

The local move drawn in Figure \ref{twistcross1au} 
is the same as the crossing change.

\begin{figure}[H]
\bigbreak
\hskip20mm\includegraphics[width=90mm]{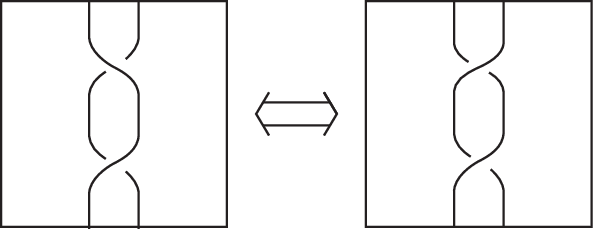} 
\bigbreak\caption{{\bf }\label{twistcross2au}}   
\bigbreak 
\end{figure}

\np
The reason is as follows. 
The left figure of Figure \ref{twistcross1au} is the same as  
the left one of Figure \ref{twistcross2au}.

\begin{figure}[H]
\bigbreak
\hskip20mm\includegraphics[width=90mm]{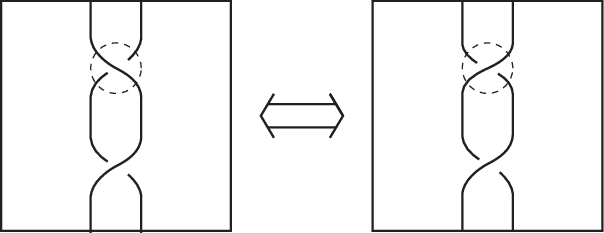} 
\bigbreak\caption{{\bf }\label{twistcross3au}}   
\bigbreak 
\end{figure}

See 
%\includegraphics[width=7mm]{B4} 
%\includegraphics[width=1cm]{dc} 
%
%
%\begin{figure}[H]
%\bigbreak
\includegraphics[width=7mm]{B4} 
%\bigbreak\caption{{\bf }\label{B4au}}   
%\bigbreak 
%\end{figure}
%
in Figure \ref{twistcross3au}.

\begin{figure}[H]
\bigbreak 
\bigbreak
\hskip20mm\includegraphics[width=90mm]{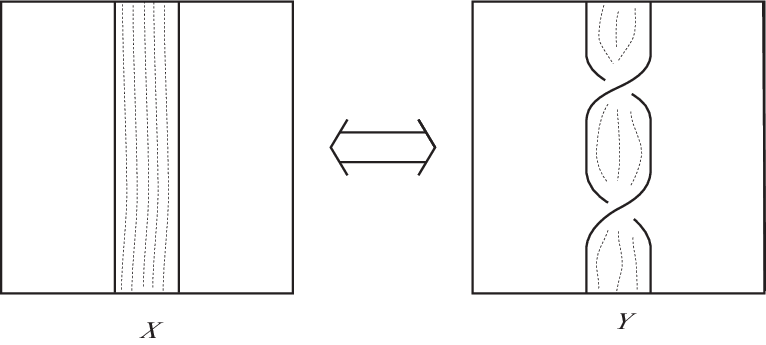} 
\bigbreak\caption{{\bf 
Bands, $X$ and $Y$.
}\label{twistcross1QQ}}   
\bigbreak 
\bigbreak 
\end{figure}

\np
Take a band $B^1\x B^1$, $X$ (resp. $Y$) 
as shown in Figure \ref{twistcross1QQ}.

\begin{figure}[H]
\bigbreak
\hskip20mm\includegraphics[width=90mm]{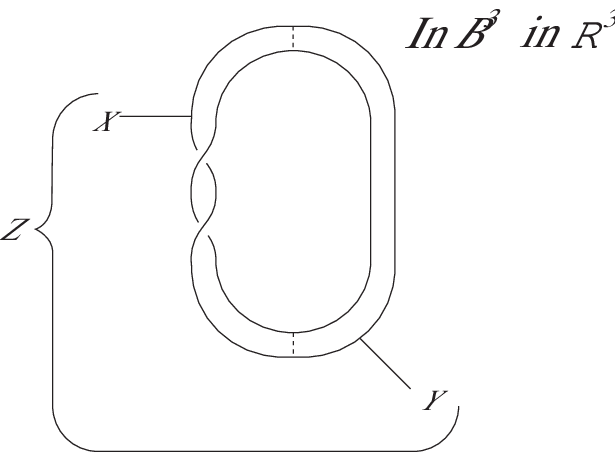}
\bigbreak\caption{{\bf }\label{twistHopfau}}   
\bigbreak 
\end{figure}

Then 
$X$ and $Y$ have the following property. 
We can assume that 
if we coincide the two 3-balls 
in Figure \ref{twistcross1au},    
$X$ and $Y$ touches each other only in the boundary of the 3-ball 
as shown in Figure  \ref{twistHopfau}.   
The union of $X$ and $Y$ is drawn in Figure \ref{twistHopfau}.  
It is an anulus which twisted one time in this case. 
The boundary is the Hopf link.

We use this property when 
we define the twist move.

\np
We define 
the {\it twist-move} on $(2m+1)$-dimensional knots in $\R^{2m+3}$. 

Note that the twist-move in the $2m+1=1$ case, i.e.,  
the twist-move on 1-knots,  is the crossing change. 

\vskip10mm

Take a $(2m+3)$-ball $B^{2m+3}$ trivially embedded in $\R^{2m+3}$.

Regard $B^{2m+3}$ as $B^1\x B^{m+1}\x B^{m+1}$. 

\begin{figure}[H]
\bigbreak
\includegraphics[width=111mm]{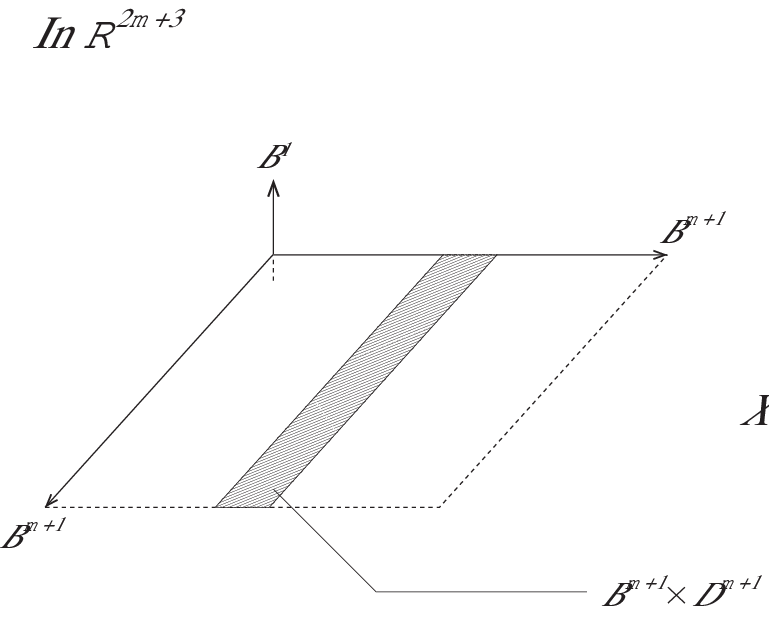} 
\bigbreak\caption{{\bf }\label{twist+au}}   
\bigbreak 
\end{figure}

Take $B^{m+1}\x D^{m+1}$ trivially embedded in $B^{2m+3}$ 
as shown in Figure \ref{twist+au}. 
Here,  $D^{m+1}$, which is defined earlier, is a smaller ball than $B^{m+1}$ . 
$D^{m+1}$ is embedded in $B^{m+1}$ . 
It is drawn very conceptually.

This embedding of $B^{m+1}\x D^{m+1}$ is called $X$.

\np  

\begin{figure}[H]
\bigbreak
\includegraphics[width=111mm]{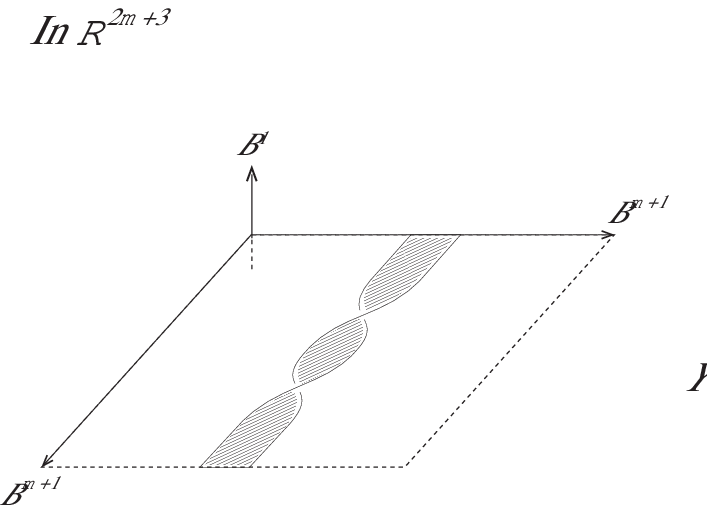} 
\bigbreak\caption{{\bf }\label{twist-au}}   
\bigbreak 
\end{figure}

Take another embedding of $B^{m+1}\x D^{m+1}$ in $B^{2m+3}$ 
as shown in Figure \ref{twist-au}.

Let this embedding be called $Y$. 
If we take $X$ and $Y$ in $B^{2m+3}$ together, 
the intersection of $X$ and the boundary of $B^{2m+3}$ coincides with 
the intersection of $Y$ and the boundary of $B^{2m+3}$
and 
the intersection of $X$ and $Y$ is included in the boundary of $B^{2m+3}$. 
The union of $X$ and $Y$ is called $Z$.

%\vskip5mm
%Note the following difference: In the case of $(p,q)$-pass-moves, $X_+$, $X_-$, $Y_+$, $Y_-$ are $B^1\x B^p\x B^q$.  Here, $X$ and $Y$ are $B^{m+1}\x D^{m+1}$ in $B^{2m+3}$ not $B^{2m+3}=B^1\x B^{m+1}\x B^{m+1}$. 
%
%The reason of this difference is just because of the convenience of the explanation. 

\np

\begin{figure}[H]
\bigbreak
\includegraphics[width=105mm]{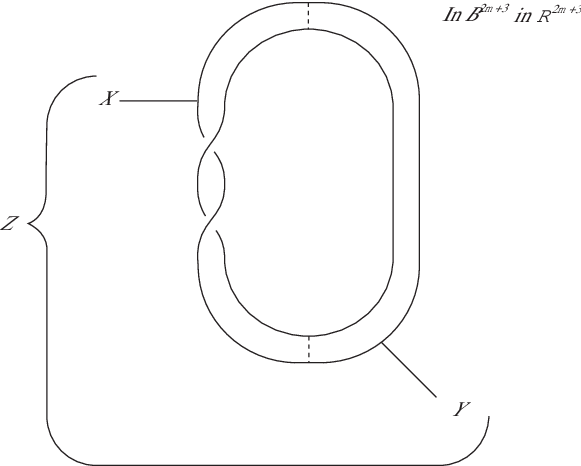}   
\bigbreak\caption{{\bf }\label{twistXYau}}   
\bigbreak 
\end{figure}

$Z$ in Figure \ref{twistXYau} 
is `twisted one time' in the following meaning. 
$X$ and $Y$ are \\$B^{m+1}\x B^{m+1}$. 
$Z$ is $S^{m+1}\x B^{m+1}$. 
$B^{m+1}\x \{c\}$ in $X$ and that in $Y$ meet at each boundary and 
make an $(m+1)$-sphere $S^{m+1}$ 
as drawn in Figure \ref{twistXYau}, 
where $\{c\}$ is the center of $D^{m+1}$. 
This  $S^{m+1}$ is, of course, included in $Z$. 
We can push off $S^{m+1}$ into the normal direction of $Z$ in $\R^{2m+3}$.

Mathematically speaking, we have the following. 
$Z$ is the total space of the $D^{m+1}$-bundle over $S^{m+1}$ 
associated with the tangent bundle of $S^{m+1}$.   
$Z$ is a codimension one smooth submanifold of $\R^{2m+3}$.  
$Z$ has the trivial normal bundle. Hence we can push off $S^{m+1}$. 

The resulting $(m+1)$-dimensional sphere by pushing off is called $\widetilde{S^{m+1}}$.

We suppose that $S^{m+1}$ and $\widetilde{S^{m+1}}$ satisfy the following condition: 
$S^{m+1}$ bounds an embedded $(m+2)$-dimensional  disc and 
the intersection of $\widetilde{S^{m+1}}$ and the  $(m+2)$-dimensional  disc is a single point.

%Note that $S^{m+1}$ and $\widetilde{S^{m+1}}$ are in $\R^{2m+3}$. 
% The problem associated with the last figur in \S\ref{S1Sn-1kotae} 
%is related to this situation. 

\np

If you know `Seifert pairing' and `handle decomposition', 
the explanation will become much clearer: 
We let $Z$ satisfy that 
the Seifert pairing of 
a codimension two submanifold $\partial Z\subset B^{n+2}$ (resp. $S^{n+2}$)  
associated with $Z$ is a $1\x1$-matrix (1). 
See the author's paper \cite{Ogasa09} 
for detail.

For now, consider the $2m+1=1$ case as shown 
in Figures  \ref{twistcross1au}-\ref{twistcross4au}   
and imagine the higher dimensional case.

\vskip1cm

The boundary of $D^{m+1}$ is an $m$-sphere $S^m$. 
Therefore,  $B^{m+1}\x S^m$ is included in the boundary of  $B^{m+1}\x D^{m+1}$. 
This  $B^{m+1}\x S^m$ in the boundary of $X$ (resp. $Y$) is called $X'$ (resp. $Y'$). 

Let $K_+$ and $K_-$ be $(2m+1)$-dimensional knots in $\R^{2m+3}$. 
Suppose that 
$K_+$ and $K_-$ differ only in $B^{2m+3}$ 
and that 
the intersection of $K_+$ (resp. $K_-$) and $B^{2m+3}$ is $X'$ (resp. $Y'$).
Then we say that 
$K_+$ (resp. $K_-$) is obtained from $K_-$ (resp. $K_+$) 
by one twist-move. 

If $K_+$ (resp. $K_-$) is obtained from $K_-$ (resp. $K_+$) 
by a sequence of a finite number of twist-moves, 
then we say that $K_+$ and $K_-$ are twist-move-equivalent.

\vskip1cm
Note:  $K_+$ and $K_-$ are sometimes diffeomorphic in some cases.
They are homeomorphic but not diffeomorphic in some cases.  
They are not homeomorphic in the other cases.

\np

We found infinitely many nontrivial $n$-knots 
which are twist-move-equivalent to the trivial knot 
(resp. which are not twist-move-equivalent to the trivial knot).

\vskip10mm

The author found some relations between some invariants of high dimensional knots and twist-moves: 
He found  a $(4k+1)$ dimensional version of 
a well-known identity 

Alex($K_+)-$Alex($K_-$)=$(t-1)\cdot$Alex($K_0$). 

\noindent 
It is associated with the twist move while the 1-dimensional case is with the crossing change. 
See 
\cite{Ogasa09, 
Ogasanu, 
OgasaZ}.

\vskip10mm

Kauffman and the author found some relations between  twist-moves  and knot products: 
A natural $(2\mu+1)$ dimensional version of 
a well-known identity 

Alex($K_+)-$Alex($K_-$)=$(t-1)\cdot$Alex($K_0$)

\noindent 
associated with the twist move is a little different, and is 

Alex($K_+)-$Alex($K_-$)=$(t+1)\cdot$Alex($K_0$).  

\noindent 
Note the left side of $=$. $(t-1)$ is replaced with $(t+1)$. 
It is a new type identity. 
See \cite{KauffmanOgasa,  KauffmanOgasa2, KauffmanOgasa3}.

\vskip10mm
The author proved that 
for any natural number $\nu$, 
there is a nontrivial $(4k+1)$-knot 
which is obtained 
from the trivial knot by $\nu$ times of twist moves,  
but which is not obtained 
from the trivial knot by $\xi$ times of twist moves if $\xi<\nu$.

\vskip1cm
There are many other exciting open problems on these local moves. 
Please try.

\np

\subsection 
{%{\bf To the advanced readers.} 
Cross-ring-moves on 2-knots}
\label{Appendix2}

 We introduce another local move, 
`cross-ring-moves on 2-dimensional knots in $\R^4$',  
after defining the concept of local moves on $n$-knots. 
%Local moves are defined not only for $n$-knots but also $n$-links and submanifolds in the same fashion.

\vskip1cm
{\bf To the advanced readers.} 
Let $K_+$ be an $n$-knot in $\R^{n+2}$. 
Let $A$ be fixed in $\R^{n+2}$. 
Fix $P$ (resp. $N$) in $A$.  
Suppose 
 that  $K_+$ and $K_-$ only differ in $A$
and 
that 
the intersection of $K_+$ (resp. $K_-$) and $A$ is $P$ (resp. $N$). 
Then we say that $K_-$ is obtained from $K_+$ by this local move.    

$A$ is not necessarily a $(n+2)$-dimensional ball. 
See Definition \ref{crossring}

\bigbreak
%\noindent  {\bf To  the advanced readers.}    
We often want to suppose that $A$ is a compact $(n+2)$-submanifold of $R^{n+2}$. 
\bigbreak

Of course this definition is modified in many ways. 
For example, we can replace 
`$n$-knot' with `$n$-dimensional submanifold of $(n+2)$-manifold'. 
Indeed the local moves introduced in this article are defined 
not only for $n$-knots 
but also for $n$-dimensional submanifolds of $(n+2)$-manifolds as well.

\np

\begin{defn}\label{crossring}
Let $K$ be a 2-knot in $\R^4$. 
Embed $S^1\x D^3$ trivially in $\R^4$, 
where 
$S^1$ is a circle and where $D^3$ is the close 3-ball. 
Suppose the following. 

\smallbreak
\noindent (1) 
$K\cap (S^1\x D^3)$ is 
$(S^1\x I)\amalg(S^1\x I)$, where $I$ is the interval. 

\smallbreak
\noindent (2) 
$K\cap (S^1\x D^3)$ is 

\bigbreak
\hskip30mm
%\begin{figure}[H]\bigbreak
\input XO+.tex   
%\bigbreak\caption{{\bf }\label{XO+au}}   \bigbreak \end{figure}
\bigbreak

\noindent 
where we have the following.

\smallbreak
\noindent (i) 
The bold line and its interior in the above figure 
represent the 3-disc $D^3$ 

\smallbreak
\noindent (ii) 
The arrows which are written by fine lines represent 
a submanifold of $K$. 

\smallbreak
\noindent (iii) 
$S^1\x$ (each of the two arrows) is  
each of the above $(S^1\x I)\amalg(S^1\x I)$. 

\vskip3mm
\noindent
Fix this $S^1\x D^3$ in $\R^4$. 

Let  $K'$ be a 2-knot with the following properties. 

\smallbreak
\noindent(1)
$K\cap (S^1\x D^3)$ is 

\bigbreak
\hskip30mm
%\begin{figure}[H]\bigbreak
\input XO-.tex   
%\bigbreak\caption{{\bf }\label{XO+au}}   \bigbreak \end{figure}
\bigbreak

\noindent
in the above chart. 

\smallbreak
\noindent(2) 
%$K'\cap\overline{S^4-(S^1\x D^3)}$=
%$K\cap\overline{S^4-(S^1\x D^3)}$. 
$K'\cap (\R^4-(\text{the interior of }(S^1\x D^3)))$

=
$K\cap (\R^4-(\text{the interior of }(S^1\x D^3)))$.

\np

Then we say that 
$K'$ is obtained from $K$ by one {\it cross-ring-move}.

If $K''$ is obtained from $K$ by a sequence of a finite number of cross-ring-moves 
then we say that 
$K''$ is {\it cross-ring-move-equivalent} to $K$. 
\end{defn}

Needless to say, cross-ring moves are defined not only for 2-knots but also for 2-links in the same fashion.

\vskip1cm 

We could prove that  
all $k$-twist spun knots are cross-ring-move-equivalent to the trivial knot.

It is natural to ask the following question. 

\begin{que}\label{crm}
Are all 2-knots cross-ring-move-equivalent to the trivial knot?
\end{que}

%We do not know the answer for now 
%as of the publishing of the first edition of the English version of this book. 
%We think that it is open. 
Please solve. 

\bigbreak
\bigbreak
Furthermore the following is an open problem. 
If $n=1$, the answer is easy. It is the crossing change in Figure \ref{Alabama}.    

\begin{que}\label{unk}
Is there an explicit unknotting operation 
on $n$-knots ($n\geqq2$)?
\end{que} 

Please solve.

\np

\section{Comments}\label{comments} 
We introduce literature on high dimensional knot theory. 
We would be sorry if we forgot citing other important papers in the following references. 

There are many other important themes associated with high dimensional knots, which are not  introduced in this article: 
knot cobordism, link cobordism, the complements of $n$-knots ($n=1$ and $n\geqq2$),  fibered $n$-knots ($n=1$ and $n\geqq2$), knot products, 
intersectional pair of knots, etc. 
Some of the following papers and articles include them. 

\bigbreak
The following is an outstanding open problem. See 
\cite{CochranOrr, Ogasa99D, OgasaOCSL, Orr}. 

\begin{que}\label{slice}
Are all  even dimensional links slice?
\end{que}

Please solve.

\bigbreak\noindent
Eiji Ogasa\newline 
 Computer Science, Meijigakuin University, Yokohama, Kanagawa, 244-8539, Japan 
\newline 
ogasa@mail1.meijigkakuin.ac.jp 
\quad
pqr100pqr100@yahoo.co.jp  

\smallbreak
\noindent
Manuscripts of the author's papers can be obtained from his website 

\smallbreak
\noindent
http://www.geocities.jp/n\_dimension\_n\_dimension/list.html

\smallbreak
\noindent
Don't forget the three \_ in this address if you type it. 

\smallbreak
\noindent
You can find his website by typing in the author's name `Eiji Ogasa' in the search engine 
although you will not type the address of the website. 
\end{document}

%% file: 0Y207.tex
%WinTpicVersion2.15
\unitlength 0.1in
\begin{picture}(47.96,18.84)(7.89,-25.53)
% CIRCLE 0 0 3 0
% 4 1710 2030 2470 2550 2470 2550 2470 2550
% 
\special{pn 20}%
\special{ar 1710 1630 921 921  0.0000000 6.2831853}%
% SPLINE 0 0 3 0
% 60 4560 2480 4380 2380 4320 2310 4290 2270 4210 2090 4130 1930 4120 1890 4080 1770 4060 1670 4060 1450 4060 1370 4110 1280 4170 1210 4250 1170 4330 1140 4370 1100 4390 1090 4450 1070 4480 1070 4580 1100 4610 1130 4710 1230 4780 1370 4800 1450 4820 1710 4830 1850 4830 1850 4800 2270 4700 2570 4580 2720 4440 2800 4330 2840 4190 2920 4010 2950 3800 2810 3760 2590 3760 2500 3760 2360 3790 2210 3910 2010 4010 1880 4110 1780 4270 1710 4390 1680 4660 1650 5110 1580 5200 1600 5360 1690 5440 1730 5530 2000 5580 2240 5580 2370 5500 2530 5360 2600 5190 2630 5060 2650 4880 2620 4800 2580 4570 2500 4550 2490
% 
\special{pn 20}%
\special{pa 4560 2080}%
\special{pa 4530 2067}%
\special{pa 4500 2054}%
\special{pa 4470 2041}%
\special{pa 4442 2025}%
\special{pa 4415 2008}%
\special{pa 4390 1989}%
\special{pa 4367 1967}%
\special{pa 4346 1943}%
\special{pa 4326 1918}%
\special{pa 4306 1893}%
\special{pa 4288 1867}%
\special{pa 4272 1839}%
\special{pa 4259 1810}%
\special{pa 4247 1781}%
\special{pa 4235 1751}%
\special{pa 4224 1721}%
\special{pa 4210 1691}%
\special{pa 4196 1662}%
\special{pa 4180 1634}%
\special{pa 4164 1606}%
\special{pa 4149 1578}%
\special{pa 4136 1549}%
\special{pa 4127 1518}%
\special{pa 4119 1487}%
\special{pa 4110 1457}%
\special{pa 4099 1426}%
\special{pa 4089 1396}%
\special{pa 4079 1366}%
\special{pa 4071 1335}%
\special{pa 4064 1303}%
\special{pa 4060 1272}%
\special{pa 4058 1240}%
\special{pa 4058 1208}%
\special{pa 4059 1176}%
\special{pa 4060 1144}%
\special{pa 4061 1112}%
\special{pa 4061 1080}%
\special{pa 4060 1048}%
\special{pa 4058 1015}%
\special{pa 4058 983}%
\special{pa 4065 953}%
\special{pa 4079 924}%
\special{pa 4097 897}%
\special{pa 4117 871}%
\special{pa 4137 845}%
\special{pa 4158 821}%
\special{pa 4182 801}%
\special{pa 4210 785}%
\special{pa 4240 773}%
\special{pa 4272 764}%
\special{pa 4303 755}%
\special{pa 4330 740}%
\special{pa 4352 716}%
\special{pa 4377 696}%
\special{pa 4406 683}%
\special{pa 4437 672}%
\special{pa 4468 669}%
\special{pa 4501 672}%
\special{pa 4533 678}%
\special{pa 4563 689}%
\special{pa 4589 708}%
\special{pa 4611 731}%
\special{pa 4634 754}%
\special{pa 4657 776}%
\special{pa 4680 798}%
\special{pa 4702 821}%
\special{pa 4722 846}%
\special{pa 4739 873}%
\special{pa 4754 901}%
\special{pa 4766 931}%
\special{pa 4777 962}%
\special{pa 4787 993}%
\special{pa 4795 1024}%
\special{pa 4801 1055}%
\special{pa 4806 1086}%
\special{pa 4809 1118}%
\special{pa 4812 1150}%
\special{pa 4814 1182}%
\special{pa 4816 1214}%
\special{pa 4817 1246}%
\special{pa 4818 1278}%
\special{pa 4820 1310}%
\special{pa 4822 1342}%
\special{pa 4825 1374}%
\special{pa 4827 1406}%
\special{pa 4829 1438}%
\special{pa 4831 1470}%
\special{pa 4832 1502}%
\special{pa 4832 1534}%
\special{pa 4832 1566}%
\special{pa 4831 1598}%
\special{pa 4830 1630}%
\special{pa 4828 1662}%
\special{pa 4825 1694}%
\special{pa 4822 1726}%
\special{pa 4818 1757}%
\special{pa 4814 1789}%
\special{pa 4809 1821}%
\special{pa 4803 1853}%
\special{pa 4797 1885}%
\special{pa 4791 1916}%
\special{pa 4783 1947}%
\special{pa 4775 1979}%
\special{pa 4766 2010}%
\special{pa 4757 2040}%
\special{pa 4746 2070}%
\special{pa 4734 2100}%
\special{pa 4721 2129}%
\special{pa 4706 2158}%
\special{pa 4691 2186}%
\special{pa 4674 2213}%
\special{pa 4655 2240}%
\special{pa 4635 2265}%
\special{pa 4613 2289}%
\special{pa 4590 2311}%
\special{pa 4566 2331}%
\special{pa 4539 2350}%
\special{pa 4512 2367}%
\special{pa 4483 2381}%
\special{pa 4454 2395}%
\special{pa 4424 2406}%
\special{pa 4393 2416}%
\special{pa 4363 2427}%
\special{pa 4334 2438}%
\special{pa 4305 2453}%
\special{pa 4278 2469}%
\special{pa 4250 2486}%
\special{pa 4223 2502}%
\special{pa 4195 2518}%
\special{pa 4165 2531}%
\special{pa 4134 2541}%
\special{pa 4103 2548}%
\special{pa 4071 2552}%
\special{pa 4038 2553}%
\special{pa 4006 2549}%
\special{pa 3973 2542}%
\special{pa 3942 2531}%
\special{pa 3912 2517}%
\special{pa 3883 2499}%
\special{pa 3857 2479}%
\special{pa 3833 2456}%
\special{pa 3813 2430}%
\special{pa 3796 2403}%
\special{pa 3783 2374}%
\special{pa 3774 2343}%
\special{pa 3768 2312}%
\special{pa 3763 2280}%
\special{pa 3761 2247}%
\special{pa 3760 2214}%
\special{pa 3760 2182}%
\special{pa 3760 2149}%
\special{pa 3760 2118}%
\special{pa 3760 2086}%
\special{pa 3759 2054}%
\special{pa 3759 2022}%
\special{pa 3759 1990}%
\special{pa 3760 1958}%
\special{pa 3763 1926}%
\special{pa 3768 1894}%
\special{pa 3774 1862}%
\special{pa 3783 1831}%
\special{pa 3793 1801}%
\special{pa 3806 1772}%
\special{pa 3821 1744}%
\special{pa 3837 1716}%
\special{pa 3854 1689}%
\special{pa 3872 1662}%
\special{pa 3891 1636}%
\special{pa 3910 1610}%
\special{pa 3929 1584}%
\special{pa 3948 1559}%
\special{pa 3968 1533}%
\special{pa 3987 1508}%
\special{pa 4008 1483}%
\special{pa 4028 1458}%
\special{pa 4050 1434}%
\special{pa 4072 1411}%
\special{pa 4096 1390}%
\special{pa 4122 1372}%
\special{pa 4150 1356}%
\special{pa 4179 1342}%
\special{pa 4208 1330}%
\special{pa 4239 1320}%
\special{pa 4270 1310}%
\special{pa 4301 1301}%
\special{pa 4332 1293}%
\special{pa 4363 1285}%
\special{pa 4394 1279}%
\special{pa 4426 1274}%
\special{pa 4457 1270}%
\special{pa 4489 1267}%
\special{pa 4521 1264}%
\special{pa 4553 1262}%
\special{pa 4585 1259}%
\special{pa 4617 1256}%
\special{pa 4649 1252}%
\special{pa 4681 1247}%
\special{pa 4712 1241}%
\special{pa 4744 1234}%
\special{pa 4776 1227}%
\special{fp}%
\special{pa 4776 1227}%
\special{pa 4807 1219}%
\special{pa 4839 1211}%
\special{pa 4870 1204}%
\special{pa 4901 1197}%
\special{pa 4933 1191}%
\special{pa 4965 1186}%
\special{pa 4996 1182}%
\special{pa 5028 1179}%
\special{pa 5060 1178}%
\special{pa 5091 1179}%
\special{pa 5123 1181}%
\special{pa 5155 1187}%
\special{pa 5186 1195}%
\special{pa 5216 1207}%
\special{pa 5244 1221}%
\special{pa 5271 1237}%
\special{pa 5298 1254}%
\special{pa 5326 1271}%
\special{pa 5354 1287}%
\special{pa 5384 1301}%
\special{pa 5413 1314}%
\special{pa 5440 1330}%
\special{pa 5462 1350}%
\special{pa 5479 1374}%
\special{pa 5492 1402}%
\special{pa 5502 1433}%
\special{pa 5510 1466}%
\special{pa 5516 1500}%
\special{pa 5521 1535}%
\special{pa 5525 1569}%
\special{pa 5531 1604}%
\special{pa 5537 1637}%
\special{pa 5544 1669}%
\special{pa 5551 1700}%
\special{pa 5558 1731}%
\special{pa 5565 1762}%
\special{pa 5572 1792}%
\special{pa 5577 1823}%
\special{pa 5582 1854}%
\special{pa 5584 1886}%
\special{pa 5585 1918}%
\special{pa 5583 1950}%
\special{pa 5578 1982}%
\special{pa 5570 2014}%
\special{pa 5558 2045}%
\special{pa 5544 2074}%
\special{pa 5526 2101}%
\special{pa 5505 2125}%
\special{pa 5481 2146}%
\special{pa 5454 2164}%
\special{pa 5425 2178}%
\special{pa 5394 2190}%
\special{pa 5363 2199}%
\special{pa 5332 2207}%
\special{pa 5300 2213}%
\special{pa 5269 2218}%
\special{pa 5237 2223}%
\special{pa 5205 2227}%
\special{pa 5174 2233}%
\special{pa 5143 2239}%
\special{pa 5111 2244}%
\special{pa 5079 2248}%
\special{pa 5047 2251}%
\special{pa 5015 2250}%
\special{pa 4983 2248}%
\special{pa 4951 2242}%
\special{pa 4919 2234}%
\special{pa 4889 2224}%
\special{pa 4860 2210}%
\special{pa 4831 2196}%
\special{pa 4803 2181}%
\special{pa 4774 2169}%
\special{pa 4743 2158}%
\special{pa 4713 2149}%
\special{pa 4682 2139}%
\special{pa 4651 2130}%
\special{pa 4621 2120}%
\special{pa 4591 2109}%
\special{pa 4562 2096}%
\special{pa 4550 2090}%
\special{sp}%
% BOX 3 5 2 0
% 2 4010 1730 4170 1910
% 
\special{pn 4}%
\special{sh 0}%
\special{pa 4010 1330}%
\special{pa 4170 1330}%
\special{pa 4170 1510}%
\special{pa 4010 1510}%
\special{pa 4010 1330}%
\special{ip}%
% BOX 3 5 2 0
% 2 4730 1540 4890 1720
% 
\special{pn 4}%
\special{sh 0}%
\special{pa 4730 1140}%
\special{pa 4890 1140}%
\special{pa 4890 1320}%
\special{pa 4730 1320}%
\special{pa 4730 1140}%
\special{ip}%
% BOX 3 5 2 0
% 2 4630 2480 4790 2660
% 
\special{pn 4}%
\special{sh 0}%
\special{pa 4630 2080}%
\special{pa 4790 2080}%
\special{pa 4790 2260}%
\special{pa 4630 2260}%
\special{pa 4630 2080}%
\special{ip}%
% LINE 0 0 3 0
% 2 4130 1930 4060 1720
% 
\special{pn 20}%
\special{pa 4130 1530}%
\special{pa 4060 1320}%
\special{fp}%
% LINE 0 0 3 0
% 2 4710 1630 4900 1600
% 
\special{pn 20}%
\special{pa 4710 1230}%
\special{pa 4900 1200}%
\special{fp}%
% LINE 0 0 3 0
% 2 4750 2430 4630 2680
% 
\special{pn 20}%
\special{pa 4750 2030}%
\special{pa 4630 2280}%
\special{fp}%
\end{picture}%

%% file: Y207.tex
%WinTpicVersion2.15
\unitlength 0.1in
\begin{picture}(47.96,21.21)(7.89,-27.90)
% CIRCLE 0 0 3 0
% 4 1710 2030 2470 2550 2470 2550 2470 2550
% 
\special{pn 20}%
\special{ar 1710 1630 921 921  0.0000000 6.2831853}%
% SPLINE 0 0 3 0
% 60 4560 2480 4380 2380 4320 2310 4290 2270 4210 2090 4130 1930 4120 1890 4080 1770 4060 1670 4060 1450 4060 1370 4110 1280 4170 1210 4250 1170 4330 1140 4370 1100 4390 1090 4450 1070 4480 1070 4580 1100 4610 1130 4710 1230 4780 1370 4800 1450 4820 1710 4830 1850 4830 1850 4800 2270 4700 2570 4580 2720 4440 2800 4330 2840 4190 2920 4010 2950 3800 2810 3760 2590 3760 2500 3760 2360 3790 2210 3910 2010 4010 1880 4110 1780 4270 1710 4390 1680 4660 1650 5110 1580 5200 1600 5360 1690 5440 1730 5530 2000 5580 2240 5580 2370 5500 2530 5360 2600 5190 2630 5060 2650 4880 2620 4800 2580 4570 2500 4550 2490
% 
\special{pn 20}%
\special{pa 4560 2080}%
\special{pa 4530 2067}%
\special{pa 4500 2054}%
\special{pa 4470 2041}%
\special{pa 4442 2025}%
\special{pa 4415 2008}%
\special{pa 4390 1989}%
\special{pa 4367 1967}%
\special{pa 4346 1943}%
\special{pa 4326 1918}%
\special{pa 4306 1893}%
\special{pa 4288 1867}%
\special{pa 4272 1839}%
\special{pa 4259 1810}%
\special{pa 4247 1781}%
\special{pa 4235 1751}%
\special{pa 4224 1721}%
\special{pa 4210 1691}%
\special{pa 4196 1662}%
\special{pa 4180 1634}%
\special{pa 4164 1606}%
\special{pa 4149 1578}%
\special{pa 4136 1549}%
\special{pa 4127 1518}%
\special{pa 4119 1487}%
\special{pa 4110 1457}%
\special{pa 4099 1426}%
\special{pa 4089 1396}%
\special{pa 4079 1366}%
\special{pa 4071 1335}%
\special{pa 4064 1303}%
\special{pa 4060 1272}%
\special{pa 4058 1240}%
\special{pa 4058 1208}%
\special{pa 4059 1176}%
\special{pa 4060 1144}%
\special{pa 4061 1112}%
\special{pa 4061 1080}%
\special{pa 4060 1048}%
\special{pa 4058 1015}%
\special{pa 4058 983}%
\special{pa 4065 953}%
\special{pa 4079 924}%
\special{pa 4097 897}%
\special{pa 4117 871}%
\special{pa 4137 845}%
\special{pa 4158 821}%
\special{pa 4182 801}%
\special{pa 4210 785}%
\special{pa 4240 773}%
\special{pa 4272 764}%
\special{pa 4303 755}%
\special{pa 4330 740}%
\special{pa 4352 716}%
\special{pa 4377 696}%
\special{pa 4406 683}%
\special{pa 4437 672}%
\special{pa 4468 669}%
\special{pa 4501 672}%
\special{pa 4533 678}%
\special{pa 4563 689}%
\special{pa 4589 708}%
\special{pa 4611 731}%
\special{pa 4634 754}%
\special{pa 4657 776}%
\special{pa 4680 798}%
\special{pa 4702 821}%
\special{pa 4722 846}%
\special{pa 4739 873}%
\special{pa 4754 901}%
\special{pa 4766 931}%
\special{pa 4777 962}%
\special{pa 4787 993}%
\special{pa 4795 1024}%
\special{pa 4801 1055}%
\special{pa 4806 1086}%
\special{pa 4809 1118}%
\special{pa 4812 1150}%
\special{pa 4814 1182}%
\special{pa 4816 1214}%
\special{pa 4817 1246}%
\special{pa 4818 1278}%
\special{pa 4820 1310}%
\special{pa 4822 1342}%
\special{pa 4825 1374}%
\special{pa 4827 1406}%
\special{pa 4829 1438}%
\special{pa 4831 1470}%
\special{pa 4832 1502}%
\special{pa 4832 1534}%
\special{pa 4832 1566}%
\special{pa 4831 1598}%
\special{pa 4830 1630}%
\special{pa 4828 1662}%
\special{pa 4825 1694}%
\special{pa 4822 1726}%
\special{pa 4818 1757}%
\special{pa 4814 1789}%
\special{pa 4809 1821}%
\special{pa 4803 1853}%
\special{pa 4797 1885}%
\special{pa 4791 1916}%
\special{pa 4783 1947}%
\special{pa 4775 1979}%
\special{pa 4766 2010}%
\special{pa 4757 2040}%
\special{pa 4746 2070}%
\special{pa 4734 2100}%
\special{pa 4721 2129}%
\special{pa 4706 2158}%
\special{pa 4691 2186}%
\special{pa 4674 2213}%
\special{pa 4655 2240}%
\special{pa 4635 2265}%
\special{pa 4613 2289}%
\special{pa 4590 2311}%
\special{pa 4566 2331}%
\special{pa 4539 2350}%
\special{pa 4512 2367}%
\special{pa 4483 2381}%
\special{pa 4454 2395}%
\special{pa 4424 2406}%
\special{pa 4393 2416}%
\special{pa 4363 2427}%
\special{pa 4334 2438}%
\special{pa 4305 2453}%
\special{pa 4278 2469}%
\special{pa 4250 2486}%
\special{pa 4223 2502}%
\special{pa 4195 2518}%
\special{pa 4165 2531}%
\special{pa 4134 2541}%
\special{pa 4103 2548}%
\special{pa 4071 2552}%
\special{pa 4038 2553}%
\special{pa 4006 2549}%
\special{pa 3973 2542}%
\special{pa 3942 2531}%
\special{pa 3912 2517}%
\special{pa 3883 2499}%
\special{pa 3857 2479}%
\special{pa 3833 2456}%
\special{pa 3813 2430}%
\special{pa 3796 2403}%
\special{pa 3783 2374}%
\special{pa 3774 2343}%
\special{pa 3768 2312}%
\special{pa 3763 2280}%
\special{pa 3761 2247}%
\special{pa 3760 2214}%
\special{pa 3760 2182}%
\special{pa 3760 2149}%
\special{pa 3760 2118}%
\special{pa 3760 2086}%
\special{pa 3759 2054}%
\special{pa 3759 2022}%
\special{pa 3759 1990}%
\special{pa 3760 1958}%
\special{pa 3763 1926}%
\special{pa 3768 1894}%
\special{pa 3774 1862}%
\special{pa 3783 1831}%
\special{pa 3793 1801}%
\special{pa 3806 1772}%
\special{pa 3821 1744}%
\special{pa 3837 1716}%
\special{pa 3854 1689}%
\special{pa 3872 1662}%
\special{pa 3891 1636}%
\special{pa 3910 1610}%
\special{pa 3929 1584}%
\special{pa 3948 1559}%
\special{pa 3968 1533}%
\special{pa 3987 1508}%
\special{pa 4008 1483}%
\special{pa 4028 1458}%
\special{pa 4050 1434}%
\special{pa 4072 1411}%
\special{pa 4096 1390}%
\special{pa 4122 1372}%
\special{pa 4150 1356}%
\special{pa 4179 1342}%
\special{pa 4208 1330}%
\special{pa 4239 1320}%
\special{pa 4270 1310}%
\special{pa 4301 1301}%
\special{pa 4332 1293}%
\special{pa 4363 1285}%
\special{pa 4394 1279}%
\special{pa 4426 1274}%
\special{pa 4457 1270}%
\special{pa 4489 1267}%
\special{pa 4521 1264}%
\special{pa 4553 1262}%
\special{pa 4585 1259}%
\special{pa 4617 1256}%
\special{pa 4649 1252}%
\special{pa 4681 1247}%
\special{pa 4712 1241}%
\special{pa 4744 1234}%
\special{pa 4776 1227}%
\special{fp}%
\special{pa 4776 1227}%
\special{pa 4807 1219}%
\special{pa 4839 1211}%
\special{pa 4870 1204}%
\special{pa 4901 1197}%
\special{pa 4933 1191}%
\special{pa 4965 1186}%
\special{pa 4996 1182}%
\special{pa 5028 1179}%
\special{pa 5060 1178}%
\special{pa 5091 1179}%
\special{pa 5123 1181}%
\special{pa 5155 1187}%
\special{pa 5186 1195}%
\special{pa 5216 1207}%
\special{pa 5244 1221}%
\special{pa 5271 1237}%
\special{pa 5298 1254}%
\special{pa 5326 1271}%
\special{pa 5354 1287}%
\special{pa 5384 1301}%
\special{pa 5413 1314}%
\special{pa 5440 1330}%
\special{pa 5462 1350}%
\special{pa 5479 1374}%
\special{pa 5492 1402}%
\special{pa 5502 1433}%
\special{pa 5510 1466}%
\special{pa 5516 1500}%
\special{pa 5521 1535}%
\special{pa 5525 1569}%
\special{pa 5531 1604}%
\special{pa 5537 1637}%
\special{pa 5544 1669}%
\special{pa 5551 1700}%
\special{pa 5558 1731}%
\special{pa 5565 1762}%
\special{pa 5572 1792}%
\special{pa 5577 1823}%
\special{pa 5582 1854}%
\special{pa 5584 1886}%
\special{pa 5585 1918}%
\special{pa 5583 1950}%
\special{pa 5578 1982}%
\special{pa 5570 2014}%
\special{pa 5558 2045}%
\special{pa 5544 2074}%
\special{pa 5526 2101}%
\special{pa 5505 2125}%
\special{pa 5481 2146}%
\special{pa 5454 2164}%
\special{pa 5425 2178}%
\special{pa 5394 2190}%
\special{pa 5363 2199}%
\special{pa 5332 2207}%
\special{pa 5300 2213}%
\special{pa 5269 2218}%
\special{pa 5237 2223}%
\special{pa 5205 2227}%
\special{pa 5174 2233}%
\special{pa 5143 2239}%
\special{pa 5111 2244}%
\special{pa 5079 2248}%
\special{pa 5047 2251}%
\special{pa 5015 2250}%
\special{pa 4983 2248}%
\special{pa 4951 2242}%
\special{pa 4919 2234}%
\special{pa 4889 2224}%
\special{pa 4860 2210}%
\special{pa 4831 2196}%
\special{pa 4803 2181}%
\special{pa 4774 2169}%
\special{pa 4743 2158}%
\special{pa 4713 2149}%
\special{pa 4682 2139}%
\special{pa 4651 2130}%
\special{pa 4621 2120}%
\special{pa 4591 2109}%
\special{pa 4562 2096}%
\special{pa 4550 2090}%
\special{sp}%
% BOX 3 5 2 0
% 2 4010 1730 4170 1910
% 
\special{pn 4}%
\special{sh 0}%
\special{pa 4010 1330}%
\special{pa 4170 1330}%
\special{pa 4170 1510}%
\special{pa 4010 1510}%
\special{pa 4010 1330}%
\special{ip}%
% BOX 3 5 2 0
% 2 4730 1540 4890 1720
% 
\special{pn 4}%
\special{sh 0}%
\special{pa 4730 1140}%
\special{pa 4890 1140}%
\special{pa 4890 1320}%
\special{pa 4730 1320}%
\special{pa 4730 1140}%
\special{ip}%
% BOX 3 5 2 0
% 2 4630 2480 4790 2660
% 
\special{pn 4}%
\special{sh 0}%
\special{pa 4630 2080}%
\special{pa 4790 2080}%
\special{pa 4790 2260}%
\special{pa 4630 2260}%
\special{pa 4630 2080}%
\special{ip}%
% LINE 0 0 3 0
% 2 4130 1930 4060 1720
% 
\special{pn 20}%
\special{pa 4130 1530}%
\special{pa 4060 1320}%
\special{fp}%
% LINE 0 0 3 0
% 2 4710 1630 4900 1600
% 
\special{pn 20}%
\special{pa 4710 1230}%
\special{pa 4900 1200}%
\special{fp}%
% LINE 0 0 3 0
% 2 4750 2430 4630 2680
% 
\special{pn 20}%
\special{pa 4750 2030}%
\special{pa 4630 2280}%
\special{fp}%
% STR 2 0 3 0
% 3 4180 3260 4180 3360 2 0
% The trefoil knot
\put(41.8000,-29.6000){\makebox(0,0)[lb]{The trefoil knot}}%
% STR 2 0 3 0
% 3 860 3110 860 3210 2 0
% The trivial knot
\put(8.6000,-28.1000){\makebox(0,0)[lb]{The trivial knot}}%
\end{picture}%

%% file: I234.tex
%WinTpicVersion2.15
\unitlength 0.1in
\begin{picture}(17.03,4.96)(6.59,-9.91)
% POLYLINE 2 0 3 0
% 7 1150 900 659 1391 1913 1391 2362 895 1150 900 1150 900 1150 900
% 
\special{pn 8}%
\special{pa 1150 500}%
\special{pa 659 991}%
\special{pa 1913 991}%
\special{pa 2362 495}%
\special{pa 1150 500}%
\special{pa 1150 500}%
\special{pa 1150 500}%
\special{fp}%
% STR 2 0 3 1
% 3 2583 1158 2583 1215 2 0
% A part of a new surface
\put(25.8300,-8.1500){\makebox(0,0)[lb]{A part of a new surface}}%
% POLYLINE 2 0 3 2
% 7 1155 2184 664 2673 1919 2673 2368 2178 1155 2184 1155 2184 1155 2184
% 
\special{pn 8}%
\special{pa 1155 1784}%
\special{pa 664 2273}%
\special{pa 1919 2273}%
\special{pa 2368 1778}%
\special{pa 1155 1784}%
\special{pa 1155 1784}%
\special{pa 1155 1784}%
\special{fp}%
% STR 2 0 3 3
% 3 2590 2440 2590 2499 2 0
% A part of a new surface
\put(25.9000,-20.9900){\makebox(0,0)[lb]{A part of a new surface}}%
% POLYLINE 0 0 3 4
% 7 1173 1548 682 2038 1936 2038 2385 1542 1173 1548 1173 1548 1173 1548
% 
\special{pn 20}%
\special{pa 1173 1148}%
\special{pa 682 1638}%
\special{pa 1936 1638}%
\special{pa 2385 1142}%
\special{pa 1173 1148}%
\special{pa 1173 1148}%
\special{pa 1173 1148}%
\special{fp}%
% STR 2 0 3 5
% 3 2608 1805 2608 1863 2 0
% A part of the Boy surface
\put(26.0800,-14.6300){\makebox(0,0)[lb]{A part of the Boy surface}}%
\end{picture}%

%% file: conce1.tex
%WinTpicVersion2.15
\unitlength 0.1in
\begin{picture}(38.99,31.61)(0.40,-33.10)
% LINE 2 0 3 0
% 2 2894 707 852 3696
% 
\special{pn 8}%
\special{pa 2894 307}%
\special{pa 852 3296}%
\special{fp}%
% BOX 2 5 2 0
% 2 2684 1974 2820 2158
% 
\special{pn 8}%
\special{sh 0}%
\special{pa 2684 1574}%
\special{pa 2820 1574}%
\special{pa 2820 1758}%
\special{pa 2684 1758}%
\special{pa 2684 1574}%
\special{ip}%
% ELLIPSE 0 0 3 0
% 4 3140 2478 3804 2982 3804 2982 4037 3204
% 
\special{pn 20}%
\special{ar 3140 2078 664 504  0.8172277 6.2831853}%
\special{ar 3140 2078 664 504  0.0000000 0.7853982}%
% BOX 2 5 2 0
% 2 2684 1974 2820 2158
% 
\special{pn 8}%
\special{sh 0}%
\special{pa 2684 1574}%
\special{pa 2820 1574}%
\special{pa 2820 1758}%
\special{pa 2684 1758}%
\special{pa 2684 1574}%
\special{ip}%
% BOX 2 5 2 0
% 2 2365 1334 2500 1519
% 
\special{pn 8}%
\special{sh 0}%
\special{pa 2365 934}%
\special{pa 2500 934}%
\special{pa 2500 1119}%
\special{pa 2365 1119}%
\special{pa 2365 934}%
\special{ip}%
% ELLIPSE 2 0 3 0
% 4 1946 2306 2783 3437 2857 3425 2857 3425
% 
\special{pn 8}%
\special{ar 1946 1906 837 1131  0.0000000 6.2831853}%
% BOX 2 5 2 0
% 2 2623 2773 2758 2958
% 
\special{pn 8}%
\special{sh 0}%
\special{pa 2623 2373}%
\special{pa 2758 2373}%
\special{pa 2758 2558}%
\special{pa 2623 2558}%
\special{pa 2623 2373}%
\special{ip}%
% LINE 2 5 3 0
% 2 2623 2798 2746 2872
% 
\special{pn 8}%
\special{pa 2623 2398}%
\special{pa 2746 2472}%
\special{ip}%
% LINE 0 0 3 0
% 2 2623 2785 2746 2884
% 
\special{pn 20}%
\special{pa 2623 2385}%
\special{pa 2746 2484}%
\special{fp}%
% BOX 2 5 2 0
% 2 1233 2995 1368 3179
% 
\special{pn 8}%
\special{sh 0}%
\special{pa 1233 2595}%
\special{pa 1368 2595}%
\special{pa 1368 2779}%
\special{pa 1233 2779}%
\special{pa 1233 2595}%
\special{ip}%
% LINE 2 5 3 0
% 2 1331 2982 1208 3179
% 
\special{pn 8}%
\special{pa 1331 2582}%
\special{pa 1208 2779}%
\special{ip}%
% LINE 2 0 3 0
% 2 1344 2970 1135 3290
% 
\special{pn 8}%
\special{pa 1344 2570}%
\special{pa 1135 2890}%
\special{fp}%
% STR 2 0 3 0
% 3 40 3757 40 3880 2 0
% The $zw$-plane
\put(0.4000,-34.8000){\makebox(0,0)[lb]{The $zw$-plane}}%
% STR 2 0 3 0
% 3 114 596 114 719 2 0
% $R^4$
\put(1.1400,-3.1900){\makebox(0,0)[lb]{$R^4$}}%
% STR 2 0 3 0
% 3 630 1961 630 2084 2 0
% $S^1$
\put(6.3000,-16.8400){\makebox(0,0)[lb]{$S^1$}}%
% STR 2 0 3 0
% 3 3939 2466 3939 2589 2 0
% $S^2$
\put(39.3900,-21.8900){\makebox(0,0)[lb]{$S^2$}}%
\end{picture}%

%% file: Figure1.tex
%WinTpicVersion2.15
\unitlength 0.1in
\begin{picture}(56.10,43.70)(8.50,-44.30)
% ELLIPSE 2 0 3 0
% 4 3510 720 4070 970 4070 970 4070 970
% 
\special{pn 8}%
\special{ar 3510 320 560 250  0.0000000 6.2831853}%
% LINE 0 0 3 0
% 2 3520 4350 3520 3130
% 
\special{pn 20}%
\special{pa 3520 3950}%
\special{pa 3520 2730}%
\special{fp}%
% ELLIPSE 2 0 3 0
% 4 3510 4480 4070 4730 4070 4730 4070 4730
% 
\special{pn 8}%
\special{ar 3510 4080 560 250  0.0000000 6.2831853}%
% LINE 2 0 3 0
% 2 4080 730 4080 4480
% 
\special{pn 8}%
\special{pa 4080 330}%
\special{pa 4080 4080}%
\special{fp}%
% LINE 2 0 3 0
% 2 2950 740 2950 4480
% 
\special{pn 8}%
\special{pa 2950 340}%
\special{pa 2950 4080}%
\special{fp}%
% ELLIPSE 2 0 3 0
% 4 1410 730 1970 980 1970 980 1970 980
% 
\special{pn 8}%
\special{ar 1410 330 560 250  0.0000000 6.2831853}%
% ELLIPSE 2 0 3 0
% 4 1410 4490 1970 4740 1970 4740 1970 4740
% 
\special{pn 8}%
\special{ar 1410 4090 560 250  0.0000000 6.2831853}%
% LINE 2 0 3 0
% 2 1980 740 1980 4490
% 
\special{pn 8}%
\special{pa 1980 340}%
\special{pa 1980 4090}%
\special{fp}%
% LINE 2 0 3 0
% 2 850 750 850 4490
% 
\special{pn 8}%
\special{pa 850 350}%
\special{pa 850 4090}%
\special{fp}%
% ELLIPSE 2 0 3 0
% 4 5890 710 6450 960 6450 960 6450 960
% 
\special{pn 8}%
\special{ar 5890 310 560 250  0.0000000 6.2831853}%
% ELLIPSE 2 0 3 0
% 4 5890 4470 6450 4720 6450 4720 6450 4720
% 
\special{pn 8}%
\special{ar 5890 4070 560 250  0.0000000 6.2831853}%
% LINE 2 0 3 0
% 2 6460 720 6460 4470
% 
\special{pn 8}%
\special{pa 6460 320}%
\special{pa 6460 4070}%
\special{fp}%
% LINE 2 0 3 0
% 2 5330 730 5330 4470
% 
\special{pn 8}%
\special{pa 5330 330}%
\special{pa 5330 4070}%
\special{fp}%
% ELLIPSE 0 0 3 0
% 4 3520 2470 4080 2720 4080 2720 4080 2720
% 
\special{pn 20}%
\special{ar 3520 2070 560 250  0.0000000 6.2831853}%
% LINE 0 0 3 0
% 2 5880 2010 5880 3100
% 
\special{pn 20}%
\special{pa 5880 1610}%
\special{pa 5880 2700}%
\special{fp}%
% STR 2 0 3 0
% 3 1200 4900 1200 5000 2 0
% t=-0.5
\put(12.0000,-46.0000){\makebox(0,0)[lb]{t=-0.5}}%
% STR 2 0 3 0
% 3 3200 4900 3200 5000 2 0
% t=0
\put(32.0000,-46.0000){\makebox(0,0)[lb]{t=0}}%
% STR 2 0 3 0
% 3 5600 4900 5600 5000 2 0
% t=0.5
\put(56.0000,-46.0000){\makebox(0,0)[lb]{t=0.5}}%
% ELLIPSE 2 0 3 0
% 4 3510 720 4070 970 4070 970 4070 970
% 
\special{pn 8}%
\special{ar 3510 320 560 250  0.0000000 6.2831853}%
% LINE 0 0 3 0
% 2 3520 4350 3520 3130
% 
\special{pn 20}%
\special{pa 3520 3950}%
\special{pa 3520 2730}%
\special{fp}%
% LINE 0 0 3 0
% 2 3500 760 3500 1990
% 
\special{pn 20}%
\special{pa 3500 360}%
\special{pa 3500 1590}%
\special{fp}%
% ELLIPSE 2 0 3 0
% 4 3510 4480 4070 4730 4070 4730 4070 4730
% 
\special{pn 8}%
\special{ar 3510 4080 560 250  0.0000000 6.2831853}%
% LINE 2 0 3 0
% 2 4080 730 4080 4480
% 
\special{pn 8}%
\special{pa 4080 330}%
\special{pa 4080 4080}%
\special{fp}%
% LINE 2 0 3 0
% 2 2950 740 2950 4480
% 
\special{pn 8}%
\special{pa 2950 340}%
\special{pa 2950 4080}%
\special{fp}%
% ELLIPSE 2 0 3 0
% 4 1410 730 1970 980 1970 980 1970 980
% 
\special{pn 8}%
\special{ar 1410 330 560 250  0.0000000 6.2831853}%
% ELLIPSE 2 0 3 0
% 4 1410 4490 1970 4740 1970 4740 1970 4740
% 
\special{pn 8}%
\special{ar 1410 4090 560 250  0.0000000 6.2831853}%
% LINE 2 0 3 0
% 2 1980 740 1980 4490
% 
\special{pn 8}%
\special{pa 1980 340}%
\special{pa 1980 4090}%
\special{fp}%
% LINE 2 0 3 0
% 2 850 750 850 4490
% 
\special{pn 8}%
\special{pa 850 350}%
\special{pa 850 4090}%
\special{fp}%
% ELLIPSE 2 0 3 0
% 4 5890 710 6450 960 6450 960 6450 960
% 
\special{pn 8}%
\special{ar 5890 310 560 250  0.0000000 6.2831853}%
% ELLIPSE 2 0 3 0
% 4 5890 4470 6450 4720 6450 4720 6450 4720
% 
\special{pn 8}%
\special{ar 5890 4070 560 250  0.0000000 6.2831853}%
% LINE 2 0 3 0
% 2 6460 720 6460 4470
% 
\special{pn 8}%
\special{pa 6460 320}%
\special{pa 6460 4070}%
\special{fp}%
% LINE 2 0 3 0
% 2 5330 730 5330 4470
% 
\special{pn 8}%
\special{pa 5330 330}%
\special{pa 5330 4070}%
\special{fp}%
% LINE 2 2 3 0
% 2 3500 2000 5880 2000
% 
\special{pn 8}%
\special{pa 3500 1600}%
\special{pa 5880 1600}%
\special{dt 0.045}%
\special{pa 5880 1600}%
\special{pa 5879 1600}%
\special{dt 0.045}%
% LINE 2 2 3 0
% 2 3520 3150 5870 3130
% 
\special{pn 8}%
\special{pa 3520 2750}%
\special{pa 5870 2730}%
\special{dt 0.045}%
\special{pa 5870 2730}%
\special{pa 5869 2730}%
\special{dt 0.045}%
\end{picture}%

%% file: Figure2.tex
%WinTpicVersion2.15
\unitlength 0.1in
\begin{picture}(56.10,43.70)(8.50,-44.30)
% ELLIPSE 2 0 3 0
% 4 3510 720 4070 970 4070 970 4070 970
% 
\special{pn 8}%
\special{ar 3510 320 560 250  0.0000000 6.2831853}%
% LINE 2 0 3 0
% 2 5330 730 5330 4470
% 
\special{pn 8}%
\special{pa 5330 330}%
\special{pa 5330 4070}%
\special{fp}%
% LINE 2 0 3 0
% 2 6460 720 6460 4470
% 
\special{pn 8}%
\special{pa 6460 320}%
\special{pa 6460 4070}%
\special{fp}%
% ELLIPSE 2 0 3 0
% 4 5890 4470 6450 4720 6450 4720 6450 4720
% 
\special{pn 8}%
\special{ar 5890 4070 560 250  0.0000000 6.2831853}%
% ELLIPSE 2 0 3 0
% 4 5890 710 6450 960 6450 960 6450 960
% 
\special{pn 8}%
\special{ar 5890 310 560 250  0.0000000 6.2831853}%
% LINE 2 0 3 0
% 2 850 750 850 4490
% 
\special{pn 8}%
\special{pa 850 350}%
\special{pa 850 4090}%
\special{fp}%
% LINE 2 0 3 0
% 2 1980 740 1980 4490
% 
\special{pn 8}%
\special{pa 1980 340}%
\special{pa 1980 4090}%
\special{fp}%
% ELLIPSE 2 0 3 0
% 4 1410 4490 1970 4740 1970 4740 1970 4740
% 
\special{pn 8}%
\special{ar 1410 4090 560 250  0.0000000 6.2831853}%
% ELLIPSE 2 0 3 0
% 4 1410 730 1970 980 1970 980 1970 980
% 
\special{pn 8}%
\special{ar 1410 330 560 250  0.0000000 6.2831853}%
% LINE 2 0 3 0
% 2 2950 740 2950 4480
% 
\special{pn 8}%
\special{pa 2950 340}%
\special{pa 2950 4080}%
\special{fp}%
% LINE 2 0 3 0
% 2 4080 730 4080 4480
% 
\special{pn 8}%
\special{pa 4080 330}%
\special{pa 4080 4080}%
\special{fp}%
% ELLIPSE 2 0 3 0
% 4 3510 4480 4070 4730 4070 4730 4070 4730
% 
\special{pn 8}%
\special{ar 3510 4080 560 250  0.0000000 6.2831853}%
% LINE 0 0 3 0
% 2 3500 760 3500 1990
% 
\special{pn 20}%
\special{pa 3500 360}%
\special{pa 3500 1590}%
\special{fp}%
% LINE 0 0 3 0
% 2 3520 4350 3520 3130
% 
\special{pn 20}%
\special{pa 3520 3950}%
\special{pa 3520 2730}%
\special{fp}%
% ELLIPSE 2 0 3 0
% 4 3510 720 4070 970 4070 970 4070 970
% 
\special{pn 8}%
\special{ar 3510 320 560 250  0.0000000 6.2831853}%
% LINE 0 0 3 0
% 2 1430 2000 1430 3090
% 
\special{pn 20}%
\special{pa 1430 1600}%
\special{pa 1430 2690}%
\special{fp}%
% ELLIPSE 0 0 3 0
% 4 3520 2470 4080 2720 4080 2720 4080 2720
% 
\special{pn 20}%
\special{ar 3520 2070 560 250  0.0000000 6.2831853}%
% LINE 2 0 3 0
% 2 5330 730 5330 4470
% 
\special{pn 8}%
\special{pa 5330 330}%
\special{pa 5330 4070}%
\special{fp}%
% LINE 2 0 3 0
% 2 6460 720 6460 4470
% 
\special{pn 8}%
\special{pa 6460 320}%
\special{pa 6460 4070}%
\special{fp}%
% ELLIPSE 2 0 3 0
% 4 5890 4470 6450 4720 6450 4720 6450 4720
% 
\special{pn 8}%
\special{ar 5890 4070 560 250  0.0000000 6.2831853}%
% ELLIPSE 2 0 3 0
% 4 5890 710 6450 960 6450 960 6450 960
% 
\special{pn 8}%
\special{ar 5890 310 560 250  0.0000000 6.2831853}%
% LINE 2 0 3 0
% 2 850 750 850 4490
% 
\special{pn 8}%
\special{pa 850 350}%
\special{pa 850 4090}%
\special{fp}%
% LINE 2 0 3 0
% 2 1980 740 1980 4490
% 
\special{pn 8}%
\special{pa 1980 340}%
\special{pa 1980 4090}%
\special{fp}%
% ELLIPSE 2 0 3 0
% 4 1410 4490 1970 4740 1970 4740 1970 4740
% 
\special{pn 8}%
\special{ar 1410 4090 560 250  0.0000000 6.2831853}%
% ELLIPSE 2 0 3 0
% 4 1410 730 1970 980 1970 980 1970 980
% 
\special{pn 8}%
\special{ar 1410 330 560 250  0.0000000 6.2831853}%
% LINE 2 0 3 0
% 2 2950 740 2950 4480
% 
\special{pn 8}%
\special{pa 2950 340}%
\special{pa 2950 4080}%
\special{fp}%
% LINE 2 0 3 0
% 2 4080 730 4080 4480
% 
\special{pn 8}%
\special{pa 4080 330}%
\special{pa 4080 4080}%
\special{fp}%
% ELLIPSE 2 0 3 0
% 4 3510 4480 4070 4730 4070 4730 4070 4730
% 
\special{pn 8}%
\special{ar 3510 4080 560 250  0.0000000 6.2831853}%
% LINE 0 0 3 0
% 2 3520 4350 3520 3130
% 
\special{pn 20}%
\special{pa 3520 3950}%
\special{pa 3520 2730}%
\special{fp}%
% STR 2 0 3 0
% 3 1200 4900 1200 5000 2 0
% t=-0.5
\put(12.0000,-46.0000){\makebox(0,0)[lb]{t=-0.5}}%
% STR 2 0 3 0
% 3 3200 4900 3200 5000 2 0
% t=0
\put(32.0000,-46.0000){\makebox(0,0)[lb]{t=0}}%
% STR 2 0 3 0
% 3 5600 4900 5600 5000 2 0
% t=0.5
\put(56.0000,-46.0000){\makebox(0,0)[lb]{t=0.5}}%
% LINE 2 2 3 0
% 2 1420 2010 3490 2000
% 
\special{pn 8}%
\special{pa 1420 1610}%
\special{pa 3490 1600}%
\special{dt 0.045}%
\special{pa 3490 1600}%
\special{pa 3489 1600}%
\special{dt 0.045}%
% LINE 2 2 3 0
% 2 1430 3100 3520 3120
% 
\special{pn 8}%
\special{pa 1430 2700}%
\special{pa 3520 2720}%
\special{dt 0.045}%
\special{pa 3520 2720}%
\special{pa 3519 2720}%
\special{dt 0.045}%
\end{picture}%

%% file: torus.tex
%WinTpicVersion2.15
\unitlength 0.1in
\begin{picture}(31.98,20.48)(36.00,-28.05)
% ELLIPSE 2 0 3 0
% 4 5150 2169 6134 2799 6134 2799 6134 2799
% 
\special{pn 8}%
\special{ar 5150 1769 984 630  0.0000000 6.2831853}%
% ELLIPSE 2 0 3 0
% 4 5199 2181 6798 3205 6798 3205 6798 3205
% 
\special{pn 8}%
\special{ar 5199 1781 1599 1024  0.0000000 6.2831853}%
% SPLINE 2 2 3 0
% 33 3612 2194 3625 2144 3661 2120 3711 2095 3748 2083 3797 2058 3834 2058 3883 2058 3944 2058 3981 2058 4030 2058 4055 2071 4104 2120 4117 2132 4141 2169 4166 2181 4178 2206 4178 2243 4178 2317 4153 2366 4080 2415 4043 2427 4006 2440 3932 2440 3858 2440 3809 2440 3748 2403 3723 2366 3674 2317 3649 2280 3612 2255 3612 2243 3612 2243
% 
\special{pn 8}%
\special{pa 3612 1794}%
\special{pa 3617 1761}%
\special{pa 3633 1736}%
\special{pa 3662 1719}%
\special{pa 3690 1704}%
\special{pa 3720 1692}%
\special{pa 3750 1682}%
\special{pa 3778 1666}%
\special{pa 3808 1657}%
\special{pa 3840 1658}%
\special{pa 3872 1658}%
\special{pa 3904 1658}%
\special{pa 3936 1658}%
\special{pa 3968 1658}%
\special{pa 4001 1656}%
\special{pa 4032 1659}%
\special{pa 4060 1675}%
\special{pa 4082 1698}%
\special{pa 4105 1721}%
\special{pa 4125 1745}%
\special{pa 4143 1771}%
\special{pa 4171 1786}%
\special{pa 4179 1818}%
\special{pa 4178 1850}%
\special{pa 4180 1882}%
\special{pa 4179 1914}%
\special{pa 4168 1944}%
\special{pa 4149 1970}%
\special{pa 4125 1992}%
\special{pa 4097 2008}%
\special{pa 4067 2019}%
\special{pa 4037 2029}%
\special{pa 4006 2040}%
\special{pa 3975 2043}%
\special{pa 3943 2041}%
\special{pa 3910 2039}%
\special{pa 3878 2039}%
\special{pa 3846 2041}%
\special{pa 3814 2041}%
\special{pa 3784 2032}%
\special{pa 3757 2013}%
\special{pa 3737 1987}%
\special{pa 3719 1961}%
\special{pa 3695 1940}%
\special{pa 3673 1916}%
\special{pa 3657 1888}%
\special{pa 3629 1870}%
\special{pa 3612 1847}%
\special{pa 3612 1843}%
\special{sp -0.045}%
% SPLINE 2 2 3 0
% 30 6146 2144 6158 2058 6183 2021 6232 1972 6281 1948 6392 1923 6478 1923 6589 1948 6675 1997 6736 2058 6773 2095 6786 2144 6798 2194 6749 2292 6724 2304 6687 2341 6650 2366 6515 2403 6454 2403 6380 2403 6318 2390 6257 2354 6232 2304 6195 2267 6171 2243 6158 2218 6158 2206 6134 2181 6134 2169 6134 2169
% 
\special{pn 8}%
\special{pa 6146 1744}%
\special{pa 6147 1711}%
\special{pa 6151 1679}%
\special{pa 6162 1650}%
\special{pa 6181 1624}%
\special{pa 6202 1599}%
\special{pa 6225 1577}%
\special{pa 6252 1560}%
\special{pa 6282 1548}%
\special{pa 6312 1538}%
\special{pa 6343 1530}%
\special{pa 6375 1525}%
\special{pa 6407 1522}%
\special{pa 6439 1521}%
\special{pa 6471 1522}%
\special{pa 6503 1526}%
\special{pa 6534 1532}%
\special{pa 6565 1540}%
\special{pa 6596 1551}%
\special{pa 6625 1564}%
\special{pa 6653 1580}%
\special{pa 6678 1600}%
\special{pa 6700 1622}%
\special{pa 6723 1645}%
\special{pa 6747 1667}%
\special{pa 6770 1690}%
\special{pa 6781 1719}%
\special{pa 6788 1751}%
\special{pa 6796 1782}%
\special{pa 6798 1814}%
\special{pa 6788 1847}%
\special{pa 6770 1875}%
\special{pa 6744 1895}%
\special{pa 6716 1910}%
\special{pa 6694 1934}%
\special{pa 6670 1955}%
\special{pa 6642 1970}%
\special{pa 6612 1983}%
\special{pa 6581 1993}%
\special{pa 6550 1999}%
\special{pa 6518 2003}%
\special{pa 6486 2003}%
\special{pa 6454 2003}%
\special{pa 6422 2004}%
\special{pa 6390 2004}%
\special{pa 6358 2000}%
\special{pa 6327 1993}%
\special{pa 6296 1982}%
\special{pa 6268 1965}%
\special{pa 6249 1940}%
\special{pa 6236 1911}%
\special{pa 6216 1886}%
\special{pa 6192 1865}%
\special{pa 6170 1842}%
\special{pa 6158 1812}%
\special{pa 6139 1788}%
\special{pa 6134 1769}%
\special{sp -0.045}%
\end{picture}%

%% file: conce2.tex
%WinTpicVersion2.15
\unitlength 0.1in
\begin{picture}(38.25,31.47)(1.00,-32.96)
% LINE 2 0 3 0
% 2 2880 707 838 3696
% 
\special{pn 8}%
\special{pa 2880 307}%
\special{pa 838 3296}%
\special{fp}%
% BOX 2 5 2 0
% 2 2671 1974 2806 2158
% 
\special{pn 8}%
\special{sh 0}%
\special{pa 2671 1574}%
\special{pa 2806 1574}%
\special{pa 2806 1758}%
\special{pa 2671 1758}%
\special{pa 2671 1574}%
\special{ip}%
% ELLIPSE 0 0 3 0
% 4 3126 2478 3790 2982 3790 2982 4024 3204
% 
\special{pn 20}%
\special{ar 3126 2078 664 504  0.8166717 6.2831853}%
\special{ar 3126 2078 664 504  0.0000000 0.7853982}%
% LINE 2 0 3 0
% 2 2695 2761 2609 2970
% 
\special{pn 8}%
\special{pa 2695 2361}%
\special{pa 2609 2570}%
\special{fp}%
% BOX 2 5 2 0
% 2 2671 1974 2806 2158
% 
\special{pn 8}%
\special{sh 0}%
\special{pa 2671 1574}%
\special{pa 2806 1574}%
\special{pa 2806 1758}%
\special{pa 2671 1758}%
\special{pa 2671 1574}%
\special{ip}%
% BOX 2 5 2 0
% 2 2351 1334 2486 1519
% 
\special{pn 8}%
\special{sh 0}%
\special{pa 2351 934}%
\special{pa 2486 934}%
\special{pa 2486 1119}%
\special{pa 2351 1119}%
\special{pa 2351 934}%
\special{ip}%
% ELLIPSE 2 0 3 0
% 4 1933 2306 2769 3437 2843 3425 2843 3425
% 
\special{pn 8}%
\special{ar 1933 1906 836 1131  0.0000000 6.2831853}%
% BOX 2 5 2 0
% 2 2609 2773 2744 2958
% 
\special{pn 8}%
\special{sh 0}%
\special{pa 2609 2373}%
\special{pa 2744 2373}%
\special{pa 2744 2558}%
\special{pa 2609 2558}%
\special{pa 2609 2373}%
\special{ip}%
% LINE 2 5 3 0
% 2 2609 2798 2732 2872
% 
\special{pn 8}%
\special{pa 2609 2398}%
\special{pa 2732 2472}%
\special{ip}%
% LINE 0 0 3 0
% 2 2609 2785 2732 2884
% 
\special{pn 20}%
\special{pa 2609 2385}%
\special{pa 2732 2484}%
\special{fp}%
% BOX 2 5 2 0
% 2 1219 2995 1355 3179
% 
\special{pn 8}%
\special{sh 0}%
\special{pa 1219 2595}%
\special{pa 1355 2595}%
\special{pa 1355 2779}%
\special{pa 1219 2779}%
\special{pa 1219 2595}%
\special{ip}%
% LINE 2 5 3 0
% 2 1318 2982 1195 3179
% 
\special{pn 8}%
\special{pa 1318 2582}%
\special{pa 1195 2779}%
\special{ip}%
% LINE 2 0 3 0
% 2 1330 2970 1121 3290
% 
\special{pn 8}%
\special{pa 1330 2570}%
\special{pa 1121 2890}%
\special{fp}%
% STR 2 0 3 0
% 3 186 3733 186 3856 2 0
% The $zw$-plane
\put(1.8600,-34.5600){\makebox(0,0)[lb]{The $zw$-plane}}%
% STR 2 0 3 0
% 3 100 596 100 719 2 0
% $R^4$
\put(1.0000,-3.1900){\makebox(0,0)[lb]{$R^4$}}%
% STR 2 0 3 0
% 3 641 2035 641 2158 2 0
% $T^2$
\put(6.4100,-17.5800){\makebox(0,0)[lb]{$T^2$}}%
% STR 2 0 3 0
% 3 3925 2466 3925 2589 2 0
% $S^2$
\put(39.2500,-21.8900){\makebox(0,0)[lb]{$S^2$}}%
\end{picture}%

%% file: 1.1.tex
%WinTpicVersion2.15
\unitlength 0.1in
\begin{picture}(56.10,43.70)(8.50,-44.30)
% ELLIPSE 2 0 3 0
% 4 3510 720 4070 970 4070 970 4070 970
% 
\special{pn 8}%
\special{ar 3510 320 560 250  0.0000000 6.2831853}%
% LINE 0 0 3 0
% 2 3660 4350 3660 3130
% 
\special{pn 20}%
\special{pa 3660 3950}%
\special{pa 3660 2730}%
\special{fp}%
% LINE 0 0 3 0
% 2 3340 4370 3340 3140
% 
\special{pn 20}%
\special{pa 3340 3970}%
\special{pa 3340 2740}%
\special{fp}%
% SPLINE 0 0 3 0
% 41 3350 4400 3352 4387 3357 4375 3366 4364 3379 4353 3394 4343 3412 4335 3432 4329 3454 4324 3477 4321 3500 4320 3523 4321 3546 4324 3568 4329 3588 4335 3606 4343 3621 4353 3634 4364 3643 4375 3648 4387 3650 4400 3648 4413 3643 4425 3634 4436 3621 4447 3606 4457 3588 4465 3568 4471 3546 4476 3523 4479 3500 4480 3477 4479 3454 4476 3432 4471 3412 4465 3394 4457 3379 4447 3366 4436 3357 4425 3352 4413 3350 4400
% 
\special{pn 20}%
\special{pa 3350 4000}%
\special{pa 3360 3970}%
\special{pa 3384 3949}%
\special{pa 3413 3935}%
\special{pa 3444 3926}%
\special{pa 3475 3921}%
\special{pa 3507 3920}%
\special{pa 3539 3923}%
\special{pa 3570 3930}%
\special{pa 3601 3940}%
\special{pa 3627 3958}%
\special{pa 3647 3983}%
\special{pa 3648 4014}%
\special{pa 3630 4040}%
\special{pa 3603 4058}%
\special{pa 3573 4070}%
\special{pa 3542 4077}%
\special{pa 3510 4080}%
\special{pa 3478 4079}%
\special{pa 3447 4075}%
\special{pa 3416 4066}%
\special{pa 3387 4053}%
\special{pa 3362 4032}%
\special{pa 3350 4003}%
\special{pa 3350 4000}%
\special{sp}%
% LINE 0 0 3 0
% 2 3340 760 3340 1980
% 
\special{pn 20}%
\special{pa 3340 360}%
\special{pa 3340 1580}%
\special{fp}%
% LINE 0 0 3 0
% 2 3660 740 3660 1970
% 
\special{pn 20}%
\special{pa 3660 340}%
\special{pa 3660 1570}%
\special{fp}%
% ELLIPSE 0 0 3 0
% 4 3500 710 3650 790 4490 1280 4490 1280
% 
\special{pn 20}%
\special{ar 3500 310 150 80  0.0000000 6.2831853}%
% ELLIPSE 2 0 3 0
% 4 3510 4480 4070 4730 4070 4730 4070 4730
% 
\special{pn 8}%
\special{ar 3510 4080 560 250  0.0000000 6.2831853}%
% LINE 2 0 3 0
% 2 4080 730 4080 4480
% 
\special{pn 8}%
\special{pa 4080 330}%
\special{pa 4080 4080}%
\special{fp}%
% LINE 2 0 3 0
% 2 2950 740 2950 4480
% 
\special{pn 8}%
\special{pa 2950 340}%
\special{pa 2950 4080}%
\special{fp}%
% ELLIPSE 2 0 3 0
% 4 1410 730 1970 980 1970 980 1970 980
% 
\special{pn 8}%
\special{ar 1410 330 560 250  0.0000000 6.2831853}%
% ELLIPSE 2 0 3 0
% 4 1410 4490 1970 4740 1970 4740 1970 4740
% 
\special{pn 8}%
\special{ar 1410 4090 560 250  0.0000000 6.2831853}%
% LINE 2 0 3 0
% 2 1980 740 1980 4490
% 
\special{pn 8}%
\special{pa 1980 340}%
\special{pa 1980 4090}%
\special{fp}%
% LINE 2 0 3 0
% 2 850 750 850 4490
% 
\special{pn 8}%
\special{pa 850 350}%
\special{pa 850 4090}%
\special{fp}%
% ELLIPSE 2 0 3 0
% 4 5890 710 6450 960 6450 960 6450 960
% 
\special{pn 8}%
\special{ar 5890 310 560 250  0.0000000 6.2831853}%
% ELLIPSE 2 0 3 0
% 4 5890 4470 6450 4720 6450 4720 6450 4720
% 
\special{pn 8}%
\special{ar 5890 4070 560 250  0.0000000 6.2831853}%
% LINE 2 0 3 0
% 2 6460 720 6460 4470
% 
\special{pn 8}%
\special{pa 6460 320}%
\special{pa 6460 4070}%
\special{fp}%
% LINE 2 0 3 0
% 2 5330 730 5330 4470
% 
\special{pn 8}%
\special{pa 5330 330}%
\special{pa 5330 4070}%
\special{fp}%
% ELLIPSE 0 0 3 0
% 4 3500 1990 3650 2070 4490 2560 4490 2560
% 
\special{pn 20}%
\special{ar 3500 1590 150 80  0.0000000 6.2831853}%
% ELLIPSE 0 0 3 0
% 4 3500 3110 3650 3190 4490 3680 4490 3680
% 
\special{pn 20}%
\special{ar 3500 2710 150 80  0.0000000 6.2831853}%
% ELLIPSE 0 0 3 0
% 4 3520 2470 4080 2720 4080 2720 4080 2720
% 
\special{pn 20}%
\special{ar 3520 2070 560 250  0.0000000 6.2831853}%
% LINE 2 2 3 0
% 2 3660 2000 5730 2000
% 
\special{pn 8}%
\special{pa 3660 1600}%
\special{pa 5730 1600}%
\special{dt 0.045}%
\special{pa 5730 1600}%
\special{pa 5729 1600}%
\special{dt 0.045}%
% LINE 2 2 3 0
% 2 3680 3120 5750 3120
% 
\special{pn 8}%
\special{pa 3680 2720}%
\special{pa 5750 2720}%
\special{dt 0.045}%
\special{pa 5750 2720}%
\special{pa 5749 2720}%
\special{dt 0.045}%
% ELLIPSE 0 0 3 0
% 4 5850 2000 6000 2080 6840 2570 6840 2570
% 
\special{pn 20}%
\special{ar 5850 1600 150 80  0.0000000 6.2831853}%
% ELLIPSE 0 0 3 0
% 4 5850 3130 6000 3210 6840 3700 6840 3700
% 
\special{pn 20}%
\special{ar 5850 2730 150 80  0.0000000 6.2831853}%
% LINE 0 0 3 0
% 2 5690 2010 5690 3100
% 
\special{pn 20}%
\special{pa 5690 1610}%
\special{pa 5690 2700}%
\special{fp}%
% LINE 0 0 3 0
% 2 6020 2040 6020 3130
% 
\special{pn 20}%
\special{pa 6020 1640}%
\special{pa 6020 2730}%
\special{fp}%
% STR 2 0 3 0
% 3 1200 4900 1200 5000 2 0
% t=-0.5
\put(12.0000,-46.0000){\makebox(0,0)[lb]{t=-0.5}}%
% STR 2 0 3 0
% 3 3200 4900 3200 5000 2 0
% t=0
\put(32.0000,-46.0000){\makebox(0,0)[lb]{t=0}}%
% STR 2 0 3 0
% 3 5600 4900 5600 5000 2 0
% t=0.5
\put(56.0000,-46.0000){\makebox(0,0)[lb]{t=0.5}}%
% ELLIPSE 2 0 3 0
% 4 3510 720 4070 970 4070 970 4070 970
% 
\special{pn 8}%
\special{ar 3510 320 560 250  0.0000000 6.2831853}%
% LINE 2 0 3 1
% 2 2950 740 2950 4480
% 
\special{pn 8}%
\special{pa 2950 340}%
\special{pa 2950 4080}%
\special{fp}%
% LINE 2 0 3 2
% 2 4080 730 4080 4480
% 
\special{pn 8}%
\special{pa 4080 330}%
\special{pa 4080 4080}%
\special{fp}%
% ELLIPSE 2 0 3 3
% 4 3510 4480 4070 4730 4070 4730 4070 4730
% 
\special{pn 8}%
\special{ar 3510 4080 560 250  0.0000000 6.2831853}%
% ELLIPSE 0 0 3 4
% 4 3500 710 3650 790 4490 1280 4490 1280
% 
\special{pn 20}%
\special{ar 3500 310 150 80  0.0000000 6.2831853}%
% LINE 0 0 3 5
% 2 3660 740 3660 1970
% 
\special{pn 20}%
\special{pa 3660 340}%
\special{pa 3660 1570}%
\special{fp}%
% LINE 0 0 3 6
% 2 3340 760 3340 1980
% 
\special{pn 20}%
\special{pa 3340 360}%
\special{pa 3340 1580}%
\special{fp}%
% SPLINE 0 0 3 7
% 41 3350 4400 3352 4387 3357 4375 3366 4364 3379 4353 3394 4343 3412 4335 3432 4329 3454 4324 3477 4321 3500 4320 3523 4321 3546 4324 3568 4329 3588 4335 3606 4343 3621 4353 3634 4364 3643 4375 3648 4387 3650 4400 3648 4413 3643 4425 3634 4436 3621 4447 3606 4457 3588 4465 3568 4471 3546 4476 3523 4479 3500 4480 3477 4479 3454 4476 3432 4471 3412 4465 3394 4457 3379 4447 3366 4436 3357 4425 3352 4413 3350 4400
% 
\special{pn 20}%
\special{pa 3350 4000}%
\special{pa 3360 3970}%
\special{pa 3384 3949}%
\special{pa 3413 3935}%
\special{pa 3444 3926}%
\special{pa 3475 3921}%
\special{pa 3507 3920}%
\special{pa 3539 3923}%
\special{pa 3570 3930}%
\special{pa 3601 3940}%
\special{pa 3627 3958}%
\special{pa 3647 3983}%
\special{pa 3648 4014}%
\special{pa 3630 4040}%
\special{pa 3603 4058}%
\special{pa 3573 4070}%
\special{pa 3542 4077}%
\special{pa 3510 4080}%
\special{pa 3478 4079}%
\special{pa 3447 4075}%
\special{pa 3416 4066}%
\special{pa 3387 4053}%
\special{pa 3362 4032}%
\special{pa 3350 4003}%
\special{pa 3350 4000}%
\special{sp}%
% LINE 0 0 3 8
% 2 3340 4370 3340 3140
% 
\special{pn 20}%
\special{pa 3340 3970}%
\special{pa 3340 2740}%
\special{fp}%
% LINE 0 0 3 9
% 2 3660 4350 3660 3130
% 
\special{pn 20}%
\special{pa 3660 3950}%
\special{pa 3660 2730}%
\special{fp}%
% ELLIPSE 2 0 3 0
% 4 1410 730 1970 980 1970 980 1970 980
% 
\special{pn 8}%
\special{ar 1410 330 560 250  0.0000000 6.2831853}%
% ELLIPSE 2 0 3 0
% 4 1410 4490 1970 4740 1970 4740 1970 4740
% 
\special{pn 8}%
\special{ar 1410 4090 560 250  0.0000000 6.2831853}%
% LINE 2 0 3 0
% 2 1980 740 1980 4490
% 
\special{pn 8}%
\special{pa 1980 340}%
\special{pa 1980 4090}%
\special{fp}%
% LINE 2 0 3 0
% 2 850 750 850 4490
% 
\special{pn 8}%
\special{pa 850 350}%
\special{pa 850 4090}%
\special{fp}%
% ELLIPSE 2 0 3 0
% 4 5890 710 6450 960 6450 960 6450 960
% 
\special{pn 8}%
\special{ar 5890 310 560 250  0.0000000 6.2831853}%
% ELLIPSE 2 0 3 0
% 4 5890 4470 6450 4720 6450 4720 6450 4720
% 
\special{pn 8}%
\special{ar 5890 4070 560 250  0.0000000 6.2831853}%
% LINE 2 0 3 0
% 2 6460 720 6460 4470
% 
\special{pn 8}%
\special{pa 6460 320}%
\special{pa 6460 4070}%
\special{fp}%
% LINE 2 0 3 0
% 2 5330 730 5330 4470
% 
\special{pn 8}%
\special{pa 5330 330}%
\special{pa 5330 4070}%
\special{fp}%
\end{picture}%

%% file: 1.2.tex
%WinTpicVersion2.15
\unitlength 0.1in
\begin{picture}(56.10,43.70)(8.50,-44.30)
% ELLIPSE 2 0 3 0
% 4 3510 720 4070 970 4070 970 4070 970
% 
\special{pn 8}%
\special{ar 3510 320 560 250  0.0000000 6.2831853}%
% LINE 0 0 3 0
% 2 3660 4350 3660 3130
% 
\special{pn 20}%
\special{pa 3660 3950}%
\special{pa 3660 2730}%
\special{fp}%
% LINE 0 0 3 0
% 2 3340 4370 3340 3140
% 
\special{pn 20}%
\special{pa 3340 3970}%
\special{pa 3340 2740}%
\special{fp}%
% SPLINE 0 0 3 0
% 41 3350 4400 3352 4387 3357 4375 3366 4364 3379 4353 3394 4343 3412 4335 3432 4329 3454 4324 3477 4321 3500 4320 3523 4321 3546 4324 3568 4329 3588 4335 3606 4343 3621 4353 3634 4364 3643 4375 3648 4387 3650 4400 3648 4413 3643 4425 3634 4436 3621 4447 3606 4457 3588 4465 3568 4471 3546 4476 3523 4479 3500 4480 3477 4479 3454 4476 3432 4471 3412 4465 3394 4457 3379 4447 3366 4436 3357 4425 3352 4413 3350 4400
% 
\special{pn 20}%
\special{pa 3350 4000}%
\special{pa 3360 3970}%
\special{pa 3384 3949}%
\special{pa 3413 3935}%
\special{pa 3444 3926}%
\special{pa 3475 3921}%
\special{pa 3507 3920}%
\special{pa 3539 3923}%
\special{pa 3570 3930}%
\special{pa 3601 3940}%
\special{pa 3627 3958}%
\special{pa 3647 3983}%
\special{pa 3648 4014}%
\special{pa 3630 4040}%
\special{pa 3603 4058}%
\special{pa 3573 4070}%
\special{pa 3542 4077}%
\special{pa 3510 4080}%
\special{pa 3478 4079}%
\special{pa 3447 4075}%
\special{pa 3416 4066}%
\special{pa 3387 4053}%
\special{pa 3362 4032}%
\special{pa 3350 4003}%
\special{pa 3350 4000}%
\special{sp}%
% LINE 0 0 3 0
% 2 3340 760 3340 1980
% 
\special{pn 20}%
\special{pa 3340 360}%
\special{pa 3340 1580}%
\special{fp}%
% LINE 0 0 3 0
% 2 3660 740 3660 1970
% 
\special{pn 20}%
\special{pa 3660 340}%
\special{pa 3660 1570}%
\special{fp}%
% ELLIPSE 0 0 3 0
% 4 3500 710 3650 790 4490 1280 4490 1280
% 
\special{pn 20}%
\special{ar 3500 310 150 80  0.0000000 6.2831853}%
% ELLIPSE 2 0 3 0
% 4 3510 4480 4070 4730 4070 4730 4070 4730
% 
\special{pn 8}%
\special{ar 3510 4080 560 250  0.0000000 6.2831853}%
% LINE 2 0 3 0
% 2 4080 730 4080 4480
% 
\special{pn 8}%
\special{pa 4080 330}%
\special{pa 4080 4080}%
\special{fp}%
% LINE 2 0 3 0
% 2 2950 740 2950 4480
% 
\special{pn 8}%
\special{pa 2950 340}%
\special{pa 2950 4080}%
\special{fp}%
% ELLIPSE 2 0 3 0
% 4 1410 730 1970 980 1970 980 1970 980
% 
\special{pn 8}%
\special{ar 1410 330 560 250  0.0000000 6.2831853}%
% ELLIPSE 2 0 3 0
% 4 1410 4490 1970 4740 1970 4740 1970 4740
% 
\special{pn 8}%
\special{ar 1410 4090 560 250  0.0000000 6.2831853}%
% LINE 2 0 3 0
% 2 1980 740 1980 4490
% 
\special{pn 8}%
\special{pa 1980 340}%
\special{pa 1980 4090}%
\special{fp}%
% LINE 2 0 3 0
% 2 850 750 850 4490
% 
\special{pn 8}%
\special{pa 850 350}%
\special{pa 850 4090}%
\special{fp}%
% ELLIPSE 2 0 3 0
% 4 5890 710 6450 960 6450 960 6450 960
% 
\special{pn 8}%
\special{ar 5890 310 560 250  0.0000000 6.2831853}%
% ELLIPSE 2 0 3 0
% 4 5890 4470 6450 4720 6450 4720 6450 4720
% 
\special{pn 8}%
\special{ar 5890 4070 560 250  0.0000000 6.2831853}%
% LINE 2 0 3 0
% 2 6460 720 6460 4470
% 
\special{pn 8}%
\special{pa 6460 320}%
\special{pa 6460 4070}%
\special{fp}%
% LINE 2 0 3 0
% 2 5330 730 5330 4470
% 
\special{pn 8}%
\special{pa 5330 330}%
\special{pa 5330 4070}%
\special{fp}%
% ELLIPSE 0 0 3 0
% 4 3500 1990 3650 2070 4490 2560 4490 2560
% 
\special{pn 20}%
\special{ar 3500 1590 150 80  0.0000000 6.2831853}%
% ELLIPSE 0 0 3 0
% 4 3500 3110 3650 3190 4490 3680 4490 3680
% 
\special{pn 20}%
\special{ar 3500 2710 150 80  0.0000000 6.2831853}%
% ELLIPSE 0 0 3 0
% 4 3520 2470 4080 2720 4080 2720 4080 2720
% 
\special{pn 20}%
\special{ar 3520 2070 560 250  0.0000000 6.2831853}%
% STR 2 0 3 0
% 3 1200 4900 1200 5000 2 0
% t=-0.5
\put(12.0000,-46.0000){\makebox(0,0)[lb]{t=-0.5}}%
% STR 2 0 3 0
% 3 3200 4900 3200 5000 2 0
% t=0
\put(32.0000,-46.0000){\makebox(0,0)[lb]{t=0}}%
% STR 2 0 3 0
% 3 5600 4900 5600 5000 2 0
% t=0.5
\put(56.0000,-46.0000){\makebox(0,0)[lb]{t=0.5}}%
% ELLIPSE 2 0 3 0
% 4 3510 720 4070 970 4070 970 4070 970
% 
\special{pn 8}%
\special{ar 3510 320 560 250  0.0000000 6.2831853}%
% LINE 2 0 3 1
% 2 2950 740 2950 4480
% 
\special{pn 8}%
\special{pa 2950 340}%
\special{pa 2950 4080}%
\special{fp}%
% LINE 2 0 3 2
% 2 4080 730 4080 4480
% 
\special{pn 8}%
\special{pa 4080 330}%
\special{pa 4080 4080}%
\special{fp}%
% ELLIPSE 2 0 3 3
% 4 3510 4480 4070 4730 4070 4730 4070 4730
% 
\special{pn 8}%
\special{ar 3510 4080 560 250  0.0000000 6.2831853}%
% ELLIPSE 0 0 3 4
% 4 3500 710 3650 790 4490 1280 4490 1280
% 
\special{pn 20}%
\special{ar 3500 310 150 80  0.0000000 6.2831853}%
% LINE 0 0 3 5
% 2 3660 740 3660 1970
% 
\special{pn 20}%
\special{pa 3660 340}%
\special{pa 3660 1570}%
\special{fp}%
% LINE 0 0 3 6
% 2 3340 760 3340 1980
% 
\special{pn 20}%
\special{pa 3340 360}%
\special{pa 3340 1580}%
\special{fp}%
% SPLINE 0 0 3 7
% 41 3350 4400 3352 4387 3357 4375 3366 4364 3379 4353 3394 4343 3412 4335 3432 4329 3454 4324 3477 4321 3500 4320 3523 4321 3546 4324 3568 4329 3588 4335 3606 4343 3621 4353 3634 4364 3643 4375 3648 4387 3650 4400 3648 4413 3643 4425 3634 4436 3621 4447 3606 4457 3588 4465 3568 4471 3546 4476 3523 4479 3500 4480 3477 4479 3454 4476 3432 4471 3412 4465 3394 4457 3379 4447 3366 4436 3357 4425 3352 4413 3350 4400
% 
\special{pn 20}%
\special{pa 3350 4000}%
\special{pa 3360 3970}%
\special{pa 3384 3949}%
\special{pa 3413 3935}%
\special{pa 3444 3926}%
\special{pa 3475 3921}%
\special{pa 3507 3920}%
\special{pa 3539 3923}%
\special{pa 3570 3930}%
\special{pa 3601 3940}%
\special{pa 3627 3958}%
\special{pa 3647 3983}%
\special{pa 3648 4014}%
\special{pa 3630 4040}%
\special{pa 3603 4058}%
\special{pa 3573 4070}%
\special{pa 3542 4077}%
\special{pa 3510 4080}%
\special{pa 3478 4079}%
\special{pa 3447 4075}%
\special{pa 3416 4066}%
\special{pa 3387 4053}%
\special{pa 3362 4032}%
\special{pa 3350 4003}%
\special{pa 3350 4000}%
\special{sp}%
% LINE 0 0 3 8
% 2 3340 4370 3340 3140
% 
\special{pn 20}%
\special{pa 3340 3970}%
\special{pa 3340 2740}%
\special{fp}%
% LINE 0 0 3 9
% 2 3660 4350 3660 3130
% 
\special{pn 20}%
\special{pa 3660 3950}%
\special{pa 3660 2730}%
\special{fp}%
% ELLIPSE 2 0 3 0
% 4 1410 730 1970 980 1970 980 1970 980
% 
\special{pn 8}%
\special{ar 1410 330 560 250  0.0000000 6.2831853}%
% ELLIPSE 2 0 3 0
% 4 1410 4490 1970 4740 1970 4740 1970 4740
% 
\special{pn 8}%
\special{ar 1410 4090 560 250  0.0000000 6.2831853}%
% LINE 2 0 3 0
% 2 1980 740 1980 4490
% 
\special{pn 8}%
\special{pa 1980 340}%
\special{pa 1980 4090}%
\special{fp}%
% LINE 2 0 3 0
% 2 850 750 850 4490
% 
\special{pn 8}%
\special{pa 850 350}%
\special{pa 850 4090}%
\special{fp}%
% ELLIPSE 2 0 3 0
% 4 5890 710 6450 960 6450 960 6450 960
% 
\special{pn 8}%
\special{ar 5890 310 560 250  0.0000000 6.2831853}%
% ELLIPSE 2 0 3 0
% 4 5890 4470 6450 4720 6450 4720 6450 4720
% 
\special{pn 8}%
\special{ar 5890 4070 560 250  0.0000000 6.2831853}%
% LINE 2 0 3 0
% 2 6460 720 6460 4470
% 
\special{pn 8}%
\special{pa 6460 320}%
\special{pa 6460 4070}%
\special{fp}%
% LINE 2 0 3 0
% 2 5330 730 5330 4470
% 
\special{pn 8}%
\special{pa 5330 330}%
\special{pa 5330 4070}%
\special{fp}%
% LINE 2 2 3 0
% 2 3320 1990 1590 1990
% 
\special{pn 8}%
\special{pa 3320 1590}%
\special{pa 1590 1590}%
\special{dt 0.045}%
\special{pa 1590 1590}%
\special{pa 1591 1590}%
\special{dt 0.045}%
% LINE 2 2 3 0
% 2 3330 3140 1620 3140
% 
\special{pn 8}%
\special{pa 3330 2740}%
\special{pa 1620 2740}%
\special{dt 0.045}%
\special{pa 1620 2740}%
\special{pa 1621 2740}%
\special{dt 0.045}%
% ELLIPSE 0 0 3 0
% 4 1420 1990 1570 2070 2410 2560 2410 2560
% 
\special{pn 20}%
\special{ar 1420 1590 150 80  0.0000000 6.2831853}%
% ELLIPSE 0 0 3 1
% 4 1420 3120 1570 3200 2410 3690 2410 3690
% 
\special{pn 20}%
\special{ar 1420 2720 150 80  0.0000000 6.2831853}%
% LINE 0 0 3 2
% 2 1260 2000 1260 3090
% 
\special{pn 20}%
\special{pa 1260 1600}%
\special{pa 1260 2690}%
\special{fp}%
% LINE 0 0 3 3
% 2 1590 2030 1590 3120
% 
\special{pn 20}%
\special{pa 1590 1630}%
\special{pa 1590 2720}%
\special{fp}%
\end{picture}%

%% file: pass.tex
%WinTpicVersion2.15
\unitlength 0.1in
\begin{picture}(51.99,21.71)(4.01,-29.42)
% ELLIPSE 2 0 3 0
% 4 1454 2286 2507 3342 1770 2729 1770 2729
% 
\special{pn 8}%
\special{ar 1454 1886 1053 1056  0.0000000 6.2831853}%
% VECTOR 2 0 3 0
% 2 636 2951 1969 1385
% 
\special{pn 8}%
\special{pa 636 2551}%
\special{pa 1969 985}%
\special{fp}%
\special{sh 1}%
\special{pa 1969 985}%
\special{pa 1911 1023}%
\special{pa 1934 1026}%
\special{pa 1941 1049}%
\special{pa 1969 985}%
\special{fp}%
% VECTOR 2 0 3 0
% 2 2263 1621 931 3194
% 
\special{pn 8}%
\special{pa 2263 1221}%
\special{pa 931 2794}%
\special{fp}%
\special{sh 1}%
\special{pa 931 2794}%
\special{pa 989 2756}%
\special{pa 965 2753}%
\special{pa 959 2730}%
\special{pa 931 2794}%
\special{fp}%
% LINE 2 0 3 0
% 2 1675 2227 1461 2072
% 
\special{pn 8}%
\special{pa 1675 1827}%
\special{pa 1461 1672}%
\special{fp}%
% LINE 2 0 3 0
% 2 1446 2508 1240 2330
% 
\special{pn 8}%
\special{pa 1446 2108}%
\special{pa 1240 1930}%
\special{fp}%
% VECTOR 2 0 3 0
% 2 3221 2441 2816 2442
% 
\special{pn 8}%
\special{pa 3221 2041}%
\special{pa 2816 2042}%
\special{fp}%
\special{sh 1}%
\special{pa 2816 2042}%
\special{pa 2883 2062}%
\special{pa 2869 2042}%
\special{pa 2883 2022}%
\special{pa 2816 2042}%
\special{fp}%
% VECTOR 2 0 3 0
% 2 2831 2271 3258 2271
% 
\special{pn 8}%
\special{pa 2831 1871}%
\special{pa 3258 1871}%
\special{fp}%
\special{sh 1}%
\special{pa 3258 1871}%
\special{pa 3191 1851}%
\special{pa 3205 1871}%
\special{pa 3191 1891}%
\special{pa 3258 1871}%
\special{fp}%
% STR 2 0 3 0
% 3 2757 2744 2757 2817 34 89
% pass-move
\put(27.5700,-24.1700){\makebox(0,0)[lb]{pass-move}}%
% ELLIPSE 2 0 3 0
% 4 4547 2227 5600 3283 4103 2544 4103 2544
% 
\special{pn 8}%
\special{ar 4547 1827 1053 1056  0.0000000 6.2831853}%
% LINE 2 0 3 0
% 2 4326 2218 4503 2013
% 
\special{pn 8}%
\special{pa 4326 1818}%
\special{pa 4503 1613}%
\special{fp}%
% VECTOR 2 0 3 0
% 2 4539 2545 4068 3165
% 
\special{pn 8}%
\special{pa 4539 2145}%
\special{pa 4068 2765}%
\special{fp}%
\special{sh 1}%
\special{pa 4068 2765}%
\special{pa 4124 2724}%
\special{pa 4100 2723}%
\special{pa 4092 2700}%
\special{pa 4068 2765}%
\special{fp}%
% LINE 2 0 3 0
% 2 4259 2330 3759 2914
% 
\special{pn 8}%
\special{pa 4259 1930}%
\special{pa 3759 2514}%
\special{fp}%
% VECTOR 2 0 3 0
% 2 4599 1909 5063 1318
% 
\special{pn 8}%
\special{pa 4599 1509}%
\special{pa 5063 918}%
\special{fp}%
\special{sh 1}%
\special{pa 5063 918}%
\special{pa 5006 958}%
\special{pa 5030 960}%
\special{pa 5038 983}%
\special{pa 5063 918}%
\special{fp}%
% LINE 2 0 3 0
% 2 4606 2470 4812 2219
% 
\special{pn 8}%
\special{pa 4606 2070}%
\special{pa 4812 1819}%
\special{fp}%
% LINE 2 0 3 0
% 2 4879 2161 5350 1561
% 
\special{pn 8}%
\special{pa 4879 1761}%
\special{pa 5350 1161}%
\special{fp}%
% VECTOR 2 0 3 0
% 2 1778 2315 2382 2781
% 
\special{pn 8}%
\special{pa 1778 1915}%
\special{pa 2382 2381}%
\special{fp}%
\special{sh 1}%
\special{pa 2382 2381}%
\special{pa 2341 2324}%
\special{pa 2340 2348}%
\special{pa 2317 2356}%
\special{pa 2382 2381}%
\special{fp}%
% LINE 2 0 3 0
% 2 1380 1983 783 1495
% 
\special{pn 8}%
\special{pa 1380 1583}%
\special{pa 783 1095}%
\special{fp}%
% LINE 2 0 3 0
% 2 2160 3061 1519 2574
% 
\special{pn 8}%
\special{pa 2160 2661}%
\special{pa 1519 2174}%
\special{fp}%
% VECTOR 2 0 3 0
% 2 1174 2264 547 1746
% 
\special{pn 8}%
\special{pa 1174 1864}%
\special{pa 547 1346}%
\special{fp}%
\special{sh 1}%
\special{pa 547 1346}%
\special{pa 586 1404}%
\special{pa 588 1380}%
\special{pa 611 1373}%
\special{pa 547 1346}%
\special{fp}%
% VECTOR 2 0 3 0
% 2 3928 1399 5483 2729
% 
\special{pn 8}%
\special{pa 3928 999}%
\special{pa 5483 2329}%
\special{fp}%
\special{sh 1}%
\special{pa 5483 2329}%
\special{pa 5445 2270}%
\special{pa 5442 2294}%
\special{pa 5419 2301}%
\special{pa 5483 2329}%
\special{fp}%
% VECTOR 2 0 3 0
% 2 5217 3032 3634 1717
% 
\special{pn 8}%
\special{pa 5217 2632}%
\special{pa 3634 1317}%
\special{fp}%
\special{sh 1}%
\special{pa 3634 1317}%
\special{pa 3673 1375}%
\special{pa 3675 1351}%
\special{pa 3698 1344}%
\special{pa 3634 1317}%
\special{fp}%
\end{picture}%

%% file: XO+.tex
%WinTpicVersion2.15
\unitlength 0.1in
\begin{picture}(23.22,15.54)(20.56,-22.52)
% STR 2 0 3 0
% 3 2686 1533 2686 1713 4 0
% $S^1\x$
\put(26.8600,-13.1300){\makebox(0,0)[rt]{$S^1\x$}}%
% CIRCLE 0 0 3 0
% 4 3526 1875 4066 2434 4066 2434 4156 2559
% 
\special{pn 20}%
\special{ar 3526 1475 777 777  0.8264709 6.2831853}%
\special{ar 3526 1475 777 777  0.0000000 0.8026849}%
% LINE 2 0 3 0
% 2 3652 1966 4120 2326
% 
\special{pn 8}%
\special{pa 3652 1566}%
\special{pa 4120 1926}%
\special{fp}%
% VECTOR 2 0 3 0
% 2 3418 1767 2950 1390
% 
\special{pn 8}%
\special{pa 3418 1367}%
\special{pa 2950 990}%
\special{fp}%
\special{sh 1}%
\special{pa 2950 990}%
\special{pa 2989 1047}%
\special{pa 2992 1023}%
\special{pa 3014 1016}%
\special{pa 2950 990}%
\special{fp}%
% VECTOR 2 0 3 0
% 2 3040 2470 4030 1318
% 
\special{pn 8}%
\special{pa 3040 2070}%
\special{pa 4030 918}%
\special{fp}%
\special{sh 1}%
\special{pa 4030 918}%
\special{pa 3971 956}%
\special{pa 3995 958}%
\special{pa 4002 982}%
\special{pa 4030 918}%
\special{fp}%
% STR 2 0 3 0
% 3 4378 2325 4378 2505 2 0
% ,
\put(43.7800,-21.0500){\makebox(0,0)[lb]{,}}%
\end{picture}%

%% file: XO-.tex
%WinTpicVersion2.15
\unitlength 0.1in
\begin{picture}(23.94,15.54)(20.56,-21.86)
% STR 2 0 3 0
% 3 2686 1545 2686 1725 4 0
% $S^1\x$
\put(26.8600,-13.2500){\makebox(0,0)[rt]{$S^1\x$}}%
% STR 2 0 3 0
% 3 4450 2229 4450 2409 2 0
% ,
\put(44.5000,-20.0900){\makebox(0,0)[lb]{,}}%
% CIRCLE 0 0 3 0
% 4 3514 1809 4077 1274 4077 1274 4203 1184
% 
\special{pn 20}%
\special{ar 3514 1409 777 777  5.5464549 6.2831853}%
\special{ar 3514 1409 777 777  0.0000000 5.5232825}%
% LINE 2 0 3 0
% 2 3075 2379 3438 1914
% 
\special{pn 8}%
\special{pa 3075 1979}%
\special{pa 3438 1514}%
\special{fp}%
% VECTOR 2 0 3 0
% 2 4119 2322 2974 1322
% 
\special{pn 8}%
\special{pa 4119 1922}%
\special{pa 2974 922}%
\special{fp}%
\special{sh 1}%
\special{pa 2974 922}%
\special{pa 3011 981}%
\special{pa 3014 957}%
\special{pa 3037 951}%
\special{pa 2974 922}%
\special{fp}%
% VECTOR 2 0 3 0
% 2 3603 1683 3972 1209
% 
\special{pn 8}%
\special{pa 3603 1283}%
\special{pa 3972 809}%
\special{fp}%
\special{sh 1}%
\special{pa 3972 809}%
\special{pa 3915 849}%
\special{pa 3939 851}%
\special{pa 3947 874}%
\special{pa 3972 809}%
\special{fp}%
\end{picture}%

%% file: higharX2.bbl
\begin{thebibliography}{ABCD}


\bibitem{Boy} 
W. Boy: \"Uber die Curvatura integra und die Topologie geschlossener Flächen, 
{\it Math. Ann.} 57 (1903) 151-184.
 

\bibitem{CappellShaneson}
S. E. Cappell and J. L. Shaneson: 
There Exist Inequivalent Knots With the Same Complement
{\it Ann. of Math. Second Series}103 (1976) 349-353. 


\bibitem{CassonGordon75} A. J. Casson and C. McA. Gordon: 
 Cobordism of classical knots
{\it Progress in Mathematics } 62, 1975. 

\bibitem{CassonGordon78} A. J. Casson and C. McA. Gordon: 
 On slice knots in dimension three
{\it Proceedings of Symposia in Pure Mathematics} 32, 1978. 


\bibitem{CochranOrr}  T. D. Cochran and K. E. Orr: 
Not all links are concordant to boundary links 
{\it Ann. of Math.}, 138, 519--554, 1993. 


\bibitem{CochranOrrTeichner}
T. D. Cochran, K. E. Orr and P. Teichner: 
Knot Concordance, Whitney Towers and L2-Signatures
{\it 
Annals of Mathematics} 
Second Series, Vol. 157, (2003), 433-519. 


\bibitem{Farber1983} M. Farber: 
Classification of simple knots. 
{\it (Russian) Uspekhi Mat. Nauk 38},   
no. 5(233), 59--106, 1983. 

\bibitem{Farber1984I}
M. Š. Farber:  
An Algebraic Classification of Some Even-Dimensional Spherical Knots. I
{\it Transactions of the American Mathematical Society} 
Vol. 281 (1984) 507-527. 


\bibitem{Farber1984II}
M. Š. Farber: 
An Algebraic Classification of Some Even-Dimensional Spherical Knots. II
{\it Transactions of the American Mathematical Society} 
Vol. 281 (1984) 529-570. 




\bibitem{Gordon}
C. McA Gordon:  Knots in the 4-sphere
{\it Commentarii Mathematici Helvetici} 51 (1976) 585-596. 



\bibitem{Giller} 
 C.   Giller:  Towards a classical knot theory for surfaces in $\Bbb R^4$
{\it Illinois.  J. Math.}  
 26 (1982) 591-631. 



\bibitem{GilmerLivingston92} P. M. Gilmer and C. Livingston: 
The Casson-Gordon invariant and link concordance
{\it Topology} 31, 475-492, 1992. 


\bibitem{Haefliger} 
A. Haefliger: 
Knotted (4k - 1)-Spheres in 6k-Space  
{\it Annals of Mathematics} 75 (1962) 452-466.



\bibitem{Kauffman}   
L. H. Kauffman:  
Products of knots, 
{\it 
Bulletin of the American mathematical society
Vol 80, No 6, November 1104-1107  1974.}  


\bibitem{KauffmanNeumann}   
L. H. Kauffman and W. D. Neumann:  
Products of knots, branched fibrations and sums of singularities, 
{\it Topology, vol 16, No 4, 369-393, 1977}



\bibitem{Kauffman1987}
L. H. Kauffman:  On Knots 
 {\it  Ann. of Math. Stud. }, vol. 115, 1987



\bibitem{Kauffman1994}
L. H. Kauffman: 
Knots and Physics,   Second Edition 
{\it World Scientific Publishing} 1994. 





\bibitem{KauffmanOgasa} L.H.Kauffman and  E. Ogasa: 
 Local moves on knots and products of knots, 
{\it  Knots in Poland III-Part III Banach Center Publications Volume} 
103 (2014), 159-209 
{\it Institute of Mathematics Warszawa} 2014  　arXiv:1210.4667[math.GT] 


\bibitem{KauffmanOgasa2} L.H.Kauffman and  E. Ogasa: 
 Local moves on knots and products of knots II,  
 arXiv:1406.5573[math.GT] 


\bibitem{KauffmanOgasa3} L.H.Kauffman and  E. Ogasa: 
 Brieskorn submanifolds, Local moves on knots, and knot products, 
 arXiv:1504.01229 




\bibitem{Kervaire}   M. Kervaire:  
Les noeudes de dimensions sup\'ereures,  
{\it  Bull.Soc.Math.France}  93, 225-271, 1965.  




\bibitem{Levineslice}  J. Levine:  Knot cobordism in codimension two, 
{\it Comment. Math. Helv.},  44,  229-244, 1969. 


\bibitem{Levinesimp} J. Levine: 
An algebraic classification of some knots of codimension two. 
{\it Comment. Math. Helv. 45},  185--198, 1970. 


\bibitem{LevineAlex} J. Levine: 
Polynomial invariants of knots of codimension two. 
{\it Ann. of Math.  84,  537--554}, 1966. 



\bibitem{LevineOrr} 
J. Levine and K. Orr: 
A survey of applications of surgery to knot and link theory. 
{\it Surveys on surgery theory: surveys presented 
in honor of C.T.C. Wall Vol. 1, 345--364, Ann. of Math. Stud.,
145, Princeton Univ. Press, Princeton, NJ, 2000.}


\bibitem{Levine94} J. Levine: 
Link invariants via the eta-invariant
{\it  Comment. Math. Helveticii} 69,  82-119, 1994. 

\bibitem{Milnorexotic}
J. W. Milnor: 
On manifolds homeomorphic to the 7-sphere {\it Annals of Mathematics} 64 (1956) 399–405. 


\bibitem{MilnorStasheff} J. W. Milnor and J. D. Stasheff:
 Characteristic classes. 
 {\it Annals of Mathematics Studies, No. 76. 
 Princeton University Press} 1974. 



\bibitem{Ogasa02} E. Ogasa:   
The intersection of spheres in a sphere and 
a new geometric meaning of the Arf invariants, 
{\it Journal of knot theory and its ramifications, 
 11(2002) 1211-1231}, 
Univ. of Tokyo preprint series UTMS 95-7, 
math.GT/0003089 in http://xxx.lanl.gov. 


\bibitem{Ogasa98n} E. Ogasa:   
Intersectional pairs of $n$-knots, local moves of $n$-knots  
   and invariants of $n$-knots, 
{\it Mathematical Research Letters, 1998, vol5, 577-582},  
Univ. of Tokyo preprint UTMS 95-50. 


\bibitem{Ogasa98SL} E. Ogasa:  
The intersection of spheres in a sphere and  
  a new application of the Sato-Levine invariant,  
{\it Proceedings of the American Mathematical Society 
126, 1998, PP.3109-3116.} UTMS95-54.



\bibitem{Ogasa98F} E. Ogasa:  
Some properties of ordinary sense slice 1-links:
    Some answers to the problem (26) of Fox, 
{\it Proceedings of the American Mathematical Society 
126, 1998, PP.2175-2182. }  UTMS96-11.


\bibitem{Ogasa01P} E. Ogasa:  
The projections of $n$-knots which are not 
    the projection of any unknotted knot, 
{\it Journal of knot theory and its ramifications, 
 10 (2001), no. 1, 121--132. }
 UTMS 97-34, math.GT/0003088.  


\bibitem{Ogasa99S} E. Ogasa:  
Singularities of projections of $n$-dimensional knots, 
{\it Mathematical Proceedings of Cambridge Philosophical Society 
126, 1999, 511-519. }  
UTMS96-39



\bibitem{Ogasa99D} E. Ogasa:  
Link cobordism and the intersection of slice discs, 
{\it Bulletin of the London Mathematical Society, 
31, 1999, 1-8.}



\bibitem{Ogasa04} E. Ogasa:  
Ribbon-moves of 2-links preserve the $\mu$-invariant of 2-links, 
{\it Journal of knot theory and its ramifications, 13 (2004), no. 5, 669--687.}
UTMS 97-35, math.GT/0004008.  




\bibitem{Ogasa09} E. Ogasa:  
 Local move identities for the Alexander polynomials of 
high-dimensional knots and inertia groups,  
{\it 
Journal of Knot Theory and Its Ramifications (JKTR) Volume: 18 No: 4 Year: 2009 pp. 531-545}    
 UTMS 97-63. math.GT/0512168. 


\bibitem{Ogasa01NR} E. Ogasa:  
Nonribbon 2-links all of whose components are trivial knots 
     and some of whose bund-sums are nonribbon knots, 
 {\it Journal of knot theory and its ramifications,  
10 (2001), no. 6, 913--922.}



\bibitem{OgasaBS} E. Ogasa:  
$n$-dimensional links, their components, and their band-sums, 
math.GT/0011163. UTMS00-65.





\bibitem{Ogasa07} E. Ogasa:  
 Ribbon-moves of 2-knots: The Farber-Levine pairing and 
the Atiyah-Pathodi-Singer-Casson-Gordon-Ruberman $\tilde\eta$ invariant of 2-knots, 
{\it Journal of Knot Theory and Its Ramifications, Vol. 16, No. 5 (2007) 523-543}
math.GT/0004007, UTMS 00-22,  math.GT/0407164. 



\bibitem{Ogasa06} E. Ogasa:  
  Supersymmetry, homology with twisted coefficients and 
$n$-dimensional knots, 
{\it International Journal of Modern Physics A,}  
Vol. 21, Nos. 19-20 (2006) 4185-4196
  hep-th/0311136. 




\bibitem{Ogasanu} E. Ogasa:  
 A new invariant associated with decompositions of manifolds,  
math.GT/0512320,  hep-th/0401217



\bibitem{OgasaT3} E. Ogasa: 
A new obstruction for ribbon-moves of 2-knots: 2-knots fibred by the punctured 3-torus and 2-knots bounded by the Poicar\'e sphere,  
arXiv:1003.2473math.GT 



\bibitem{OgasaBoy} E. Ogasa: 
Make your Boy surface, 
arXiv:1303.6448 math.GT



\bibitem{OgasaZ} E. Ogasa: 
 Local-move-identities for the $Z[t,t^{-1}]$-Alexander polynomials of 2-links, the alinking number, and high dimensional analogues  arXiv:1602.07775 



\bibitem{OgasaOCSL} E. Ogasa: 
 A new pair of non-cobordant surface-links which the Orr invariant, the Cochran sequence, the Sato-Levine invariant, and the alinking number cannot find  arXiv:1605.06874 




\bibitem{Ogasanikai} E. Ogasa: 
Ribbon-move-unknotting-number-two 2-knots, pass-move-unknotting-number-two 1-knots, and high dimensional analogue,  arXiv:1612.03325 




\bibitem{Oi} E. Ogasa: Ijigen e no tobira  (In Japanese) 
{\it Nippon Hyoron Sha Co., Ltd.} 2009.  

\bibitem{Oy} E. Ogasa: Yojigen ijou no k\^ukan ga mieru (In Japanese)
{Beret Shuppan}  2006. 




\bibitem{Orr} 
K. E. Orr: 
New link invariants and applications  
{\it Comment. Math. Helv.} 
62 (1987) 542-560.




\bibitem{Rolfsen} 
D.   Rolfsen:  Knots and links 
{\it Publish or Perish, Inc.}  1976.

\bibitem{Ruberman}  
D. Ruberman:  
Doubly slice knots and the Casson-Gordon invariants  
{\it Trans. Amer. Math. Soc.}  279 (1983) 569-588.  

\bibitem{Sato}  
N. Sato: `Cobordisms of semi-boundary links'  
{\it Topology and its application}18 (1984) 225-234. 


\bibitem{Suciu}
A. I. Suciu: Inequivalent frame-spun knots with the same complement
{\it Commentarii Mathematici Helvetici} 67 (1992) 47-63. 



\bibitem{Zeeman} E. Zeeman: Twisting spun knots 
{\it Trans. AMS} 115 (1965) 471-495. 	 





\end{thebibliography}
